\documentclass[a4paper, 12pt]{article}
\usepackage{acronym}
\usepackage{algorithm}
\usepackage{algpseudocode}
\usepackage{amsmath}
\usepackage{amsthm}
\usepackage{amssymb}
\usepackage{authblk}
\usepackage[backend=biber,
defernumbers=false,
maxbibnames=5,
style=numeric-comp,
isbn=false,
bibencoding=utf8,
safeinputenc, 
url=false,
doi=true,
giveninits=true,sorting=none]{biblatex}
\usepackage{blindtext}
\usepackage{caption}
\usepackage{enumerate}
\usepackage[margin=2.5cm]{geometry}
\usepackage{graphicx} 
\usepackage{halloweenmath}
\usepackage{hyperref}
\usepackage[capitalise]{cleveref} 
\usepackage{mathtools}
\usepackage{multirow}
\usepackage{subcaption}
\usepackage{tikz}
\usetikzlibrary{decorations.pathreplacing}
\usepackage{xcolor}

\hypersetup{
breaklinks=true,
bookmarksopen=true,
pdftitle={Multigrid Ghost FEM},    
pdfauthor={Dilip.H, Coco.A},     
colorlinks=true, 
linkcolor = blue,
citecolor=blue,  
filecolor=blue,      
urlcolor=blue           
}
\newtheorem*{prop}{Proposition}

\acrodef{amg}[AMG]{algebraic multigrid}
\acrodef{xfem}[XFEM]{extended finite element method}
\acrodef{fe}[FE]{finite element}
\acrodef{fem}[FEM]{finite element method}
\acrodef{fdm}[FDM]{finite difference method}
\acrodef{gmg}[GMG]{geometric multigrid}
\acrodef{hpc}[HPC]{high performance computing}
\acrodef{jit}[JIT]{just-in-time}
\acrodef{mg}[MG]{multigrid}
\acrodef{pde}[PDE]{partial differential equation}
\acrodef{tgcs}[TGCS]{Two-Grid Correction Scheme}

\renewcommand{\vec}[1]{\mathbf{#1}}
\makeatletter
\def\thanks#1{\protected@xdef\@thanks{\@thanks
        \protect\footnotetext{#1}}}
\makeatother

\title{Multigrid methods for the ghost finite element approximation of elliptic problems}

\author{Hridya Dilip$^{1,\ast}$ and Armando Coco$^{1}$ }

\thanks{$^{1}$ Department of Mathematics and Computer Science, University of Catania, 95125, Catania, Italy.\
$^\ast$ Corresponding author.\
E-mails: {\tt hridya.dilip@unict.it},
{\tt armando.coco@unict.it}
}
\date{}

\begin{document}
\maketitle

\begin{abstract}
We present multigrid methods for solving elliptic partial differential equations on arbitrary domains using the nodal ghost finite element method, an unfitted boundary approach where the domain is implicitly defined by a level-set function. This method achieves second-order accuracy and offers substantial computational advantages over both direct solvers and finite-difference-based multigrid methods. A key strength of the ghost finite element framework is its variational formulation, which naturally enables consistent transfer operators and avoids residual splitting across grid levels.
We provide a detailed construction of the multigrid components in both one and two spatial dimensions, including smoothers, transfer operators, and coarse grid operators. The choice of the stabilization parameter plays a crucial role in ensuring well-posedness and optimal convergence of the multigrid method. We derive explicit algebraic expressions for this parameter based on the geometry of cut cells. In the two-dimensional setting, we further improve efficiency by performing additional smoothing exclusively on cut cells, reducing computational cost without compromising convergence. Numerical results validate the proposed method across a range of geometries and confirm its robustness and scalability.
\end{abstract}

\section{Introduction} \label{sec:intro}
\par \Ac{mg} methods are powerful tools for enhancing the computational efficiency of numerical methods. They work by solving the problem at multiple levels of resolution, smoothing out high-frequency errors on coarser grids and then refining the solution on finer grids \cite{Briggs2000,trottenberg2000multigrid}.  Originally developed for elliptic equations, \ac{mg} methods have been widely adopted to solve a variety of problems due to their optimal computational efficiency, which scales linearly with the number of computational nodes for sparse matrices, outperforming many other numerical solvers. A comprehensive survey of robust \ac{mg} methods for solving second-order elliptic \acp{pde} can be found in \cite{Chan2000}.
\par
While \ac{mg} methods are highly efficient for rectilinear domains, their optimal performance is not always guaranteed for more complex geometries. With the growing interest in solving PDEs on arbitrary domains, \ac{mg} techniques for irregularly shaped domains have attracted increasing attention. In such cases, boundary effects can significantly reduce \ac{mg} efficiency~\cite{Brandt1984, Coco2023}.

Modeling complex-shaped domains is commonly done using fitted-boundary methods, with the \acp{fem}~\cite{shaidurov2013multigrid, dePrenter2019} being the most widely adopted approach. The popularity of \acp{fem} relies on the use of isoparametric elements, which allow for flexible and precise domain discretization. However, \acp{fem} present significant difficulties when dealing with highly intricate geometries: (i) generating meshes that conform to irregular boundary curvatures can be cumbersome, and (ii) implementing parallel solvers requires considerable effort to partition the mesh efficiently for balanced computational workload distribution.

To overcome these challenges, unfitted-boundary methods have gained popularity in recent years.
In the context of \acp{fdm} and level-set approaches~\cite{gibou2018review},
unfitted-boundary methods are particularly advantageous in \ac{hpc} environments, as they naturally facilitate the design of parallel solvers and eliminate the need for mesh generation by embedding the domain within a fixed grid.
The typical approach of \acp{fdm} for curved domains is the Ghost-Point method~\cite{Peskin:IBM, LeVequeLi:IIM, Fedkiw:GFM, Gibou:Ghost, Gibou:fourth_order}.

\par \ac{mg} methods have been successfully applied to Ghost-Point \acp{fdm} on arbitrary domains for a range of problems. These include elliptic problems \cite{COCO2013464, Coco2023}, elliptic interface problems \cite{COCO2018299}, the Navier-Stokes equation \cite{COCO2020109623}, multi-scale modeling of sorption kinetics of a surfactant past an oscillating bubble \cite{ASTUTO2023111880}, and higher-order problems \cite{CocoHighOrder}. The stability analysis for the finite difference approximation of elliptic problems has also been explored in \cite{CocoStissi2023}. However, \acp{fdm} often have limitations, particularly when dealing with complex geometries and irregular domains,
since a stability and convergence analysis is cumbersome~\cite{COCO2013464}.

\par In contrast, unfitted \acp{fem} utilize simple background meshes, such as Cartesian meshes, eliminating the need for a mesh that conforms to the geometry. This makes them a promising alternative to finite differences.  However, they are prone to ill-conditioning when the intersection between a background cell and the domain is small \cite{dePrenter2017}. Several techniques have been proposed to mitigate this issue, including cell-aggregation methods \cite{Kummer2016,Badia2018}, the addition of stabilization terms on the boundary \cite{Burman2014}, and the use of a snapping-back-to-grid mechanism \cite{Astuto2024Ghost}. The ghost-\ac{fem} is a second-order accurate method that leverages a snapping-back-to-grid algorithm to address ill-conditioning problems \cite{Astuto2024Ghost}.

\par The existing literature on \ac{mg} methods for unfitted \ac{fe} approximations is sparse. \ac{mg} methods for the numerical approximation of elliptic interface problems have been explored using \ac{xfem} in \cite{Kothari2021} and using Cut-\ac{fem} in \cite{Ludescher2020}. Additionally, a \ac{gmg} preconditioner for higher-order immersed \acp{fem} has been introduced in \cite{dePrenter2019}. This technique has a linear computational cost and is robust to cut elements. To the best of our knowledge, no \ac{mg} methods have been developed for the ghost \ac{fe} framework. In this approach, Dirichlet boundary conditions are enforced weakly using Nitsche's method. A crucial aspect of this formulation is the choice of the stabilization parameter, which ensures the well-posedness of the problem. However, this parameter also significantly affects the conditioning of the \ac{mg} solver \cite{Saberi2023}. Thus, selecting an optimal stabilization parameter is essential for maintaining the well-posedness of the discrete formulation and achieving optimal \ac{mg} convergence.

\par The key contributions of this work are: (i) introducing \ac{mg} methods for the ghost-\ac{fem}, (ii) providing algebraic expressions for the stabilization parameter based on the configuration of the cut cell, and (iii) presenting an alternative approach to improve \ac{mg} efficiency with lower computational cost. 

\par The outline of this paper is as follows. In \cref{sec:ghost_fem}, we briefly explain the nodal ghost-\ac{fem} for approximating elliptic \acp{pde}. We describe the 1D two-grid correction scheme in \cref{sec:2_grid_scheme_1D}, followed by the V-cycle and W-cycle methods in \cref{sec:1D_v_w_cycle}. The \ac{mg} methods for one-dimensional problems are discussed in \cref{sec:1D_mg}, where we first present the smoother in \cref{sec:1D_smoother}, followed by the transfer operators in \cref{sec:1D_transfer_op} and coarse grid operators in \cref{sec:1D_coarse_op}. The selection of the stabilization parameter is detailed in \cref{sec:optimal_par_1D}, and numerical simulations for one-dimensional problems are provided in \cref{sec:1D_tests}. Next, we extend our discussion to the two-dimensional setting in \cref{sec:2D_mg}. We describe the transfer operators for the two-grid correction scheme in \cref{sec:2D_transfer_op}, and discuss the choice of the stabilization parameter in \cref{sec:2D_optimal_lambda}.  Numerical tests for two-dimensional problems are presented in \cref{sec:2D_tests}. Finally, we summarize our findings in \cref{sec:conclusions}.

\section{The ghost finite element method} \label{sec:ghost_fem}
\par Let us consider the Poisson equation with mixed boundary conditions as the model problem: 
\begin{equation}\label{eq:poisson_model_prob}
    -\boldsymbol{\Delta }u = f \text{~ in ~}\Omega, \quad u = g_D \text{~ on ~} \Gamma_D, \text{~ and } \quad  \boldsymbol{n}\cdot \boldsymbol{\nabla}u = g_N \text{~ on ~} \Gamma_N.
\end{equation}
Here, $\Omega \subset \mathbb{R}^d, \text{~with~} d \in \{1,2\}$ represents the domain, $\Gamma_D$ and $\Gamma_N$ denote the Dirichlet and Neumann boundaries, respectively, with $\Gamma_D \cup \Gamma_N =\partial \Omega$ and 
$\Gamma_D \cap \Gamma_N =\emptyset$. The vector $\boldsymbol{n}$ is the outward unit normal to $\partial \Omega$. The function $u:\mathbb{R}^d \rightarrow \mathbb{R}$ represents the solution, $f \in H^{-1}(\Omega)$ is the source term, $g_D \in H^{1/2}(\Gamma_D)$ is the Dirichlet boundary condition, and $g_N \in H^{-1/2}(\Gamma_N)$ is the Neumann boundary condition. 

\par The variational formulation corresponding to \cref{eq:poisson_model_prob} is to seek $u \in H^1(\Omega)$ with a trace $g_D \in H^{1/2}(\Gamma_D)$ such that 
\begin{equation}\label{eq:var_formulation}
    a(u,v) = l(v) \quad \forall ~v \in H^1(\Omega),
\end{equation}
where the bilinear form $a(\cdot,\cdot)$ and the linear functional $l(\cdot)$ are defined as
\begin{align}\label{eq:bil_forms}
    a(u,v) &\doteq \int_{\Omega} \boldsymbol{\nabla}u \cdot \boldsymbol{\nabla} v \mathrm{~d}\boldsymbol{x} - \int_{\Gamma_D} \left(\boldsymbol{n} \cdot \boldsymbol{\nabla} u \right) v \mathrm{~d}\boldsymbol{S}, \\ 
    l(v) &\doteq \int_{\Omega} f ~v \mathrm{~d}\boldsymbol{x} + \int_{\Gamma_N} g_N ~ v \mathrm{~d}\boldsymbol{S}.
\end{align}

\par Before introducing the discrete approximation to \cref{eq:var_formulation}, we first define some key terms and notations. Let $\Omega_{art}$ be a rectilinear domain that can be easily meshed using a Cartesian grid, with $\Omega \subset \Omega_{art}$. We denote by $\mathcal{T}_h$ a conforming, shape regular and quasi-uniform partition of $\Omega_{art}$, where each cell is represented by $K \in \mathcal{T}_h$. The size of a cell is given by $h_K \doteq meas(K)$, and the overall mesh size of the partition $\mathcal{T}_h$ is defined as $h \doteq \max_{K \in \mathcal{T}_h} h_K$. The cells of the background mesh are categorized as follows: A cell $K$ is called an internal cell if $K \subset \Omega$. A cell is an external cell if $K \cap \Omega = \emptyset$. Otherwise $K$ is called a cut cell. The set of internal, external and cut cells are denoted by $\mathcal{T}_{h,in},\mathcal{T}_{h,ext} $ and $\mathcal{T}_{h,cut}$, respectively. The internal and cut cells are collectively referred to as active cells, and the corresponding active mesh is defined as $\mathcal{T}_{h,act} \doteq \mathcal{T}_{h,in}\cup \mathcal{T}_{h,cut}$. The active domain is defined as $\Omega_{h,act} = \bigcup_{K\in \mathcal{T}_{h,act}} K$.

\par Let $\mathcal{V}_{h,act} \subset H^1(\Omega_{h,act})$ be a nodal Lagrangian \ac{fe} space defined on the active domain, that is, 
\begin{equation}\label{eq:Vhact}
\mathcal{V}_{h,act} \doteq \{v_h \in H^1(\Omega_{h,act}) : v_h|_K \in \mathcal{Q}_1 \text{ for any } K \in \mathcal{T}_{h,act} \},
\end{equation}
where $\mathcal{Q}_1$ denotes the space of polynomials of degree atmost one in each spatial coordinate. The discrete approximation of \cref{eq:var_formulation} is to find $u_h \in \mathcal{V}_{h,act}$ such that 
\begin{align}\label{eq:disc_formulation}
        a_h(u_h,v_h) = l_h(v_h) \quad \forall ~v_h \in \mathcal{V}_{h,act},
\end{align}
where 
\begin{align}
    &a_h(u,v) \doteq \sum_{K \in \mathcal{T}_{h,act}} \left( \int_{K \cap \Omega} \boldsymbol{\nabla}u \cdot \boldsymbol{\nabla} v \mathrm{~d}\boldsymbol{x} + \int_{K \cap\Gamma_D} \lambda_K u ~v -\left(\boldsymbol{n} \cdot \boldsymbol{\nabla} u \right) v  -\left(\boldsymbol{n} \cdot \boldsymbol{\nabla} v \right) u \mathrm{~d}\boldsymbol{S} \right), \label{eq:disc_bil_form}\\ 
    &l_h(v) \doteq \sum_{K \in \mathcal{T}_{h,act}} \left( \int_{K \cap \Omega} f ~v \mathrm{~d}\boldsymbol{x} + \int_{K \cap\Gamma_D} \lambda_K g_D ~v   -\left(\boldsymbol{n} \cdot \boldsymbol{\nabla} v \right) g_D \mathrm{~d}\boldsymbol{S}  + \int_{K \cap \Gamma_N} g_N ~ v \mathrm{~d}\boldsymbol{S} \right). \label{eq:disc_rhs}
\end{align}
Here, the cell-wise stabilization parameter is defined as $\lambda_K = \gamma_K ~ h_K^{-\beta}$, where $\gamma_K \in \mathbb{R}^+$ and $\beta \ge 1$. The parameter must be chosen sufficiently large to guarantee the coercivity of the bilinear form $a_h$ on $\mathcal{V}_{h,act}$ \cite{Nitsche1971}.

\par The linear system corresponding to the discrete problem \cref{eq:disc_formulation} becomes severely ill-conditioned in the presence of small cut cells. To mitigate this issue, the ghost-\ac{fem} employs a snapping-back-to-grid mechanism \cite{Astuto2024Ghost}. Before describing this mechanism, we introduce the following preliminaries. We assume that the domain $\Omega$ is defined using a level-set. Let $\psi : \mathbb{R}^d \rightarrow \mathbb{R}$ be a function such that $\Omega = \{\boldsymbol{x} \in \mathbb{R}^d : \psi(\boldsymbol{x}) < 0\}$. The boundary of $\Omega$ is given by $\partial \Omega = \{\boldsymbol{x}: \psi(\boldsymbol{x}) = 0\}$. We choose a snapping threshold as $h^\alpha$, for some $\alpha \in \mathbb{R}^{+}$. If a point $\boldsymbol{x} \in \Omega$, satisfies $| \psi(\boldsymbol{x}) | < h^\alpha$, we enforce $\psi(\boldsymbol{x}) = 0$, effectively snapping it to the closest boundary. Let $\mathcal{T}_{h,snap} \subseteq \mathcal{T}_{h,act}$ denote the set of active cells after snapping-back-to-grid and the corresponding active domain is $\Omega_{h,snap} \doteq \bigcup_{K\in \mathcal{T}_{h,snap}} K$. The \ac{fe} space on $\Omega_{h,snap}$ is defined as 
\[ \mathcal{V}_{h} = \{v_h \in H^1(\Omega_{h,snap}) : v_h|_K \in \mathcal{Q}_1 \text{ for any } K \in \mathcal{T}_{h,snap} \}.\] 
The linear system corresponding to the discrete approximation \cref{eq:disc_formulation} on $\mathcal{V}_{h}$ is well-conditioned.

\par For elliptic problems, the ghost-\ac{fem} is second-order accurate \cite{Astuto2024Ghost}. However, the numerical experiments in \cite{Astuto2024Ghost} use a direct solver, which becomes computationally expensive for large systems. Employing a \ac{mg} solver enhances computational efficiency.

\section{Multigrid method}\label{sec:mg_method}
We briefly summarize the main idea behind the \ac{mg} method and refer the reader to more detailed sources for an in-depth discussion (see, for example,~\cite{trottenberg2000multigrid}).  

We denote by $\Omega_{h}$ the set of grid nodes that are associated with $\mathcal{T}_{h}$, namely the vertices of the cells of $\mathcal{T}_{h}$.
Let us denote by \( A_h \boldsymbol{u}_h = \boldsymbol{F}_h \) the linear system arising from the ghost-\ac{fem} of Section~\ref{sec:ghost_fem}.

\subsection{Two-grid correction scheme in 1D} \label{sec:2_grid_scheme_1D}

The process begins by applying a relaxation scheme to the linear system \( A_h \boldsymbol{u}_h = \boldsymbol{F}_h \) for a small number of iterations, typically \(\nu_1\), known as \textit{pre-smoothing} steps. This yields an approximate solution, denoted as \(\hat{\vec{u}}_h\). A commonly used relaxation method is a Richardson-type iteration:  

\begin{equation}\label{RichRel}
\textbf{u}_h^{(m+1)} = \textbf{u}_h^{(m)} + P^{-1} (\textbf{F}_h-A_h \textbf{u}_h^{(m)})
\end{equation}  
where \( P \) is a suitable preconditioner matrix and $m$ is the iteration parameter.
Next, the residual \(\vec{r}_h = \vec{F}_h - A_h \hat{\vec{u}}_h\) is computed on the fine grid \(\Omega_h\) and then transferred to a coarser grid \(\Omega_{H}\), with $H>h$, using a restriction operator $ \vec{r}_{H} = \mathcal{I}^h_{H} \vec{r}_h$.

On this coarse grid, the residual equation \( A_{H} \vec{e}_{H} = \vec{r}_{H} \) is solved exactly, and the computed error correction \(\vec{e}_{H}\) is interpolated back to the fine grid using the interpolation operator $\vec{e}_{h} = \mathcal{I}^{H}_{h} \vec{e}_{H}$.

The initial approximation is then updated as \( \hat{\vec{u}}_h \coloneqq \hat{\vec{u}}_h + \vec{e}_h \), followed by \(\nu_2\) additional \textit{post-smoothing} iterations of \eqref{RichRel} on the fine grid. This completes one \textit{two-grid iteration (cycle)}, which is repeated until the residual on the fine grid meets a specified tolerance.  

This approach is referred to as the {\em \ac{tgcs}} because it involves only two grids: \(\Omega_h\) and \(\Omega_{H}\). 

\par For simplicity, we use uniform Cartesian grids $\Omega_h$ and $\Omega_{2h}$ as the fine and coarse grids (then $H=2h$). See Fig.~\ref{fig:ghostFem1D} for an illustration of the geometrical quantities associated with $\Omega_h$. 
The key steps of one cycle of the \ac{tgcs} are presented in \cref{alg:2Grid_scheme}. The data transfer operations between the two grids, as well as the matrix on the coarse grid, are explained in detail in \cref{sec:1D_transfer_op,sec:1D_coarse_op}.

\begin{algorithm}
\caption{One cycle of the \ac{tgcs} (from $\boldsymbol{u}^{(m)}_{h}$ to $\boldsymbol{u}^{(m+1)}_{h}$)}
\label{alg:2Grid_scheme}
\begin{algorithmic}[1]
\State Relax $\nu_1$ times $A_h \boldsymbol{u}_h =  \boldsymbol{F}_h$  on $\Omega_h$ with an initial guess $\boldsymbol{u}^{(m)}_{h}$. Let $\hat{\vec{u}}_h$ be the solution after $\nu_1$ iterations.
\State Compute the residual $\boldsymbol{r}_h \doteq \boldsymbol{F}_h- A_h \hat{\vec{u}}_h$ on $\Omega_h$ and its restriction $\boldsymbol{r}_{2h} \doteq \mathcal{I}^h_{2h} \boldsymbol{r}_{h} $ on $\Omega_{2h}$.
\State Solve $A_{2h} \boldsymbol{e}_{2h} =  \boldsymbol{r}_{2h}$ on $\Omega_{2h}$ exactly. \label{alg:state:reseq}
\State Compute the interpolation of the error $\boldsymbol{e}_{h} \doteq \mathcal{I}^{2h}_h \boldsymbol{e}_{2h} $ on $\Omega_h$.
\State The fine grid correction is calculated as $\boldsymbol{u}_h = \hat{\vec{u}}_h + \boldsymbol{e}_{h}$.
\State Relax $\nu_2$ times $A_h \boldsymbol{u}_h =  \boldsymbol{F}_h$ on $\Omega_h$ with an initial guess $\boldsymbol{u}_{h}$, obtaining $\boldsymbol{u}^{(m+1)}_h$. 
\end{algorithmic}
\end{algorithm}

\par After each iteration of the \ac{tgcs}, the residual is computed by $\boldsymbol{r}^{(m)}_h \doteq \boldsymbol{F}_h- A_h \boldsymbol{u}^{(m)}_h$. The performance of the scheme is measured by its convergence factor $\rho^{(m)}$, computed as 
\[
\rho^{(m)} \doteq 
\frac{\| \boldsymbol{r}^{(m)}_h \|}{\| \boldsymbol{r}^{(m-1)}_h \|}
\]
where $\|\cdot\|_\infty$ represents the $l_\infty$-norm.

\subsection{V-cycle and W-cycle}\label{sec:1D_v_w_cycle}
When the coarse-grid equation \( A_{2h} \vec{e}_{2h} = \vec{r}_{2h} \) is not solved exactly but rather approximated by applying \ac{tgcs} on $\Omega_{2h}$ and $\Omega_{4h}$ a few times, say \(\gamma^*\), and so on recursively, the method extends into a \textit{\ac{mg} scheme}. The recursion stops when the grid becomes coarse enough that using an exact solver is computationally inexpensive.

The number of recursive cycles \(\gamma^*\) determines the type of \ac{mg} cycle used: {\em V-cycle} if \(\gamma^* = 1\), meaning a single recursive application, and {\em W-cycle} if \(\gamma^* = 2\), where the process is repeated twice per level.  
Higher values of \(\gamma^*\) are rarely used, as they increase computational cost without providing significant efficiency gains~\cite{trottenberg2000multigrid}.  
Finally, we note that it is sufficient to compute the convergence factor for the \ac{tgcs}, as the W-cycle exhibits similar convergence properties under natural assumptions, as detailed in Section 3.2 of~\cite{trottenberg2000multigrid}.

\section{One-dimensional case} \label{sec:1D_mg}
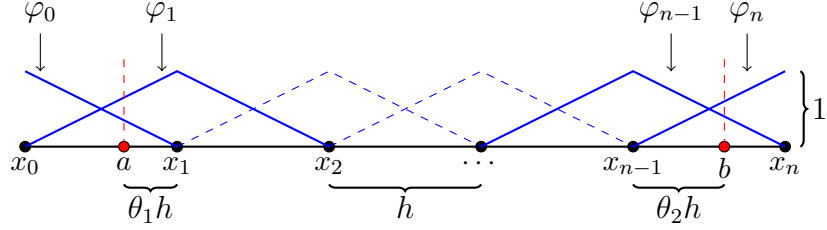
\begin{figure}
    \centering
  \begin{tikzpicture}
    \draw[thick] (0,0) -- (10,0); 

    \foreach \i in {0,2,4,6,8,10}
    {
        \draw[fill=black] (\i,0) circle (2pt);
    }

    \draw[fill=red] (1.3,0) circle (2pt);
    \draw[fill=red] (9.2,0) circle (2pt);

    \draw[blue, thick]  (0,1) -- (2,0);
    \foreach \i in {2,8}
    {
        \draw[blue, thick] (\i-2,0) -- (\i,1) -- (\i+2,0);
    }
    
    \foreach \i in {4,6}
    {
        \draw[blue, dashed] (\i-2,0) -- (\i,1) -- (\i+2,0);
    }
    \draw[blue, thick]  (8,0) -- (10,1);

    \draw[red, dashed] (1.3,0) -- (1.3,1.2);
    \draw[red, dashed] (9.2,0) -- (9.2,1.2);

    \node[below] at (0,0) {\small $x_0$};
    \node[below] at (1.3,0) {\small $a$};
    \node[below] at (2,0) {\small $x_1$};
    \node[below] at (4,0) {\small $x_2$};
    \node[below] at (6,0) {\small $\cdots$};
    \node[below] at (8,0) {\small $x_{n-1}$};
    \node[below] at (9.2,0) {\small $b$};
    \node[below] at (10,0) {\small $x_n$};

    \node[above] at (0.2,1.5) {$\varphi_0$};
    \draw[->] (0.2,1.5) -- (0.2,1);
    \node[above] at (1.8,1.5) {$\varphi_1$};
    \draw[->] (1.8,1.5) -- (1.8,1);
    \node[above] at (8.5,1.5) {$\varphi_{n-1}$};
    \draw[->] (8.5,1.5) -- (8.5,1);
    \node[above] at (9.5,1.5) {$\varphi_n$};
    \draw[->] (9.5,1.5) -- (9.5,1);

    \draw[decorate,decoration={brace,mirror},thick] (1.3,-0.5) -- (2,-0.5) node[midway,below] {$\theta_1 h$};
    \draw[decorate,decoration={brace,mirror},thick] (4,-0.5) -- (6,-0.5) node[midway,below] {$h$};
    \draw[decorate,decoration={brace,mirror},thick] (8,-0.5) -- (9.2,-0.5) node[midway,below] {$\theta_2 h$};
    \draw[decorate,decoration={brace,mirror},thick] (10.2,0) -- (10.2,1) node[midway,right] {1};
\end{tikzpicture}

    \caption{Setup of a one-dimensional grid: boundary points and basis functions.}
    \label{fig:ghostFem1D}
\end{figure}

Let us formulate the linear system for the one-dimensional case. Assume that $\Omega_{art}=[0,1]$ is partitioned using a uniform grid with $n$ equal-length intervals. The spatial step is then given by $h=1/n$.
We take \(\Omega = (a,b)\), with \(a = (1 - \theta_1)h\) and \(b = 1 - (1 - \theta_2)h\), where \(\theta_1, \theta_2 \in (0,1]\).
A basis for the vector space $\mathcal{V}_{h,act}$ defined in \eqref{eq:Vhact} is given by the $n+1$ functions $\varphi_i = 1- |x-x_i|/h$ for $i=0, \ldots, n$ (see Fig.~\ref{fig:ghostFem1D}).
These functions are referred to as {\em shape functions}.

We prescribe a Dirichlet boundary condition at \( x = a \) and a Neumann boundary condition at \( x = b \), so that \( \Gamma_D = \{ a \} \) and \( \Gamma_N = \{ b \} \).

The linear system corresponding to~\cref{eq:disc_formulation} is $A_h \vec{u}_h = \vec{F}_h$ with $A_h =A_h^I+A_h^{B} + (A_h^{B})^T + A_h^{\lambda}$, where:
\begin{gather}\label{eq:ls1dA}
A^I_h = 
\frac{1}{h}
\begin{pmatrix}
 \theta_1 & -\theta_1 & 0 & 0 & 0 & \cdots & 0 \\
 -\theta_1 & 1+\theta_1 & -1 & 0 & 0 & \cdots & 0 \\
 0 & -1 & 2 & -1 & 0 & \cdots & 0 \\
 \vdots & & \ddots & \ddots & \ddots &  & \vdots \\
 0 & \cdots & 0 & -1 & 2 & -1 & 0 \\
0 & \cdots & 0 & 0 & -1 & 1+\theta_2 & -\theta_2 \\
0 & \cdots & 0 & 0  & 0 & -\theta_2 & \theta_2 \\
\end{pmatrix}, \\ 
A^B_h = 
\frac{1}{h}
\begin{pmatrix}
 -\theta_1 &\theta_1 & 0 & \cdots & 0 \\
 \theta_1-1 & 1-\theta_1 & 0 & \cdots & 0 \\
 0 & 0 & 0 & \cdots & 0 \\
 \vdots & \vdots & \vdots & \ddots  & \vdots \\
0 & 0 & 0 & \cdots & 0 \\
\end{pmatrix}, \quad
A^\lambda_h = 
\lambda
\begin{pmatrix}
 \theta_1^2 &\theta_1(1-\theta_1) & 0 & \cdots & 0 \\
 \theta_1(1-\theta_1)  & (1-\theta_1)^2 & 0 & \cdots & 0 \\
 0 & 0 & 0 & \cdots & 0 \\
 \vdots & \vdots & \vdots & \ddots  & \vdots \\
0 & 0 & 0 & \cdots & 0 \\
\end{pmatrix}.
\end{gather}
Similarly $\boldsymbol{F}_h = \boldsymbol{F}_h^I + \boldsymbol{F}_h^B +\boldsymbol{F}_h^\lambda+\boldsymbol{F}_h^N$, where 
\begin{align}\label{eq:ls1dF}
    \boldsymbol{F}^I_h =  
    \begin{pmatrix}
        \tilde{f}_0 \\
        \tilde{f}_1 \\
        \tilde{f}_2 \\
        \vdots \\  
        \tilde{f}_n \\
    \end{pmatrix}, \quad
    \boldsymbol{F}^B_h = \frac{1}{h}
    \begin{pmatrix}
        -g_a \\
        g_a \\
        0 \\ 
        \vdots \\
        0
    \end{pmatrix}, \quad
    \boldsymbol{F}^\lambda_h = 
    \lambda 
    \begin{pmatrix}
        \theta_1 ~g_a\\
        (1-\theta_1)g_a \\ 
        0 \\
        \vdots \\ 
        0   
    \end{pmatrix}, \quad
     \boldsymbol{F}^N_h =  
    \begin{pmatrix}
        0\\
        \vdots \\
        0\\
        (1-\theta_2)g_b \\ 
       \theta_2 ~ g_b  
    \end{pmatrix},
\end{align}
with $\tilde{f}_i = \int_\Omega f \varphi_i dx$ for $i=0, \ldots, n$.

\subsection{Smoother}\label{sec:1D_smoother}
For a \ac{mg} strategy to be effective, the relaxation scheme \eqref{RichRel} must efficiently dampen high-frequency components of the error \( \vec{e}_h \), a property known as the \textit{smoothing property}. This ensures that even if low-frequency components are not immediately reduced, the overall convergence of the method remains robust. A relaxation method that possesses this characteristic is referred to as a \textit{smoother}.  
It is well established that the Jacobi method, where \( P = \text{diag}(A_h) \) represents the diagonal of \( A_h \), does not function effectively as a smoother. Instead, relaxation methods that satisfy the smoothing property include the \textit{weighted-Jacobi} scheme (with a proper weight) and the \textit{Gauss-Seidel} method.  

The weighted Jacobi smoother is defined by choosing the preconditioner as $P = \omega \, \text{diag}(A_h) + (1 - \omega) {I}$,
where \( {I} \) is the identity matrix. The optimal smoother in the one-dimensional case is obtained when \( \omega = 2/3 \).
The Gauss-Seidel method, another widely used smoother, is obtained by choosing \( P = L_h \), where \( L_h \) represents the lower triangular part of \( A_h \), including its diagonal. This approach effectively reduces high-frequency error components, improving the performance of the \ac{mg} method.  
In this paper, we use the Gauss-Seidel method as a smoother.

In some ghost-point approaches used in \acp{fdm}, standard Jacobi or Gauss-Seidel schemes may fail to converge due to the presence of ghost equations, which do not guarantee a positive definite linear system. This issue arises in the methods presented in~\cite{COCO2013464, COCO2018299, ASTUTO2023111880}, whose main idea is briefly summarized below.  

\subsubsection*{Ghost-point \ac{fdm} and limitations}

We briefly recall the linear system that arises from the finite-difference discretization of the problem. The internal equations are discretized using the standard second-order finite-difference scheme:
\[
\frac{-u_{i-1} + 2u_i - u_{i+1}}{h^2} = f_i, \quad \text{for } i=1,\ldots,n-1.
\]

To achieve second-order accuracy, the Dirichlet boundary condition is discretized via linear interpolation:
\[
\theta_1 u_{0} + (1-\theta_1) u_{1} = g_a,
\]
while the Neumann boundary condition is discretized using quadratic interpolation:
\[
\frac{(2\theta_2 -1) u_{n-2} - 4\theta_2 u_{n-1} + (2\theta_2 +1) u_{n}}{2h} = g_b.
\]

The resulting linear system is given by:
\begin{equation}\label{eq:linsysfdm}
\frac{1}{h^2}
\begin{pmatrix}
 \theta_1 h^2 & (1-\theta_1)h^2 & 0 & 0 & 0 & \cdots & 0 \\
 -1 & 2 & -1 & 0 & 0 & \cdots & 0 \\
 0 & -1 & 2 & -1 & 0 & \cdots & 0 \\
 \vdots & \vdots & \ddots & \ddots & \ddots &  & \vdots \\
0 & 0 &  \cdots & -1 & 2 & -1 & 0 \\
0 & 0 & \cdots & 0 & -1 & 2 & -1 \\
0 & 0 & \cdots & 0 & (\theta_2-1/2)h & -2\theta_2h & (\theta_2+1/2)h \\
\end{pmatrix}
\begin{pmatrix}
 u_0 \\
 u_1 \\
 u_2 \\
 \vdots \\
 u_{n-2} \\
 u_{n-1} \\
 u_{n}
\end{pmatrix}
=
\begin{pmatrix}
 g_a \\
 f_1 \\
 f_2 \\
 \vdots \\
 f_{n-2} \\
 f_{n-1} \\
 g_b
\end{pmatrix}.
\end{equation}

The Gauss-Seidel scheme fails to converge for this linear system when \( \theta_1 < 1/2 \), as it is not guaranteed to be positive definite. We observe that for Dirichlet boundary conditions with linear interpolation, a symmetric positive definite linear system can be obtained by solving for the ghost value  
\[
u_0 = \frac{g_a - (1 - \theta_1) u_1}{\theta_1}
\]  
and substituting it into the internal equation for \( i = 1 \), as proposed in~\cite{Gibou:Ghost}. However, this approach cannot be applied to the Neumann boundary condition or to higher-order discretizations of the Dirichlet boundary condition.

To address the issue of Gauss-Seidel failing to converge, two relaxation parameters, \( \tau_D \) and \( \tau_N \), are introduced:

\begin{equation}\label{eq:relFDM}
    \begin{aligned}
        u_0^{(m+1)} &= u_0^{(m)} + \tau_D \left( g_a - \left(\theta_1 u^{(m)}_{0} + (1-\theta_1) u^{(m)}_{1} \right) \right), \\
        u_i^{(m+1)} &= u_i^{(m)} + \frac{h^2}{2} \left( f_i - \frac{-u_{i-1}^{(m+1)}+2u_i^{(m)}-u_{i+1}^{(m)}}{h^2} \right), \quad \text{ for } \; i=1,\ldots, n-1 \\
        u_n^{(m+1)} &= u_n^{(m)} + \tau_N \left( g_b - \frac{(2\theta_2-1) u^{(m+1)}_{n-2} -4\theta_2 u^{(m+1)}_{n-1} + (2\theta_2 +1) u^{(m)}_{n}}{2h} \right).
    \end{aligned}
\end{equation}
This iterative scheme is equivalent to replace the diagonal entries of the preconditioner \( P \) related to the ghost values by the inverse of the relaxation parameters, namely \( P_{00} = \tau_D^{-1} \) and \( P_{nn} = \tau_N^{-1} \).
These parameters serve a dual purpose: first, they must be sufficiently small to ensure the convergence of the smoother ($\tau_D<1$ and $\tau_N<h$, see~\cite{COCO2013464}); second, they must be optimized to prevent boundary effects from degrading the efficiency of \ac{mg}~\cite{Coco2023}.  

\subsubsection*{Advantages of the ghost-\ac{fem}}
The ghost-\ac{fem} method used in this paper offers two advantages over the \ac{fdm}. 
First, the equations of the linear system scale with the same power of \( h \), except for the stabilization term. However, this term would scale similarly to the other equations when \( \lambda \propto h^{-1} \), which, as we will demonstrate, is the optimal scaling for the \ac{mg} approach.
The uniform scaling arises from the variational formulation, which naturally incorporates both internal equations and boundary conditions. In contrast, a rescaling can be applied to the linear system arising from \ac{fdm}~\eqref{eq:linsysfdm} by multiplying the first equation by \( h \) and the last equation by \( h^2 \). However, this rescaling is not recommended for \ac{mg} implementations, as the transfer between grids must account for the rescaling, potentially degrading \ac{mg} efficiency.

Second, the ghost-\ac{fem} method results in a symmetric and positive definite linear system, ensuring the convergence of weighted-Jacobi and Gauss-Seidel schemes. However, the stabilization parameter $\lambda$ plays a role similar to the relaxation parameter $\tau_D^{-1}$ in~\ac{fdm}. Specifically, \( \lambda \) must be sufficiently large to guarantee that the linear system remains positive definite (thus enabling the convergence of weighted-Jacobi and Gauss-Seidel), while also being optimized to maintain the efficiency of \ac{mg}.

\subsection{Transfer operators}\label{sec:1D_transfer_op}
\par In this section, we define the operators used to transfer data between the fine and coarse grids. Let the fine grid \( \Omega_{h} \) be partitioned into \( n \) intervals (where \( n \) is an even number for practical purposes, possibly a power of 2). The coarse grid \( \Omega_{2h} \) then consists of \( \left(\frac{n}{2} +1 \right) \) nodes.

\subsubsection{Brief overview of the Splitting Strategy in FDM}
First, we recall the restriction operator \( \vec{r}_{2h} = \stackrel{FDM}{\mathcal{I}^h_{2h}} \vec{r}_{h} \) used in the \ac{fdm} context. Several approaches exist, but the two most commonly used ones are the injection operator, where \( \vec{r}_{2h,i} = \vec{r}_{h,2i} \), and the Full-Weighting operator~\cite{trottenberg2000multigrid}:
\begin{equation*}
\begin{aligned}
\vec{r}_{2h,i} &= 0.25 (\vec{r}_{h,2i-1} + 2\vec{r}_{h,2i} + \vec{r}_{h,2i+1}), \quad \text{for } i=1,\ldots,n/2-1, \\
\vec{r}_{2h,0} &= 0.25 (2\vec{r}_{h,0} + \vec{r}_{h,1}), \\
\vec{r}_{2h,n/2} &= 0.25 (\vec{r}_{h,n-1} + 2\vec{r}_{h,n}).
\end{aligned}
\end{equation*}
We focus on the Full-Weighting operator, which provides better accuracy and smoother coarse-grid corrections compared to the simple injection operator. The latter merely copies values from the fine to the coarse grid, potentially leading to information loss and poor approximations.

The matrix representation of Full-Weighting operator is:
\[
\stackrel{FDM}{\mathcal{I}^h_{2h}} =
\frac{1}{4}
\begin{pmatrix}
 2 & 1 & 0 & 0 & 0 & 0 & \cdots & 0 \\
 0 & 1 & 2 & 1 & 0 & 0 & \cdots & 0 \\
 0 & 0 & 0 & 1 & 2 & 1 & \cdots & 0 \\
 \vdots & \vdots & \vdots & \vdots & \vdots & \vdots & \ddots & \vdots \\
0 & 0 & 0 & \cdots & 1 & 2 & 1 & 0 \\
0 & 0 & 0 & \cdots & 0 & 0 & 1 & 2 \\
\end{pmatrix}.
\]
From Eq.~\eqref{eq:relFDM}, we observe that the residuals of internal equations and boundary conditions scale with different powers of \( h \). As a result, the standard Full-Weighting approach can lead to convergence issues in the \ac{mg} method. To address this, the authors in~\cite{COCO2013464} proposed to split the contributions from the residuals of the internal equations \( \vec{r}^I \) and the boundary conditions \( \vec{r}^D \) and \( \vec{r}^N \). 
We call it the {\em splitting strategy}.
Specifically, the residuals of the boundary conditions are defined as \( r^D = r_0 \) and \( r^N = r_n \), while the residual vector of internal equations is given by \( r^I_i = r_i \) for \( 1 \leq i \leq n-1 \). The ghost values are then extrapolated using constant extrapolation: \( r^I_0 = r_1 \) and \( r^I_n = r_{n-1} \). 

Finally, the defect of internal equation is restricted from the fine to the coarse grid using the standard Full-Weighting approach:
\[
\vec{r}_{2h,i} = 0.25 (\vec{r}_{h,2i-1} + 2\vec{r}_{h,2i} + \vec{r}_{h,2i+1}), \quad \text{for } i=1,\ldots,n/2-1.
\]
The defect vector \( \vec{r}_{2h} \) is then completed by incorporating the residuals of the boundary conditions, \( r_{2h,0} = r^D \) and \( r_{2h,n/2} = r^N \). The residual equation is subsequently solved (Step~\ref{alg:state:reseq} of Alg.~\ref{alg:2Grid_scheme}):
\begin{equation}\label{eq:reseq2h}
A_{2h} \vec{e}_{2h} = \vec{r}_{2h},
\end{equation}
where \( A_{2h} \) is the matrix obtained from the discretization of the Poisson problem on the coarse grid \( \Omega_{2h} \). 

For a given approximation \( \hat{u} \), the residual equation~\eqref{eq:reseq2h} can be interpreted as the discretization on $\Omega_{2h}$ of the continuous problem:
\begin{equation}\label{eq:respoisson}
    -\boldsymbol{\Delta }e = r \text{~ in ~}\Omega, \quad e = r^D \text{~ on ~} \Gamma_D, \quad  \boldsymbol{n}\cdot \boldsymbol{\nabla}e = r^N \text{~ on ~} \Gamma_N.
\end{equation}
where \( r = f + \Delta \hat{u} \), \( r^D = g_a - \hat{u}(a) \), and \( r^N = g_b - \hat{u}'(b) \). Then, the solution of the original problem~\eqref{eq:poisson_model_prob} is given by \( u = \hat{u} + e \). 
This strategy ensures optimal convergence properties, as demonstrated in~\cite{Coco2023}.

\subsubsection{Transfer operators in ghost-FEM}\label{sect:restr_fem}
Now, we derive the restriction operator in the context of ghost-\ac{fem}. Let $\phi_{h,i}$, for $i \in \{0,1,\dots,n\}$, and $\phi_{2h,i}$, for $i \in \{0,1,\dots,\frac{n}{2}\}$, be the shape functions of $\Omega_{h}$ and $\Omega_{2h}$, respectively. The relationship between the shape functions on the two grids is given by
\begin{equation}\label{eq:shape_fun_relation}
    \phi_{2h,i} = \frac{1}{2} \phi_{h, 2i -2} + \phi_{h, 2i-1} + \frac{1}{2} \phi_{h, 2i}, \qquad i = 0,1,\dots,\frac{n}{2}.
\end{equation}
We observe that \eqref{eq:shape_fun_relation} is an exact representation of the shape function in $\Omega_{2h}$, not an approximation as in the \ac{fdm} context.

Thus, the matrix corresponding to the restriction operator $\mathcal{I}^h_{2h} : \mathbb{R}^{n+1} \to \mathbb{R}^{\frac{n}{2}+1}$ is defined as
\begin{equation}\label{eq:resop}
    \mathcal{I}^h_{2h} = \frac{1}{2}
    \begin{pmatrix}
        2 & 1 & 0 & 0 & 0 & 0 & \cdots & 0 \\
        0 & 1 & 2 & 1 & 0 & 0 & \cdots & 0 \\
        0 & 0 & 1 & 2 & 1 & 0 & \cdots & 0 \\
        \vdots & \vdots & \vdots & \vdots & \vdots & \vdots & \ddots & \vdots \\
        0 & 0 & 0 & \cdots & 0 & 1 & 2 & 1  \\
        0 & 0 & 0 & \cdots & 0 & 0& 1 & 2 
    \end{pmatrix}.
\end{equation}

Here, the splitting strategy that was adopted in the \ac{fdm} context, namely separating the contributions of internal and boundary residuals and then deriving the residual equation from \eqref{eq:respoisson}, is unnecessary. In fact, we will prove that the splitting strategy is equivalent to the standard full-weighting approach
$\vec{r}_{2h} = \mathcal{I}^h_{2h} \vec{r}_h$,
with $\mathcal{I}^{h}_{2h}$ given by \eqref{eq:resop}. This is a significant advantage of ghost methods in \ac{fem} over \ac{fdm}, as the discretization is derived from a variational formulation that inherently incorporates boundary conditions, rather than requiring separate treatment as in \ac{fdm}.

First, we derive the splitting strategy for ghost-\ac{fem}.
Given a function $\hat{u} \in \mathcal{V}_{h,act}$, the residual equation of~\eqref{eq:disc_formulation} is to find $e_h \in \mathcal{V}_{h,act}$ such that:
\begin{align}\label{eq:disc_formulation_res}
        a_h(e_h,v_h) = l^\text{res}_h(v_h) \quad \forall ~v_h \in \mathcal{V}_{h,act},
\end{align}
where 
\begin{multline}
    l^\text{res}_h(v) \doteq l_h(v) - a_h(\hat{u},v) \\
    = \sum_{K \in \mathcal{T}_{h,act}} \int_{K \cap \Omega} r ~v \mathrm{~d}\boldsymbol{x} + \int_{K \cap\Gamma_D} \lambda_K r_D ~v   -\left(\boldsymbol{n} \cdot \boldsymbol{\nabla} v \right) r_D \mathrm{~d}\boldsymbol{S}  + \int_{K \cap \Gamma_N} r_N ~ v \mathrm{~d}\boldsymbol{S}. \label{eq:disc_rhs_res}
\end{multline}
Here, the residual boundary conditions are defined as $r_D \doteq g_D - \hat{u}(a)$ and $r_N \doteq g_N - \hat{u}'(b)$, while the residual function $r$ is defined in such a way that, for any $v \in \mathcal{V}_{h,act}$, we have: 
\[
\sum_{K \in \mathcal{T}_{h,act}}
\int_{K \cap \Omega} r ~v \mathrm{~d}\boldsymbol{x} =
\sum_{K \in \mathcal{T}_{h,act}} \left(
\int_{K \cap \Omega} f ~v \mathrm{~d}\boldsymbol{x} - 
\int_{K \cap \Omega} \boldsymbol{\nabla} \hat{u} \cdot \boldsymbol{\nabla} v \mathrm{~d}\boldsymbol{x} + 
\int_{K \cap \partial \Omega} \left(\boldsymbol{n} \cdot \boldsymbol{\nabla} \hat{u} \right) v \mathrm{~d}\boldsymbol{S}
\right).
\]
In matrix form:
\begin{equation}\label{eq:rhI}
\vec{r}^I_{h} \doteq \vec{F}^I_h - A_h^I \hat{\vec{u}}_h - A_h^{B} \hat{\vec{u}}_h - A_h^{B,N} \hat{\vec{u}}_h
\end{equation}
where
\[
A^{B,N}_h = 
\frac{1}{h}
\begin{pmatrix}
 0 & \ldots & 0 & 0 & 0 \\
 \vdots & \ddots & \vdots & \vdots & \vdots \\
 0 & \ldots & 0 & 0 & 0 \\
 0 & \ldots & 0 & 1-\theta_2  & \theta_2-1 \\
0 & \ldots & 0 & \theta_2 & -\theta_2 \\
\end{pmatrix}.
\]
Therefore, the discrete residual problem on $\Omega_{2h}$ is obtained by taking $\vec{r}^I_{2h} = \mathcal{I}^h_{2h} \vec{r}^I_h$ and solving $A_{2h} \vec{e}_{2h} = \vec{r}^\text{SPLIT}_{2h}$, where $A_{2h}$ is the coarse grid operator (described in Section~\ref{sec:1D_coarse_op}) and  
\begin{equation}\label{eq:rsplit}
\vec{r}^\text{SPLIT}_{2h} = 
\vec{r}^I_{2h} + \vec{r}^B_{2h} + \vec{r}^\lambda_{2h} + 
\vec{r}^N_{2h}, 
\end{equation}
where $\vec{r}^B_{2h}$, $\vec{r}^\lambda_{2h}$ and $\vec{r}^N_{2h}$ are obtained as in \eqref{eq:ls1dA}-\eqref{eq:ls1dF} by replacing $\theta_1$ and $\theta_2$ with the corresponding coarse grid values, that are
$(1+\theta_1)/2$ and $(1+\theta_2)/2$, respectively. Then:
\begin{align*}
    \boldsymbol{r}^B_{2h} = \frac{1}{2h}
    \begin{pmatrix}
        -r_a \\
        r_a \\
        0 \\ 
        \vdots \\
        0
    \end{pmatrix}, \quad
    \boldsymbol{r}^\lambda_{2h} = 
    \frac{\lambda}{2} 
    \begin{pmatrix}
        (\theta_1+1) ~r_a\\
        (1-\theta_1) ~r_a \\ 
        0 \\
        \vdots \\ 
        0   
    \end{pmatrix}, \quad
     \boldsymbol{r}^N_{2h} =  
     \frac{1}{2}
    \begin{pmatrix}
        0\\
        \vdots \\
        0\\
        (1-\theta_2) ~ r_b \\ 
       (\theta_2+1) ~ r_b  
    \end{pmatrix}.
\end{align*} 

\begin{prop}
In the context of ghost-\ac{fem}, the splitting and the standard strategies are equivalent, namely $\vec{r}^\text{SPLIT}_{2h} = \vec{r}_{2h}$, where 
$\vec{r}^\text{SPLIT}_{2h}$ is given by \eqref{eq:rsplit} and $\vec{r}_{2h} = \mathcal{I}^h_{2h} \vec{r}_h$ corresponds to the standard (non-splitting) approach.
\end{prop}
{\em Proof.} From \eqref{eq:ls1dA}-\eqref{eq:ls1dF}, we have:
\[
\vec{r}_h = \vec{F}_h-A_h \vec{\hat{u}}_h
= \boldsymbol{F}_h^I + \boldsymbol{F}_h^B +\boldsymbol{F}_h^\lambda+\boldsymbol{F}_h^N
- (A_h^I+A_h^{B} + (A_h^{B})^T + A_h^{\lambda}) \vec{\hat{u}}_h
= 
\vec{r}_h^I + \vec{\tilde{r}}_h^B + \vec{\tilde{r}}_h^\lambda + \vec{\tilde{r}}_h^N,
\]
where $\vec{r}_h^I$ is given by \eqref{eq:rhI} and
$\vec{\tilde{r}}_h^B \doteq \vec{F}_h^B - (A_h^B)^T \hat{u}_h$,
$\vec{\tilde{r}}_h^\lambda \doteq \vec{F}_h^\lambda - A_h^\lambda \hat{u}_h$,
$\vec{\tilde{r}}_h^N \doteq \vec{F}_h^N + A_h^{B,N} \hat{u}_h$.
First, note that
\begin{equation}\label{eq:rtilde}
\vec{\tilde{r}}_h^B = \vec{r}_h^B, \quad
\vec{\tilde{r}}_h^\lambda = \vec{r}_h^\lambda, \quad 
\vec{\tilde{r}}_h^N = \vec{r}_h^N,
\end{equation}
where:
\begin{align*}
    \boldsymbol{r}^B_{h} = \frac{1}{h}
    \begin{pmatrix}
        -r_a \\
        r_a \\
        0 \\ 
        \vdots \\
        0
    \end{pmatrix}, \quad
    \boldsymbol{r}^\lambda_{h} = 
    \lambda
    \begin{pmatrix}
        \theta_1 ~r_a\\
        (1-\theta_1) ~r_a \\ 
        0 \\
        \vdots \\ 
        0   
    \end{pmatrix}, \quad
     \boldsymbol{r}^N_{h} =  
    \begin{pmatrix}
        0\\
        \vdots \\
        0\\
        (1-\theta_2) ~ r_b \\ 
       \theta_2 ~ r_b  
    \end{pmatrix}.
\end{align*} 
For example, the first one is obtained by observing that:
$
(A_h^B)^T \hat{u}_h =
(
        -\hat{u}'(a), 
        \hat{u}'(a),
        0, 
        \ldots,
        0
 )^T
$.
Then, we observe that
\begin{equation}\label{eq:rhto2h}
\vec{r}^B_{2h} = \mathcal{I}^h_{2h} \vec{r}^B_h, \quad
\vec{r}^\lambda_{2h} = \mathcal{I}^h_{2h} \vec{r}^\lambda_h, \quad
\vec{r}^N_{2h} = \mathcal{I}^h_{2h} \vec{r}^N_h.
\end{equation}
This holds because equation~\eqref{eq:shape_fun_relation} provides an exact representation of the shape function in the coarse grid through shape functions in the fine grid. Combining \eqref{eq:rtilde} and \eqref{eq:rhto2h}, we obtain $\vec{r}^\text{SPLIT}_{2h} = \vec{r}_{2h}$. $\square$

It is worth noting that the restriction operators in \ac{fem} and \ac{fdm} differ by a constant scaling factor, namely $\stackrel{\text{FDM}}{\mathcal{I}^h_{2h}} = 2 \mathcal{I}^h_{2h}$. This reflects the differences in how residuals are formulated in the two methods. In \ac{fdm}, residuals are computed at discrete points, and restriction is typically defined through simple averaging schemes such as full-weighting. In contrast, the \ac{fem} residual arises from a variational formulation, where test functions are integrated over the domain. This process introduces a mesh-dependent scaling, and as a result, the \ac{fem} restriction operator reflects this scaling when transferring residuals to the coarse grid.

Finally, the prolongation operator $\mathcal{I}_h^{2h} : \mathbb{R}^{\frac{n}{2}+1} \to \mathbb{R}^{n+1}$ is given by the transpose of the restriction operator matrix $\mathcal{I}_h^{2h} = \left(\mathcal{I}^h_{2h}\right)^T$.
This operator is identical to the classical linear interpolation from the coarse grid to the fine grid.

\subsection{Coarse grid operator} \label{sec:1D_coarse_op}
\par The operator $A_{2h}$ on the coarse grid can be determined in two ways. The first approach utilizes the transfer operators introduced in \cref{sec:1D_transfer_op}. In this case, the coarse grid operator $A_{2h}$ is defined by the Galerkin condition:
\begin{equation}\label{eq:galerkin}
A_{2h} = \mathcal{I}_{h}^{2h} A_h ~\mathcal{I}^{h}_{2h} .
\end{equation}

\par Alternatively, $A_{2h}$ can be constructed directly using the shape functions on the coarse grid, that is,
\begin{align*}
    A_{2h,i,j} = &\int_{a}^{b} \boldsymbol{\nabla} \phi_{2h,i} \cdot \boldsymbol{\nabla} \phi_{2h,j} + \lambda ~\phi_{2h,i}(a)\phi_{2h,j}(a) + \phi_{2h,i}^\prime(a)\phi_{2h,j}(a) + \phi_{2h,i}(a)\phi_{2h,j}^\prime(a) .
\end{align*}
However, regardless of the approach used, the matrix $A_{2h}$ remains the same due to the relationship between the shape functions of the fine and coarse grids, as given in \cref{eq:shape_fun_relation}.

\subsection{Optimal stabilization parameter for MG efficiency} \label{sec:optimal_par_1D}

\par The stabilization parameter $\lambda >0$ has to be chosen sufficiently large to ensure the well-posedness of the weak formulation of the ghost-\ac{fem}. However, if the parameter is too large, then the convergence of the \ac{mg} deteriorates significantly, as we will show in Section \ref{sec:1D_tests}. Therefore, an optimal choice for the stabilization parameter has to be made to guarantee the coercivity of the method and to ensure optimal convergence of the \ac{mg}. 

\par The bilinear form \cref{eq:disc_bil_form} is coercive if there exists a positive constant $\tilde{C}>0$ such that  $a_h(v_h,v_h) \geq \tilde{C}\|v_h\|_{\mathcal{V}_h}$ for any $v_h \in \mathcal{V}_h$. Here, the norm on the \ac{fe} space $\mathcal{V}_h$ is defined as 
\begin{equation}\label{eq:norm_vh}
    \|v_h\|_{\mathcal{V}_h}^2 \doteq \|\boldsymbol{\nabla }v_h \|^2_{L^2(\Omega)} + \lambda \|v_h\|^2_{L^2(\Gamma_D)}.
\end{equation} 

To establish a lower bound on the bilinear form $a_h(v_h, v_h)$, we apply Young’s inequality, which states that for any two real numbers $a$ and $b$ and a positive parameter $\varepsilon$ we have
$ab \leq a^2/(2\varepsilon) + \varepsilon b^2/2$.
Applying it to the boundary integral term, we obtain: 
\[
\int_{\Gamma_D} 2 v_h(\boldsymbol{n\cdot\nabla} v_h) \mathrm{~d} \boldsymbol{S} \leq 
\frac{1}{\varepsilon}\int_{\Gamma_D} v_h v_h \mathrm{~d} \boldsymbol{S}
+ \varepsilon \int_{\Gamma_D} (\boldsymbol{n\cdot\nabla} v_h) (\boldsymbol{n\cdot\nabla} v_h) \mathrm{~d} \boldsymbol{S}.
\]
Therefore:
\begin{align*}
    a_h(v_h,v_h) &= \int_{\Omega} \boldsymbol{\nabla} v_h\cdot\boldsymbol{\nabla} v_h \mathrm{~d}\boldsymbol{x} + \int_{\Gamma_D} \lambda v_h v_h - 2 v_h(\boldsymbol{n\cdot\nabla} v_h) \mathrm{~d} \boldsymbol{S} \\ 
    &\ge \int_{\Omega} \boldsymbol{\nabla} v_h\cdot\boldsymbol{\nabla} v_h \mathrm{~d}\boldsymbol{x}  -  \varepsilon \int_{\Gamma_D}(\boldsymbol{n}\cdot \boldsymbol{\nabla} v_h) (\boldsymbol{n}\cdot \boldsymbol{\nabla} v_h)  \mathrm{~d} \boldsymbol{S} + \left(\lambda -\frac{1}{\varepsilon}\right)\int_{\Gamma_D}  v_h v_h  \mathrm{~d} \boldsymbol{S}\\
    &\ge \left( 1- C\varepsilon \right)\int_{\Omega} \boldsymbol{\nabla} v_h\cdot\boldsymbol{\nabla} v_h \mathrm{~d}\boldsymbol{x} + \left(\lambda -\frac{1}{\varepsilon}\right)\int_{\Gamma_D}  v_h v_h  \mathrm{~d} \boldsymbol{S}, 
\end{align*}
where in the last expression we assume that there exists a constant $C>0$ such that 
\begin{equation}\label{eq:lowC}
C \int_{\Omega} \boldsymbol{\nabla} v_h\cdot\boldsymbol{\nabla} v_h \mathrm{~d}\boldsymbol{x} \ge \int_{\Gamma_D}(\boldsymbol{n}\cdot \boldsymbol{\nabla} v_h) (\boldsymbol{n}\cdot \boldsymbol{\nabla} v_h)  \mathrm{~d} \boldsymbol{S}.
\end{equation}
The bilinear form is coercive if $1-C\varepsilon > 0 $ and $\lambda - \frac{1}{\varepsilon} > 0$. Hence, a necessary condition for the stabilization parameter is $\lambda > C$.

\par The first step is to find an estimate for the smallest value of $C$ such that \eqref{eq:lowC} is satisfied.
This is obtained by solving the generalized global eigenvalue problem 
\begin{equation}\label{eq:geneigprob}
\boldsymbol{K} v = \Lambda \boldsymbol{M} v,
\end{equation}
where $\boldsymbol{K}_{ij} = \int_{\Gamma_D} (\boldsymbol{n}\cdot \boldsymbol{\nabla} \phi_i)(\boldsymbol{n}\cdot \boldsymbol{\nabla} \phi_j)~\mathrm{d}\boldsymbol{S}$ and $\boldsymbol{M}_{ij} = \int_{\Omega}  \boldsymbol{\nabla} \phi_i \cdot  \boldsymbol{\nabla} \phi_j~\mathrm{d}\boldsymbol{x}$ for $i,j \in \{0,1,\dots,n\}$. The lower bound of $C$ is then estimated by taking $C = \max \Lambda$~\cite{Saberi2023, Astuto2024Ghost}. The stabilization parameter is then $\lambda = \gamma C$, for some constant $\gamma>1$.

In the one-dimensional case addressed in this section, there is only one cut-cell that intersects the Dirichlet boundary, that is the cell $[a,x_1]$. Therefore, the estimate for the lower bound of $C$ can be obtained by solving the $2 \times 2$ generalized local eigenvalue problem \eqref{eq:geneigprob}
where $\boldsymbol{K}_{ij} = \phi_i^\prime(a) \phi_j^\prime(a)$ and $\boldsymbol{M}_{ij} = \int_{a}^{x_1} \phi_i^\prime(x)  \phi_j^\prime(x) ~\mathrm{d}{x}$ for $i,j \in \{0,1\}$, that are:
\begin{align*}
    M= \frac{\theta_1}{h}
    \begin{pmatrix}
       1  & -1 \\ 
       -1 & 1
    \end{pmatrix},
    \quad K= \frac{1}{h^2}
    \begin{pmatrix}
        1 & -1 \\
        -1 & 1
    \end{pmatrix}.
\end{align*}

We observe that $\det(K-\lambda M)=0$ for any $\lambda \in \mathbb{C}$, which prevents the use of the classical definition of generalized eigenvalues.
Instead, we adopt the definition proposed in~\cite{hochstenbach2019solving}, which states that
$\lambda_0 \in \mathbb{C}$ is a generalized eigenvalue if $\text{rank} (K-\lambda_0 M) < \displaystyle \max_{\xi \in \mathbb{C}} \text{rank}(K - \xi M)$.
By this definition we have $C = \frac{1}{\theta_1 h}$ and, consequently, the stabilization parameter is $\lambda = \gamma \cdot C = \frac{\gamma}{\theta_1 h}$, where $\gamma > 1$.

\subsection{Numerical tests in 1D}\label{sec:1D_tests}

\par As per the problem set up at the beginning of Section \ref{sec:1D_mg}, we choose $\Omega_{art} = [0,1]$ and $\Omega = (a,b)$, where $a = (1-\theta_1)h$ and $b = 1-(1-\theta_2)h$. The Dirichlet and Neumann boundaries are $\Gamma_D = \{a\}$ and $\Gamma_N = \{b\}$ respectively. We fix $\theta_2 = 10^{-2}$ and $\theta_1 \in \{0.0099,0.05,0.1,0.2,0.3,0.4,0.5,0.75,0.9,0.99,1.0\}$. The spatial mesh size is $h =2^{-m},~ m \in \{7,8,9,10\}$. The snapping threshold is $h^{-\alpha}$ with $\alpha = 2$.

We perform $\nu_1 = 2$ pre-smoothing iterations and  $\nu_2 = 1$ post-smoothing iteration.
The convergence factor is calculated as $\rho = \frac{1}{10}\Sigma_{m=41}^{50} \rho^{(m)}$.
For the chosen values of $\nu_1$ and $\nu_2$, the predicted convergence factor of the \ac{mg} method for the 2D \ac{fdm} applied to the case $\Omega=\Omega_{art}$ (namely the standard case where the boundaries coincide with Cartesian lines) is $\rho \approx 0.119$~\cite[Table 5.2]{trottenberg2000multigrid}. For the 1D case, the convergence factor is $\rho \approx 0.1$~\cite{COCO2013464}.
In order to prevent numerical issues that arise when the residual reaches the machine precision, we solve the homogeneous problem, that is, $f=g_a=g_b=0$, then 
the actual solution is $u(x) = 0$. Initial guess of the \ac{mg} iteration scheme must be different from zero. We choose $\vec{u}^{(0)}_h=\vec{1}_h$. 
 
\par
\begin{figure}
 \centering
  \begin{subfigure}[t]{.5\textwidth}
   \centering
   \includegraphics[width=\textwidth]{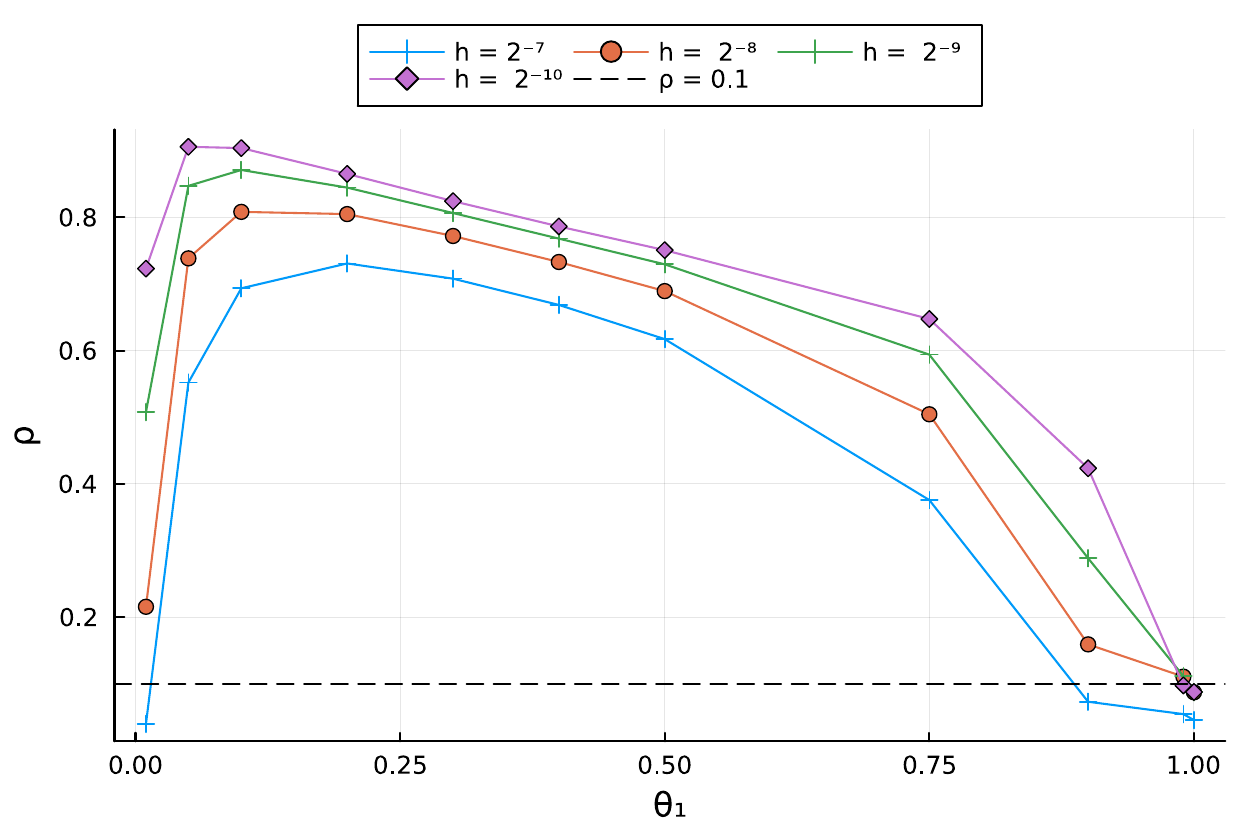}
   \caption{\ac{tgcs} with $\lambda = h^{-2}$.}
   \label{fig:1D_two_grid_not_opt_lam}
 \end{subfigure}%
 \begin{subfigure}[t]{.5\textwidth}
   \centering
   \includegraphics[width=\textwidth]{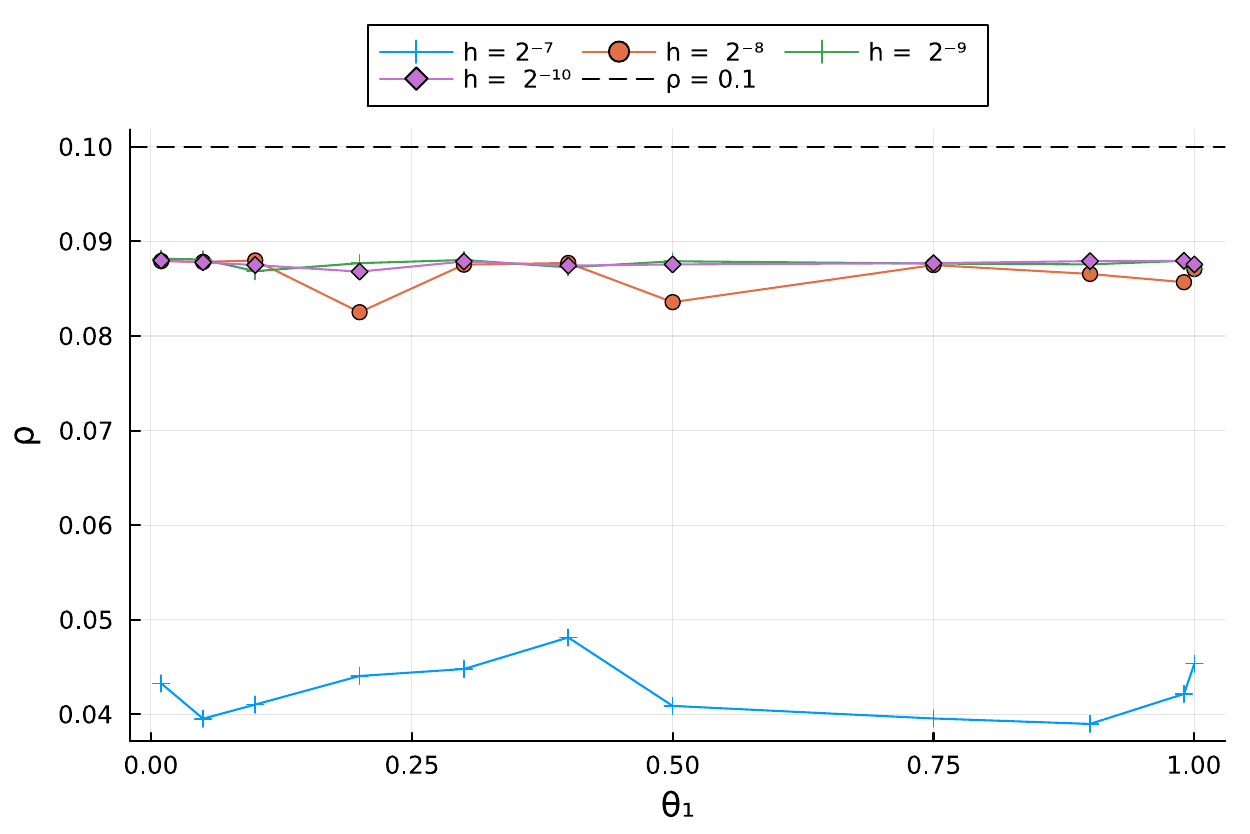}
    \caption{\ac{tgcs} with  $\lambda = 1.1~ (\theta_1 h)^{-1}$.}
   \label{fig:1D_two_grid}
 \end{subfigure}%
  \vskip\baselineskip
 \begin{subfigure}[t]{.5\textwidth}
   \centering
   \includegraphics[width=\textwidth]{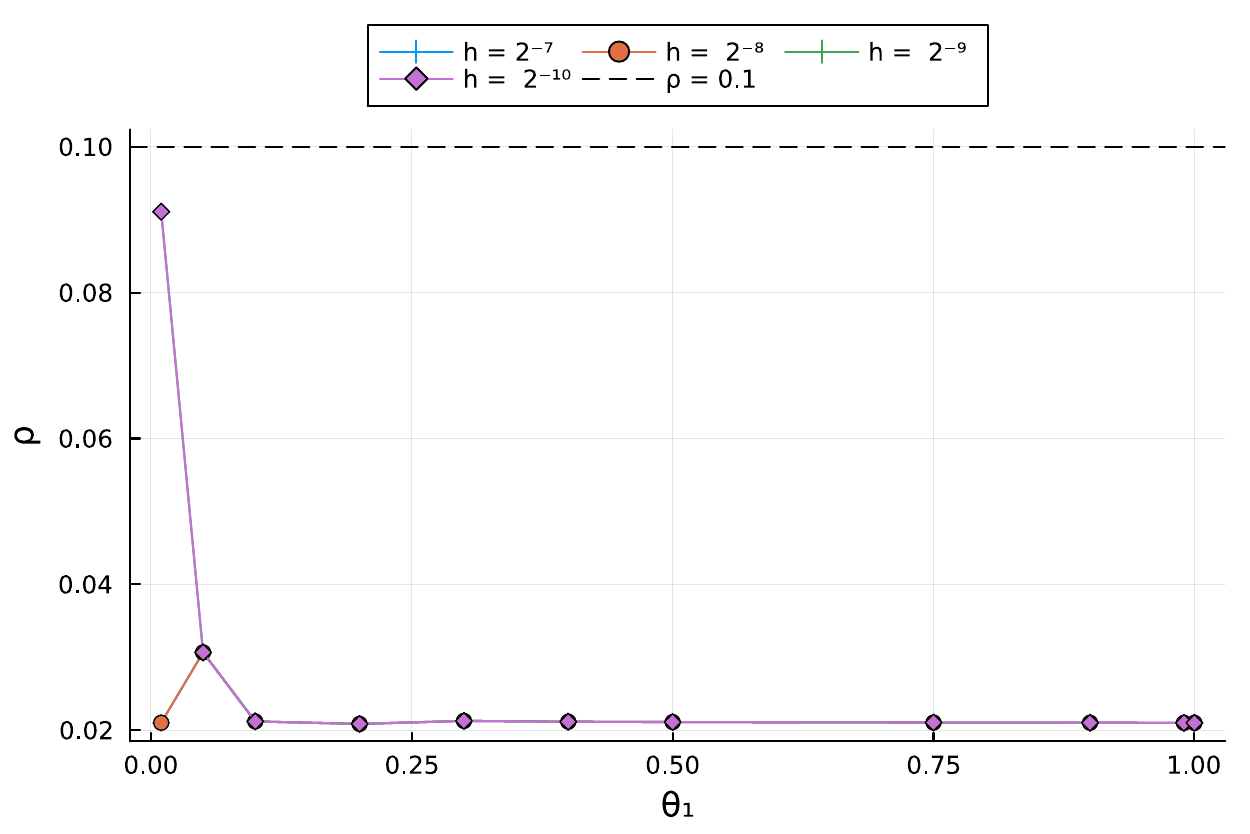}
    \caption{W-cycle with $\lambda = 1.1~ (\theta_1 h)^{-1}$.}
    \label{fig:1D_w_cycle}
 \end{subfigure}
  \caption{Plots of the convergence factor $\rho$ against the position of the Dirichlet boundary $\theta_1$ for different mesh sizes $h$ using the \ac{tgcs} and W-cycle.}
  \label{fig:1D_tests}
  \end{figure}

 In the first numerical test we set $\lambda = h^{-2}$ and compute the convergence factor of a \ac{tgcs}. The plots of the convergence factor against different Dirichlet boundary positions for various mesh sizes are illustrated in \cref{fig:1D_two_grid_not_opt_lam}. This choice of $\lambda$ ensures the coercivity of the bilinear form \eqref{eq:disc_bil_form} and the ghost-\ac{fem} remains second-order accurate \cite{Astuto2024Ghost}. However, the \ac{mg} convergence is not optimal, since the convergence factor is far from the optimal one, that is, $\rho \approx 0.1$. This highlights the importance of carefully selecting the stabilization parameter to ensure both the accuracy of the method and the efficiency of the \ac{mg} solver.

\par In the next numerical simulation, we compute the convergence factor for a \ac{tgcs} using the optimal stabilization parameter suggested in Section~\ref{sec:optimal_par_1D}, namely $\lambda = \frac{1.1}{\theta_1 h}$. \cref{fig:1D_two_grid} shows the plot of the convergence factor ($\rho$) against the position of the Dirichlet boundary ($\theta_1$) for different mesh sizes. 

\par
\begin{figure}
    \centering
    \begin{subfigure}[t]{.5\textwidth}
   \centering
   \includegraphics[width=\textwidth]{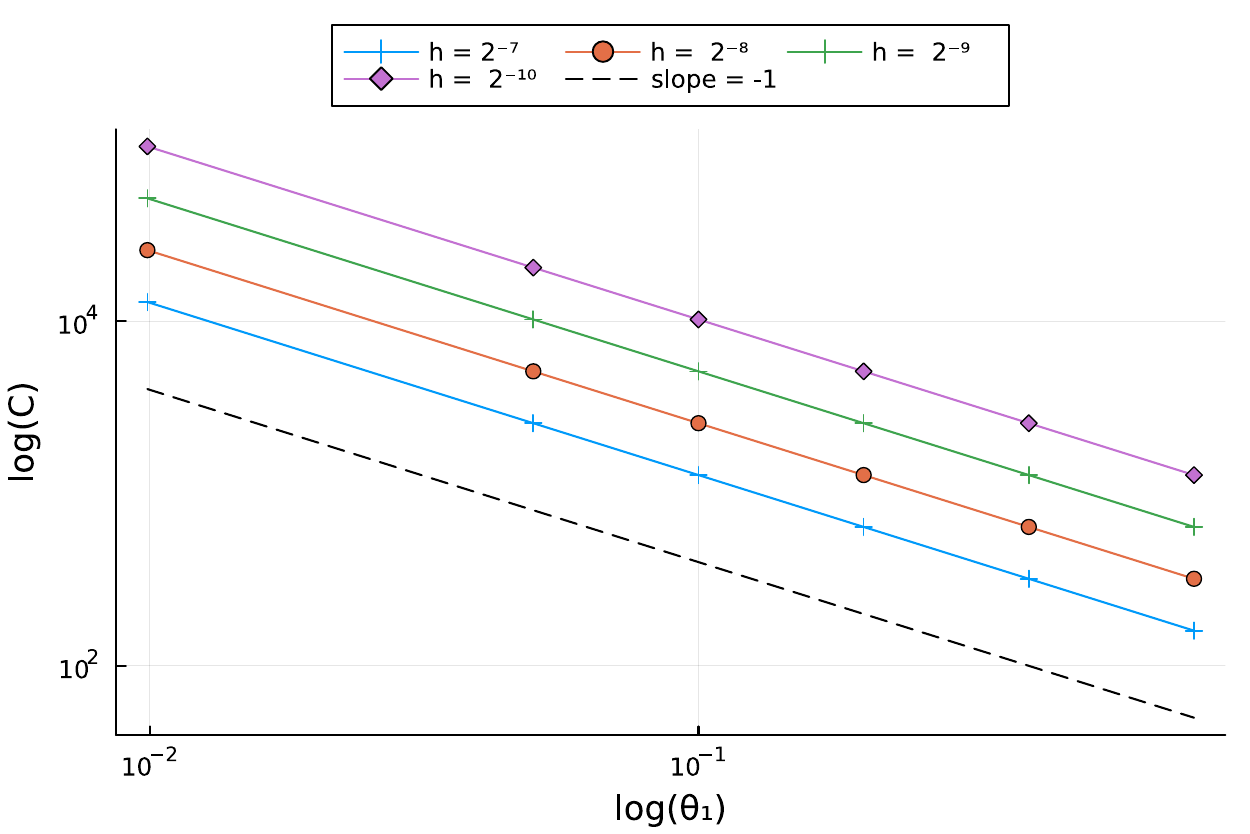}
   \caption{$C$ as a function of $\theta_1$}
   \label{fig:1D_global_lambda_theta}
 \end{subfigure}%
 \begin{subfigure}[t]{.5\textwidth}
   \centering
   \includegraphics[width=\textwidth]{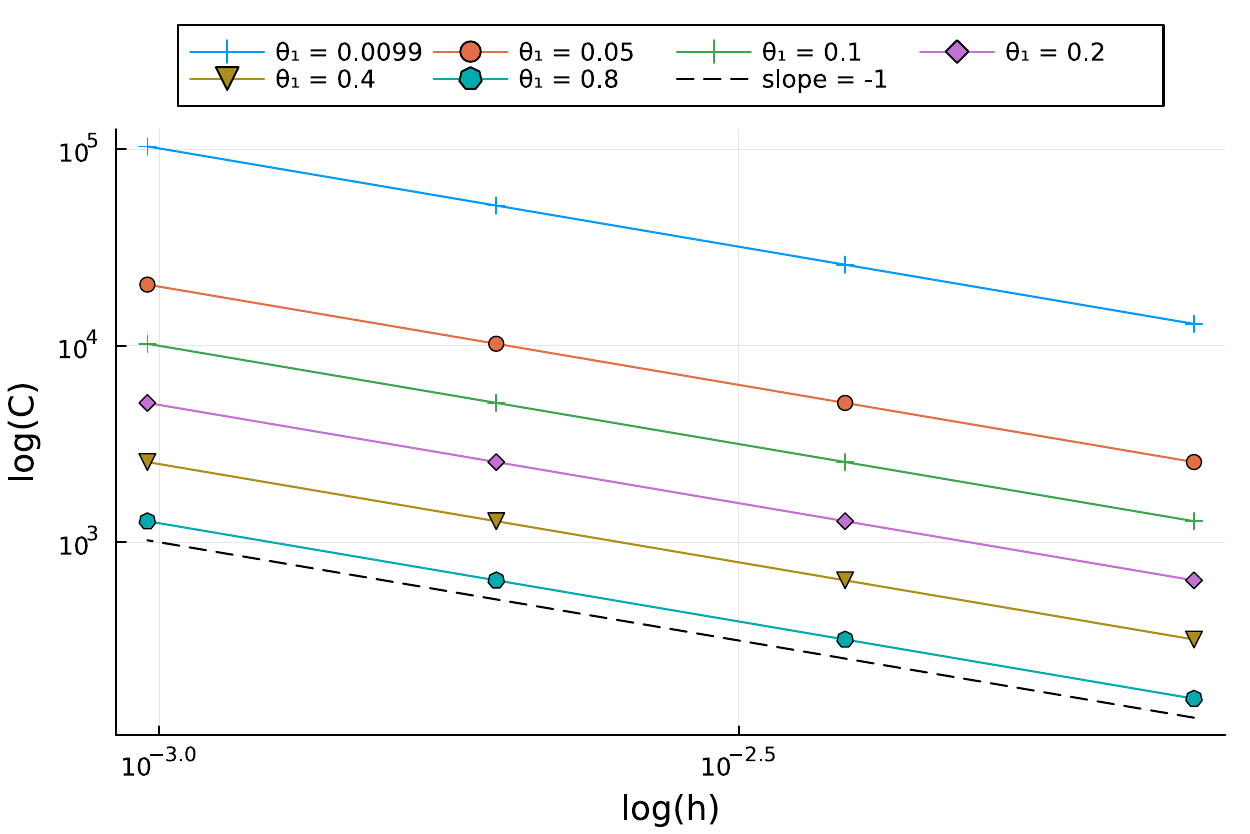}
   \caption{$C$ as a function of $h$}
   \label{fig:1D_global_lambda_h}
 \end{subfigure}%
    \caption{Plots of $C$ for a global generalized eigenvalue problem against the Dirichlet boundary configuration $\theta_1$ and mesh size $h$ on a log-log scale.}
    \label{fig:1D_global_lambda}
\end{figure}

 In \cref{fig:1D_global_lambda_theta,fig:1D_global_lambda_theta} we compute the lower bound for $C$ by numerically solving the global eigenvalue problem \eqref{eq:geneigprob}. The numerical value of $C$ is plotted as a function of the Dirichlet boundary position $\theta_1$ and mesh size $h$ on a log-log scale, confirming that $C = \frac{1}{\theta_1 h}$. This result implies that, for the one-dimensional case, the global and local eigenvalue problems yield identical solutions.

\par Finally, we compute the convergence factor for a W-cycle. The coarsest grid corresponds to a spatial step of $h=2^{-3}$, and the stabilization parameter is $\lambda = \frac{1.1}{\theta_1 h}$. \cref{fig:1D_w_cycle} illustrates the plot of the convergence factor ($\rho$) against the position of the Dirichlet boundary ($\theta_1$) for different mesh sizes for the W-cycle.

\par We observe that selecting the optimal value for the stabilization parameter is essential to ensure optimal \ac{mg} convergence, that is, $\rho \approx 0.1$.

\section{Two-dimensional case} \label{sec:2D_mg}

\par In this section, we introduce key components of a \ac{mg} method for 2D geometries. As before, we use the Gauss-Seidel method as the smoother. The operators for transferring data between fine and coarse grids are defined in \cref{sec:2D_transfer_op}. As in \cref{sec:1D_coarse_op}, the operator on the coarse grid is defined using the restriction and prolongation operators by the Galerkin condition~\eqref{eq:galerkin}. 
While in 1D this formulation was identical to discretizing the problem from scratch in the coarse grid $\Omega_{2h}$, 
using the Galerkin condition to construct coarse-level operators in higher dimensions does not generally yield the same result as assembling the operators directly on each grid level. This discrepancy arises from differences in how the geometry (and particularly the boundary) is represented at different resolutions. The Galerkin approach projects the fine-grid geometry onto the coarse level, which may not coincide exactly with the cut configuration defined by the coarse grid level-set representation. In contrast, rebuilding the operator from scratch ensures consistency with the geometry at each level but sacrifices Galerkin consistency. However, the resulting differences are accurate to second-order and to do not contrast to the accuracy of the overall method. The choice of the stabilization parameter is described in \cref{sec:2D_optimal_lambda}.

\subsection{Transfer operators}\label{sec:2D_transfer_op}

\par In this section, we define the operators used to transfer data between the coarse and fine grids. Let $\Omega_{h}$ denote the fine grid, which consists of  $n \times n$ cells. Consequently, there are $(n+1)^2$ shape functions on the fine grid and $(\frac{n}{2} + 1)^2$ shape functions on the coarse grid. To facilitate the data transfer, we first define the operator $\mathcal{N}^{2h}_h : \Omega_{2h} \rightarrow  \Omega_{h} $ that maps the nodes of the coarse grid to those of the fine grid. In particular, $\mathcal{N}^h_{2h}((x_{i},y_{j})) = (x_{2i},y_{2j})$, for $i,j \in \{0,1,\dots,\frac{n}{2}\}$. The single index of a node $(x_i,y_j)$ on the coarse grid is computed as $k = i+j~(\frac{n}{2}+1)$, where $i,j \in \{0,1,\dots,\frac{n}{2}\}$.
Similarly, the single index of a node $(x_i,y_j)$ on the fine grid is computed as $k = i+j~(n+1)$, where $i,j \in \{0,1,\dots,n\}$. 
Using this notation, we denote $\boldsymbol{x}_k \doteq(x_i,y_j)$.

\par Let $\phi_{h,k}$ and  $\phi_{2h,k}$ denote the shape functions on the fine and coarse grids, respectively. The exact relationship between them is given by
\begin{align*}
    \phi_{2h,k} &= \frac{1}{4}\phi_{h,\mathcal{N}^h_{2h}(k)-(n+2)} + \frac{1}{2}\phi_{h,\mathcal{N}^h_{2h}(k)-(n+1)} + \frac{1}{4}\phi_{h,\mathcal{N}^h_{2h}(k) - n }  \\
    &+\frac{1}{2}\phi_{h,\mathcal{N}^h_{2h}(k)-1} 
    + \phi_{h,\mathcal{N}^h_{2h}(k)}  + \frac{1}{2}\phi_{h,\mathcal{N}^h_{2h}(k)+1} \\
    &+ \frac{1}{4}\phi_{h,\mathcal{N}^h_{2h}(k) + n} + \frac{1}{2}\phi_{h,\mathcal{N}^h_{2h}(k)+(n+1)} 
       + \frac{1}{4}\phi_{h,\mathcal{N}^h_{2h}(k)+ (n+2)},
\end{align*}
for $k \in \{0,1,\dots,\left(\frac{n}{2} \right)^2+n\}$. 

\par The restriction operator is denoted by $\mathcal{I}^{h}_{2h} : \mathbb{R}^{(n+1)^2} \rightarrow \mathbb{R}^{(\frac{n}{2} + 1)^2}$, while the prolongation operator is $\mathcal{I}_{h}^{2h} : \mathbb{R}^{(\frac{n}{2} + 1)^2} \rightarrow \mathbb{R}^{(n+1)^2}$. The matrices corresponding to these operators satisfy $\mathcal{I}^{h}_{2h} = \left(\mathcal{I}_{h}^{2h}\right)^T$. For example, if $n = 2$, then  the restriction matrix $\mathcal{I}^{h}_{2h} \in \mathbb{R}^{4 \times 9}$ is given by  
\begin{align*}
    \mathcal{I}^h_{2h} = 
    \begin{pmatrix}
    1 & 0.5 & 0 & 0.5 & 0.25 & 0 & 0 & 0 & 0 \\
    0 & 0.5 & 1 & 0 & 0.25 & 0.5 & 0 & 0 & 0 \\
    0 & 0 & 0 & 0.5 & 0.25 & 0 & 1 & 0.5 & 0 \\
    0 & 0 & 0 & 0 & 0.25 & 0.5 & 0 & 0.5 & 1 
    \end{pmatrix}
    .
\end{align*}
As shown in Section~\ref{sect:restr_fem} for the one-dimensional case, the splitting and standard strategies in ghost-\ac{fem} are equivalent and the generalization from 1D to higher dimensions is straightforward.

\subsection{Global and local stabilization parameters for MG efficiency} \label{sec:2D_optimal_lambda}
\par The stabilization parameter $\lambda$ can be defined globally for all cut cells or locally for each cut cell. The estimate for the global lower bound of $C$ is obtained by solving the generalized eigenvalue problem 
\[ \boldsymbol{K} v = \Lambda \boldsymbol{M} v \quad \forall v \in \mathcal{V}_h,\]
where $\boldsymbol{K}_{ij} = \int_{\Gamma_D} (\boldsymbol{n}\cdot \boldsymbol{\nabla} \phi_i)(\boldsymbol{n}\cdot \boldsymbol{\nabla} \phi_j)~\mathrm{d}\boldsymbol{S}$ and $\boldsymbol{M}_{ij} = \int_{\Omega}  \boldsymbol{\nabla} \phi_i \cdot  \boldsymbol{\nabla} \phi_j~\mathrm{d}\boldsymbol{x}$ for $i,j \in \{0,1,\dots,n^2\}$. Here, $\phi_i, ~i \in \{0,1,\dots,n^2\}$, denote the shape functions of the \ac{fe} space $\mathcal{V}_h$. The lower bound for $C$ of \eqref{eq:lowC} and the stabilization parameter are, respectively, $C = \max \Lambda$ and $\lambda = \gamma C$, for some constant $\gamma>1$.
\par This is not the optimal choice for the penalty parameter because cut cells in different geometrical configurations may have varying sizes. While some cut cells are small and may require a large penalty parameter, others are well-conditioned and do not need an excessive value. Using a uniform penalty parameter may result in unnecessarily large values in well-conditioned cells, leading to degraded efficiency and loss of \ac{mg} convergence. A more effective approach would be adapting the penalty parameter locally based on the conditioning of each cut cell to balance stability and convergence. To this end, on each cut cell $K$ we solve the $4 \times 4$ generalized eigenvalue problem 
\[ \boldsymbol{K} v = \Lambda \boldsymbol{M} v \quad \forall v \in \mathcal{V}_h,\]
where $\boldsymbol{K}_{ij} = \int_{\Gamma_D \cap K} (\boldsymbol{n}\cdot \boldsymbol{\nabla} \phi_i)(\boldsymbol{n}\cdot \boldsymbol{\nabla} \phi_j)~\mathrm{d}\boldsymbol{S}$ and $\boldsymbol{M}_{ij} = \int_{K \cap \Omega}  \boldsymbol{\nabla} \phi_i \cdot  \boldsymbol{\nabla} \phi_j~\mathrm{d}\boldsymbol{x}$ for $i,j \in \{1,2,3,4\}$. Here $\phi_i, ~i \in \{1,2,3,4\}$, denote the non-zero shape functions on the cut cell $K$.  Then, $C(K) = \max \Lambda(K)$ and $\lambda(K) = \gamma C(K)$, for some constant $\gamma>1$.

\par
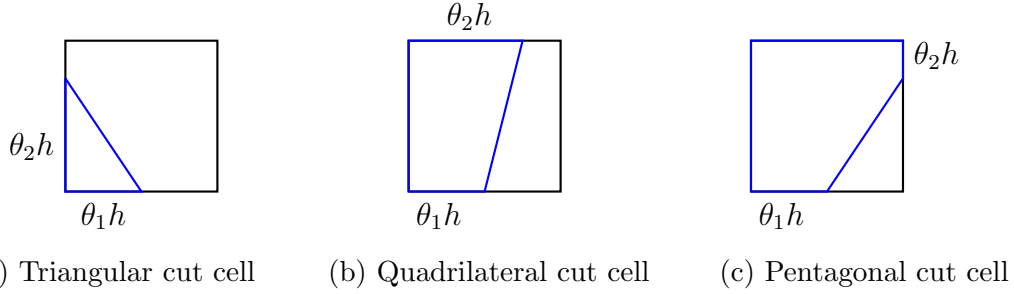
\begin{figure}
    \centering
    \begin{subfigure}[b]{0.3\textwidth}
    \centering
        \begin{tikzpicture}
    \draw[thick] (0,0) rectangle (2,2);


    \draw[thick,blue] (0,1.5) -- (1,0) -- (0,0)--(0,1.5);

    \node[below] at (0.5,0) {$\theta_1 h$};
    \node[left] at (0,0.6) {$\theta_2 h$};

\end{tikzpicture}
\label{fig:tri_cut_cell}
          \caption{Triangular cut cell}
    \end{subfigure}
    \begin{subfigure}[b]{0.3\textwidth}
    \centering
      \begin{tikzpicture}
    \draw[thick] (0,0) rectangle (2,2);


    \draw[thick,blue] (0,0) -- (1,0) -- (1.5,2)--(0,2)--(0,0);

    \node[below] at (0.4,0) {$\theta_1 h$};
    \node[above] at (0.8,2) {$\theta_2 h$};

\end{tikzpicture}
        \label{fig:quad_cut_cell}
        \caption{Quadrilateral cut cell}
    \end{subfigure}
    \begin{subfigure}[b]{0.3\textwidth}
    \centering
    \begin{tikzpicture}
    \draw[thick] (0,0) rectangle (2,2);


    \draw[thick,blue] (0,0) -- (1,0) -- (2,1.5)--(2,2)--(0,2)--(0,0);

    \node[below] at (0.4,0) {$\theta_1 h$};
    \node[right] at (2,1.8) {$\theta_2 h$};

\end{tikzpicture}
        \label{fig:pent_cut_cell}
        \caption{Pentagonal cut cell}
    \end{subfigure}
    \caption{Different configurations of cut regions for 2D geometries.}
    \label{fig:cut_geos}
\end{figure}

 For 2D geometries, a second-order accurate method approximates the boundary within each cut cell using line segments. As a result, cut cells typically take the shape of triangles, rectangles, or pentagons, depending on how the background Cartesian grid intersects with the embedded geometry. The geometrical configuration of the cut cells can be identified by evaluating the level set function $\psi$ at the vertices. Let $\{\boldsymbol{x}_i : i =1,2,3,4\}$ denote the vertices of the cut cell $K$. If $\psi(\boldsymbol{x}_i)< 0$, then it is an internal node, else it is an external node. If a cut cell has exactly one internal node, then its intersection with the domain results in a triangular region. Likewise, quadrilateral and pentagonal regions have exactly two and three internal nodes respectively. \cref{fig:cut_geos} illustrates various cut configurations along with their respective cut edge lengths.
\par The next task would be to determine the lengths of the cut sides. For a second-order accurate method, this can be achieved by using linear interpolation. For instance, let us consider a triangular cut region with $\boldsymbol{x}_1$ as the internal node. Then, the lengths $\theta_1 h$ and $\theta_2 h$ of the cut sides can be determined as 
\begin{align*}
    &h(1-\theta_1)\psi(\boldsymbol{x}_1) + h\theta_1\psi(\boldsymbol{x}_2) = 0, \qquad  h(1-\theta_2)\psi(\boldsymbol{x}_1) + h\theta_2\psi(\boldsymbol{x}_3) = 0 \\
    &\implies\theta_1 = \frac{\psi(\boldsymbol{x}_1)}{\psi(\boldsymbol{x}_1)-\psi(\boldsymbol{x}_2)} , 
    \qquad \theta_2 = \frac{\psi(\boldsymbol{x}_1)}{\psi(\boldsymbol{x}_1)-\psi(\boldsymbol{x}_3)}.
\end{align*}
The lengths of the sides quadrilateral and pentagonal cut cells can be determined using a similar approach.

\par The general expression for the constant $C=C(K)$ was derived using Mathematica for triangular and quadrilateral regions, assuming the lengths of the cut cell sides to be $\theta_1 h$ and $\theta_2h$. Pentagonal cells are well-posed and do not introduce ill-conditioning problems, as they cannot be arbitrarily small, so it is sufficient to set 
$C$ equal to the value used for a triangular cell that occupies half of the background cell. The lower bound for $C$ for triangular cut regions is 
\begin{align}\label{eq:c_triangle}
    C = \frac{\sqrt{\frac{\theta_1^2 + \theta_2^2 }{\theta_1^2} } ~\Big(3  \theta_1^4 + 3  \theta_2^4 +\sqrt{3(\theta_1^8 +10  \theta_1^4 \theta_2^4+\theta_2^8)} \Big)}{h~ \theta_2  (\theta_1^2 + \theta_2^2)^2 }.
\end{align}
Setting $\theta_1 = \theta_2 = 1$ in $\cref{eq:c_triangle}$, yields 
\begin{align}\label{eq:c_pentagon}
    C = \frac{3 \sqrt{2}}{h}.
\end{align}
This value is chosen as lower bound for $C$ for pentagonal cut regions. For quadrilateral cut regions, the lower bound for $C$ was derived as 
\begin{align*}
    C = &3  h  \sqrt{1 + \theta_1^2 - 2  \theta_1 \theta_2 + \theta_2^2} ~ (\theta_1^6 + 2  \theta_1^3  \theta_2 + \theta_2^2 - \theta_1^4  \theta_2^2 + \theta_2^6 + 2  \theta_1  (\theta_2 + \theta_2^3)) \\
    &- \theta_1^2  (-1 + 2 \theta_2^2 + \theta_2^4) + (3 h^2 (\theta_1^{14} + 2 * \theta_1^{13}  \theta_2 + \theta_1^{12}  (1 - 5  \theta_2^2) - 4  \theta_1^{11}  \theta_2   \\
    &(-4 + 7  \theta_2^2) - 2  \theta_1^9  \theta_2  (-6 + 20  \theta_2^2 + 17  \theta_2^4) + \theta_1^{10}  (10 - 22  \theta_2^2 + 73  \theta_2^4) + \theta_1^8   \\
    &(10 - 38  \theta_2^2 + 63  \theta_2^4 - 69  \theta_2^6))+ 8  \theta_1^7  \theta_2  (4 - 4  \theta_2^2 + 3  \theta_2^4 + 15  \theta_2^6) + \theta_1^6    \\
    & (1 - 32\theta_2^2 + 136  \theta_2^4 - 84  \theta_2^6 - 69  \theta_2^8) + 2  \theta_1  \theta_2^5  (5 + 16  \theta_2^2 + 6  \theta_2^4 + 8  \theta_2^6 + \theta_2^8) \\
    &- 4 \theta_1^3  \theta_2^3  (5 + 14  \theta_2^2 + 8  \theta_2^4 + 10  \theta_2^6 + 7  \theta_2^8) - 2  \theta_1^5  \theta_2  (-5 + 28  \theta_2^2 + 88  \theta_2^4\\
    &- 12  \theta_2^6 + 17  \theta_2^8) 
    + \theta_2^4  (1 + \theta_2^2 + 10  \theta_2^4 + 10  \theta_2^6 + \theta_2^8 + \theta_2^{10}) - \theta_1^2 \theta_2^2(2 + \theta_2^2 \\
    &+ 32  \theta_2^4 + 38  \theta_2^6 + 22  \theta_2^8 + 5  \theta_2^{10})  + \theta_1^4  (1 - \theta_2^2 + 104  \theta_2^4 + 136  \theta_2^6 + 63  \theta_2^8 \\
    &+ 73  \theta_2^{10}))^{1/2} 
    (h^2  (\theta_1^7 + \theta_1^6  \theta_2 + \theta_1^5  (2 - 3  \theta_2^2) + \theta_1^3  (-1 + \theta_2^2)^2 + \theta_2^3  (1 + \theta_2^2)^2 \\
    & + \theta_1^4  \theta_2  (6 + \theta_2^2) + \theta_1  \theta_2^2  (5 + 6  \theta_2^2 + \theta_2^4) + \theta_1^2  (5 \theta_2 - 2  \theta_2^3 - 3  \theta_2^5)))^{-1}.
\end{align*}

\subsection{Numerical tests} \label{sec:2D_tests}



\par
\begin{figure}
    \centering
    \begin{subfigure}[t]{.3\textwidth}
   \centering
   \includegraphics[width=\textwidth]{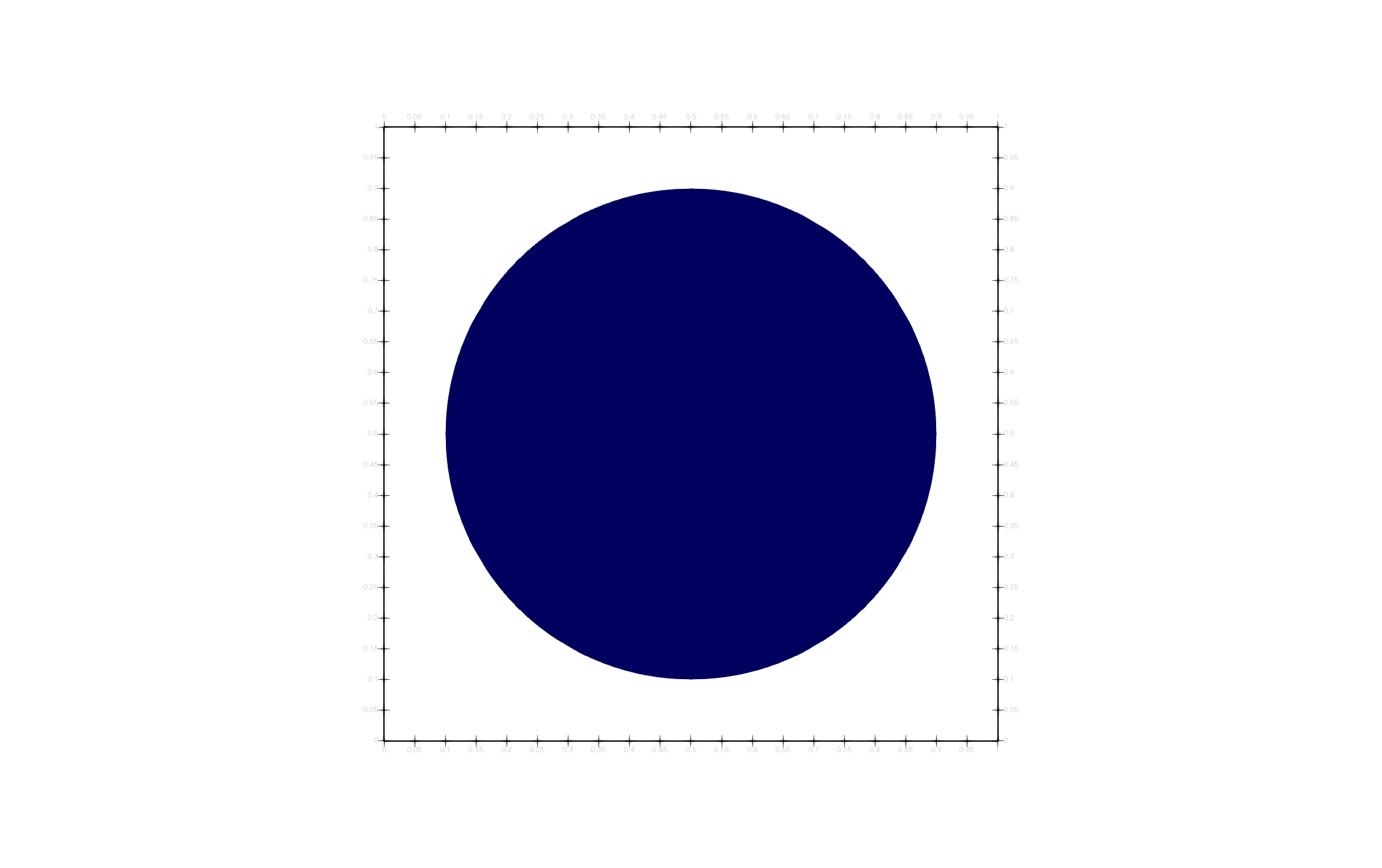}
   \caption{Disk}
   \label{fig:circle}
 \end{subfigure}%
 \begin{subfigure}[t]{.3\textwidth}
   \centering
   \includegraphics[width=\textwidth]{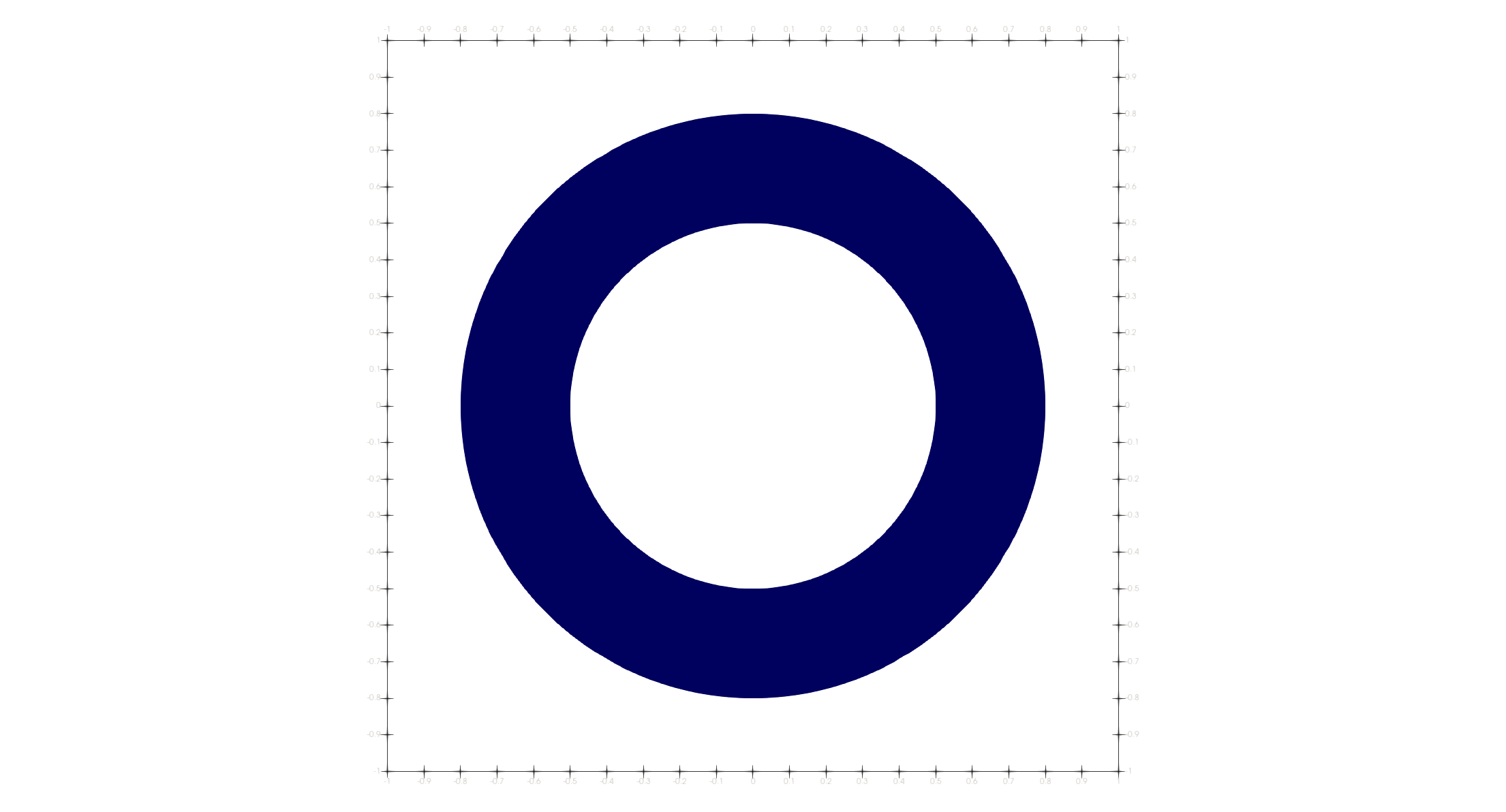}
   \caption{Annulus}
   \label{fig:annulus}
 \end{subfigure}%
 \begin{subfigure}[t]{.3\textwidth}
   \centering
   \includegraphics[width=\textwidth]{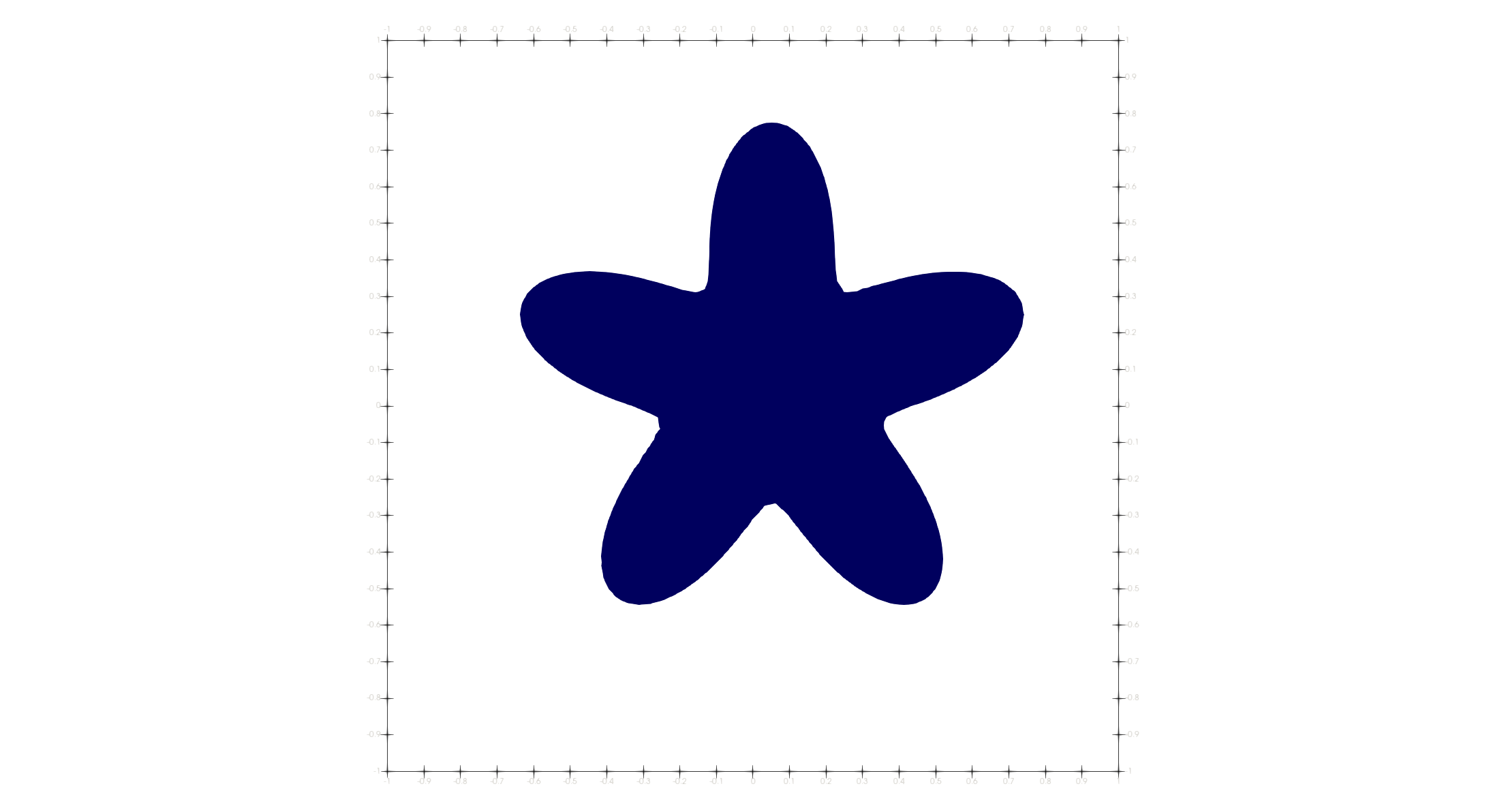}
   \caption{Flower}
   \label{fig:flower}
 \end{subfigure}%
 \vskip\baselineskip

 \begin{subfigure}[t]{.3\textwidth}
   \centering
   \includegraphics[width=\textwidth]{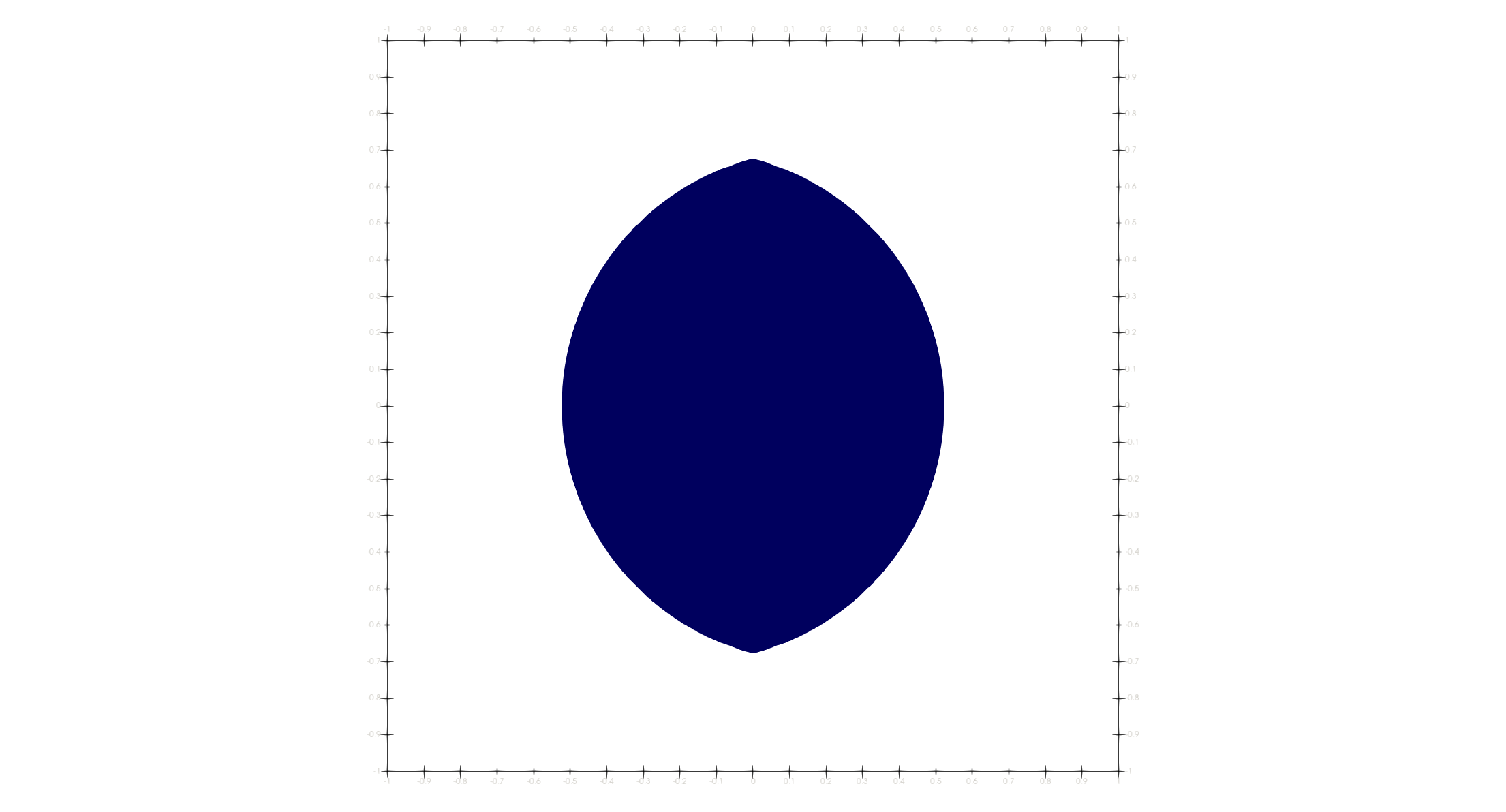}
   \caption{Leaf}
   \label{fig:leaf}
 \end{subfigure}%
 \begin{subfigure}[t]{.3\textwidth}
   \centering
   \includegraphics[width=\textwidth]{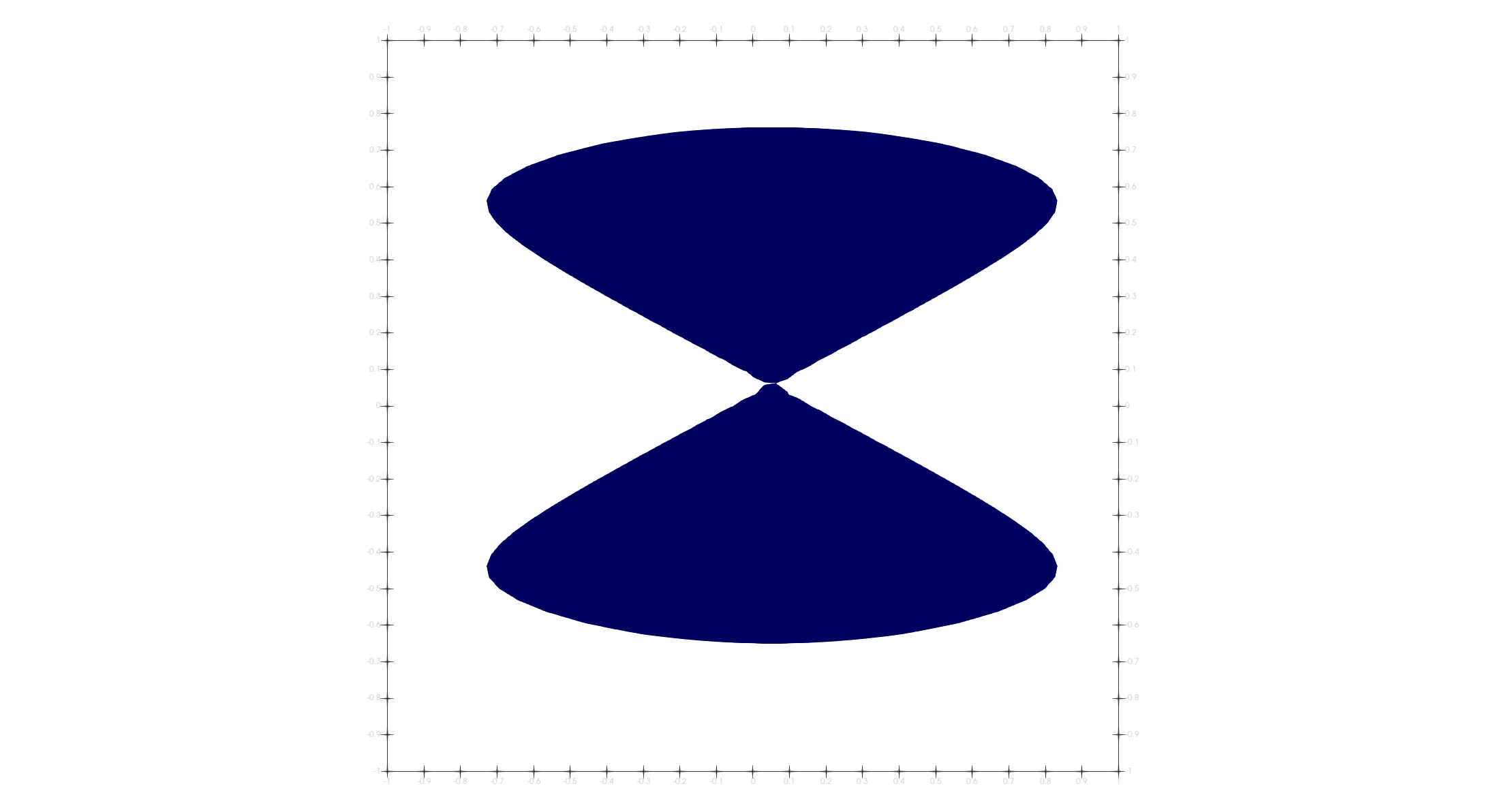}
   \caption{Hourglass}
   \label{fig:hourglass}
 \end{subfigure}%
    \caption{2D domains used for the numerical experiments.}
    \label{fig:2D_geos}
\end{figure}

 In this section, we assess the performance of the \ac{mg} on various 2D geometries, shown in \cref{fig:2D_geos}. The artificial domain is $\Omega_{art} = [0,1]^2$. As for the 1D case, we solve the homogeneous Poisson equation with the exact solution $u(x,y) = 0$, since for a nonzero solution the computed residual would remain constant once machine precision is reached~\cite{Briggs2000}. Each iteration of the scheme consists of two pre-smoothing steps and one post-smoothing step, that is, $\nu_1= 2$ and $\nu_2 =1$. The snapping threshold is set as $h^\alpha$ with $\alpha = 1.75$. Both the \ac{tgcs} and the W-cycle are executed for 30 iterations, and the convergence factor is determined using 
\[\rho = \frac{1}{10}\sum_{m = 21}^{30} \rho^{(m)} .\] 

\par All numerical experiments are performed using the open-source software \texttt{Gridap}. \texttt{Gridap} is a Julia-based \ac{fe} framework designed for approximating the solution of \acp{pde} \cite{Badia2020}. \texttt{Gridap} leverages the Julia \ac{jit} compiler and has a user-interface that closely resembles that whiteboard-style notation of the weak formulation \cite{Verdugo2022}.

\subsubsection*{Rectangular domain (with an arbitrary side)}
\par
\begin{figure}
    \centering
    \begin{subfigure}[t]{.5\textwidth}
   \centering
   \includegraphics[width=\textwidth]{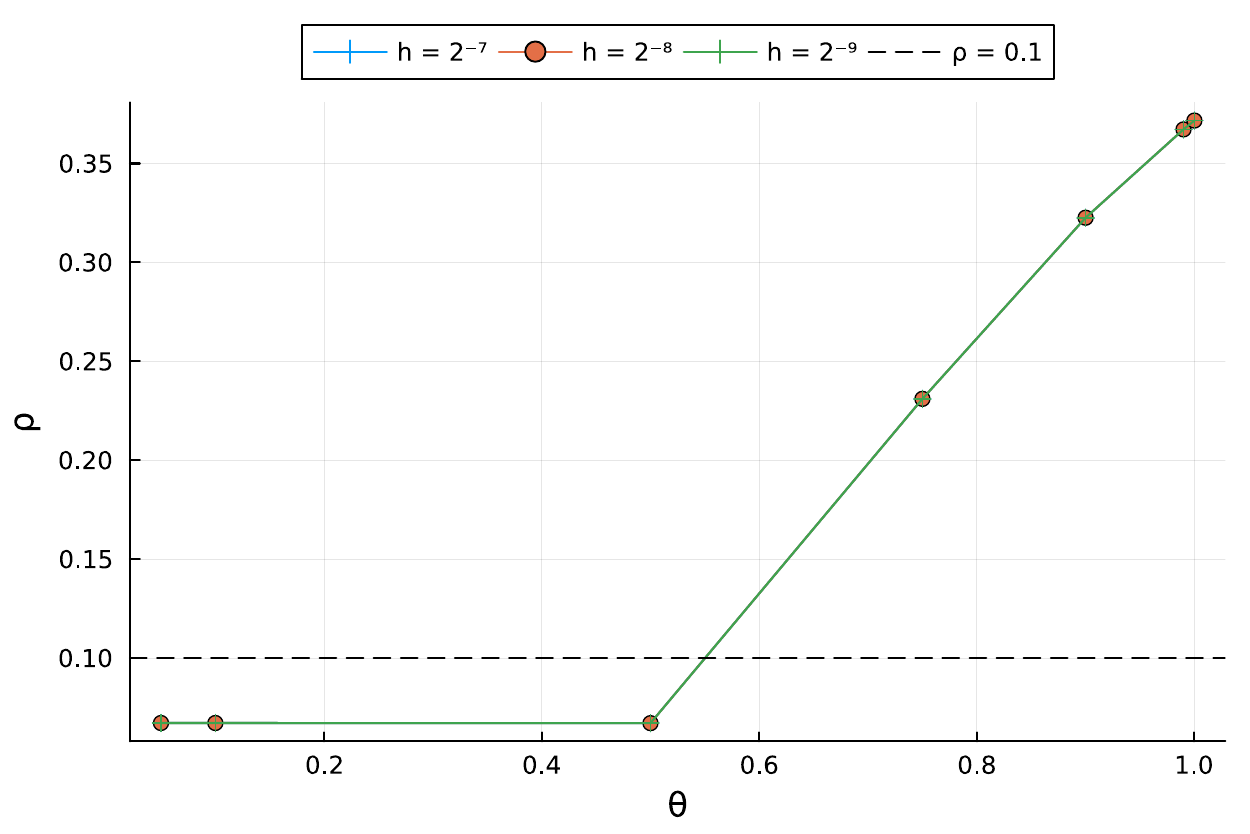}
   \caption{$\gamma = 1.1$}
   \label{fig:2D_vertical_line_gamma_1.1}
 \end{subfigure}%
 \begin{subfigure}[t]{.5\textwidth}
   \centering
   \includegraphics[width=\textwidth]{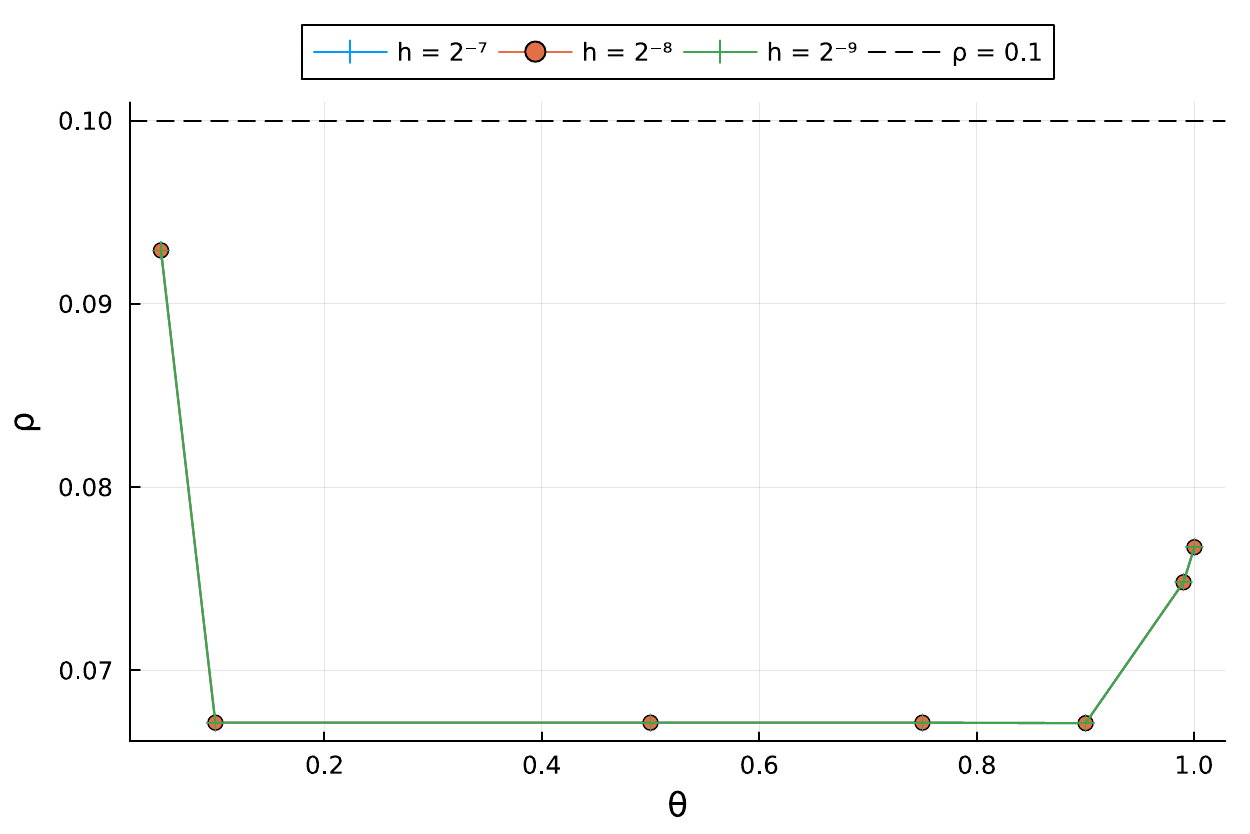}
   \caption{$\gamma = 1.3$}
   \label{fig:2D_vertical_line_gamma_1.3}
 \end{subfigure}%
  \vskip\baselineskip
  \begin{subfigure}[t]{.5\textwidth}
   \centering
   \includegraphics[width=\textwidth]{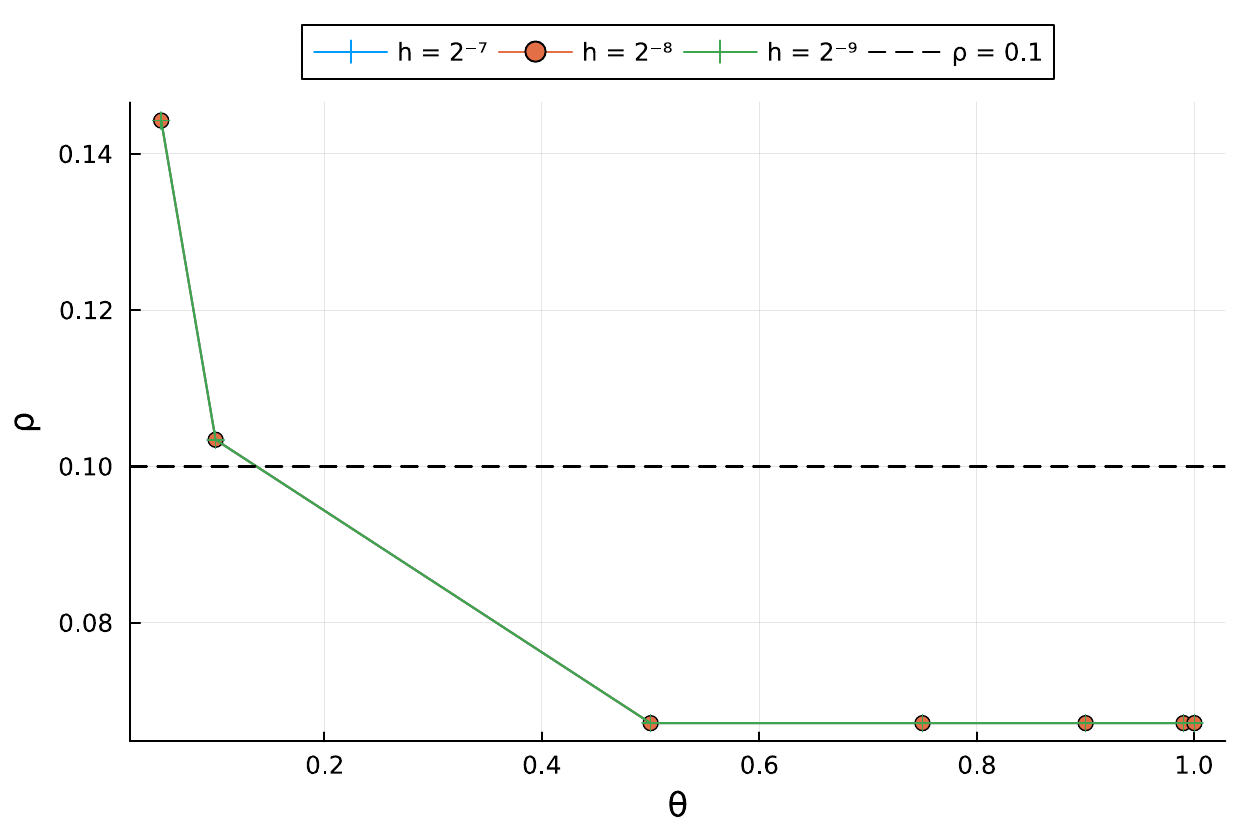}
   \caption{$\gamma = 1.5$}
   \label{fig:2D_vertical_line_gamma_1.5}
 \end{subfigure}%
 \begin{subfigure}[t]{.5\textwidth}
   \centering
   \includegraphics[width=\textwidth]{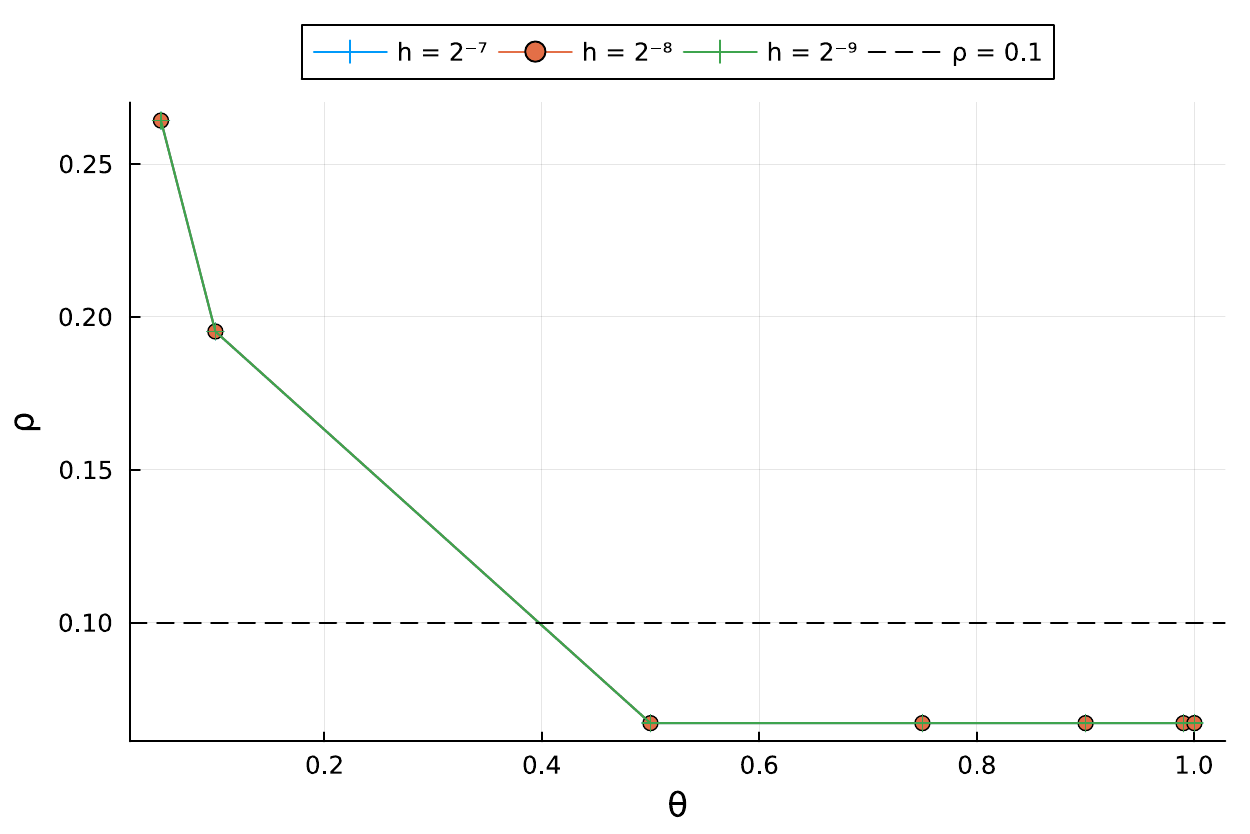}
   \caption{$\gamma = 2$}
   \label{fig:2D_vertical_line_gamma_2}
 \end{subfigure}%
 \vskip\baselineskip
  \begin{subfigure}[t]{.5\textwidth}
   \centering
   \includegraphics[width=\textwidth]{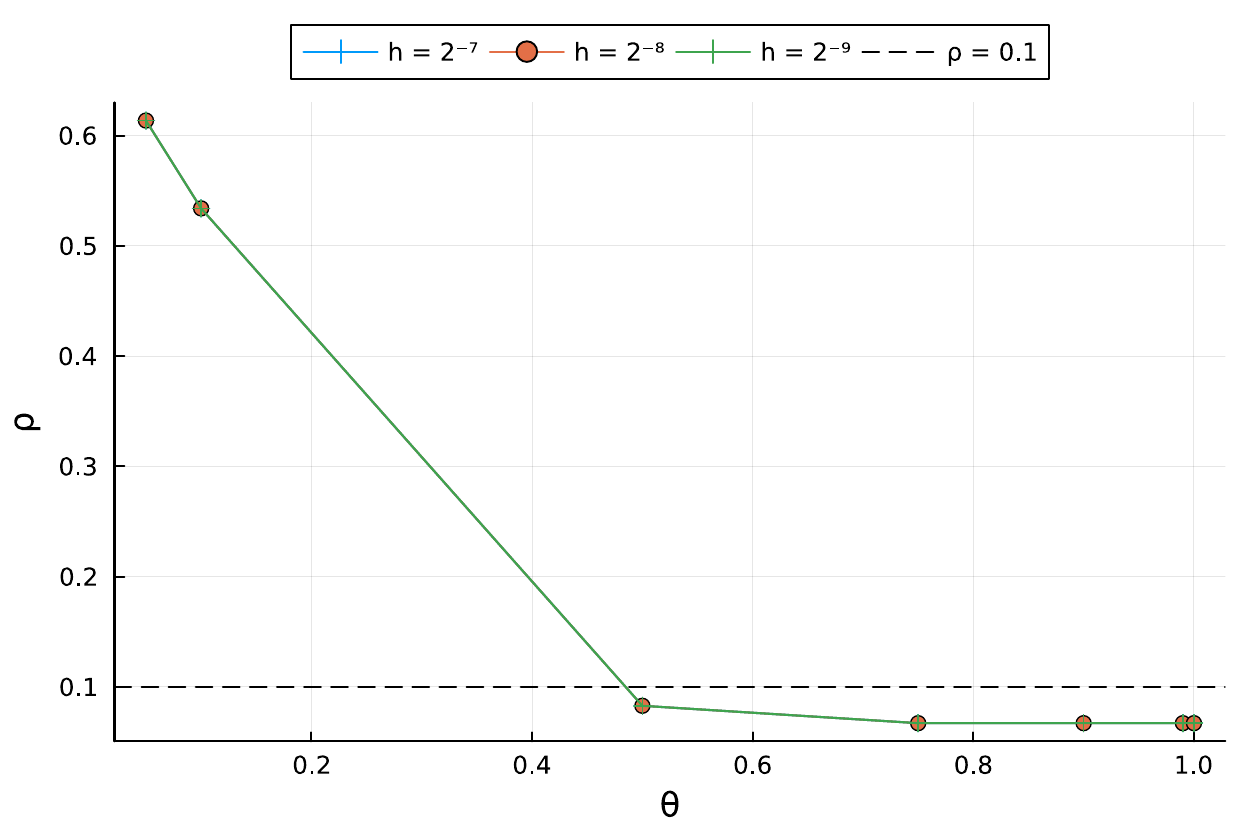}
   \caption{$\gamma = 5$}
   \label{fig:2D_vertical_line_gamma_5}
 \end{subfigure}%
 \begin{subfigure}[t]{.5\textwidth}
   \centering
   \includegraphics[width=\textwidth]{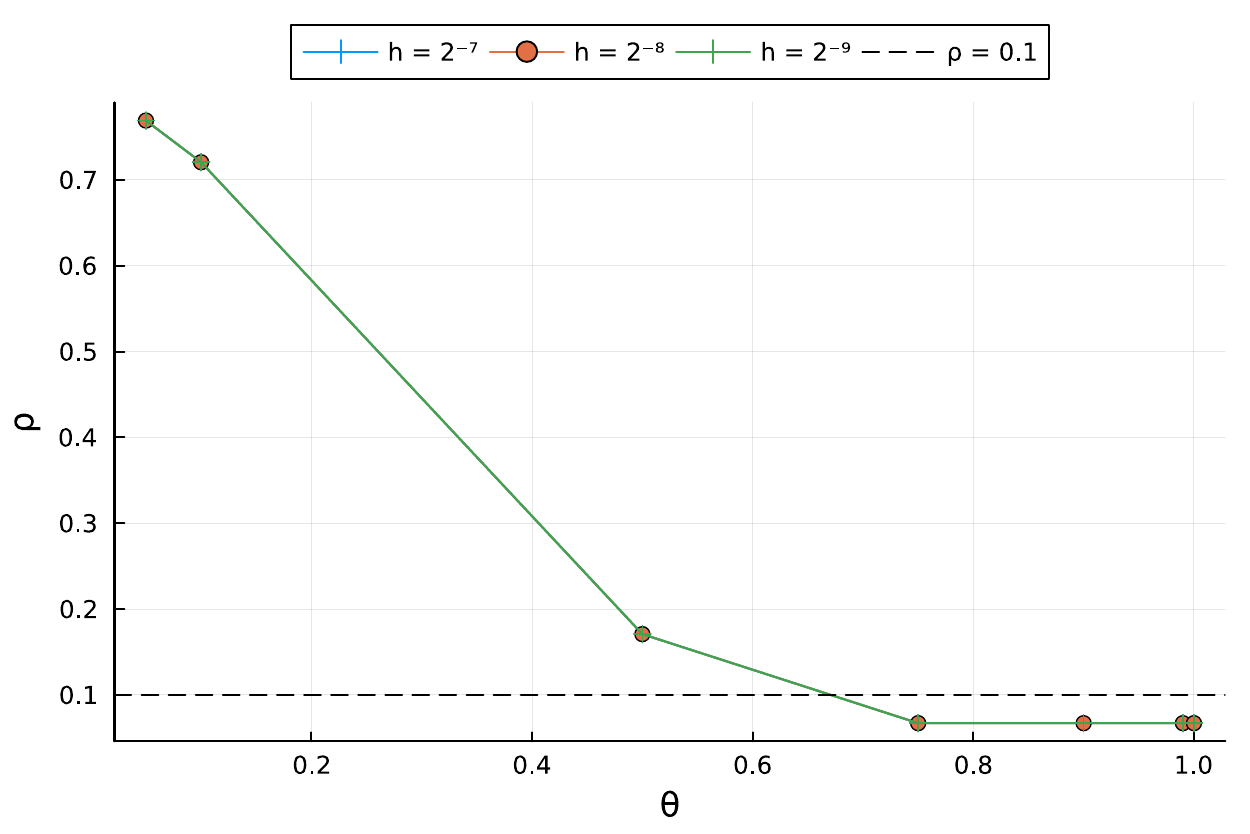}
   \caption{$\gamma = 10$}
   \label{fig:2D_vertical_line_gamma_10}
 \end{subfigure}%
    \caption{Plot of convergence factor of a \ac{tgcs} $\rho$ against $\theta$ for various mesh sizes and choices of $\gamma$ on a rectangular (with an arbitrary side).}
    \label{fig:2D_vertical_line}
\end{figure}

 In the first test, we consider a rectangular domain $\Omega$ formed by intersecting $\Omega_{art}=[0,1]^2$ with a vertical line, defined by the level-set function $\psi(x,y) = x - 1 + (1 - \theta)h$, where $\theta \in [0,1]$. Then $\Omega = \left\{ (x,y) \in \Omega_{art} \colon \psi(x,y)<0 \right\}$  This results in cut-cells that are quadrilaterals of equal size, with the length of the cut sides precisely equal to $\theta h$. The mesh size is $h = 2^{-m}, ~ m \in \{7,8,9\}$. The Dirichlet boundary condition is weakly enforced on $\Gamma_{Dw} = \{(x,y) \in \Omega_{art}:\psi(x,y) = 0\} $ and  strongly enforced on $\Gamma_{Ds} = \partial \Omega \backslash \Gamma_{Dw} = \{(x,0) : 0\le x \le 1 - (1-\theta)h\} \cup   \{(x,1) : 0\le x \le 1 - (1-\theta)h\} \cup  \{(y,0) : 0\le y \le 1 \}$. The stabilization parameter for quadrilaterals with cut lengths $\theta_1h = \theta_2 h = \theta h$ is given by $\lambda = \frac{\gamma}{\theta h}$. We choose $\gamma \in \{1.1,1.3,1.5,2,5,10 \}$ for our computations. The convergence factor is evaluated for different level-set positions and mesh sizes. The plot of the convergence factor $\rho$ against level-set position $\theta$ is shown in \cref{fig:2D_vertical_line} for various mesh sizes and choices of $\gamma$. We observe that choosing $\gamma = 1.3$ results in optimal \ac{mg} convergence. Additionally, as $\gamma$ increases, the efficiency of the \ac{mg} scheme decreases.

\subsubsection*{Circular domain} 
\par We perform several tests on a circular domain defined by the level set $\psi(x,y) = (x-x_c)^2 + (y-y_c)^2 - r^2$, centered at $(x_c,y_c) = (0.5,0.5)$ with a radius $r = 0.4$. We impose purely Dirichlet boundary conditions on $\Gamma_D = \{(x,y):\psi(x,y) = 0\} $. The artificial domain is $\Omega_{art} = [0,1]^2$ and the mesh size is $h = 2^{-m}, ~ m \in \{6,7,8,9\}$.

\par
\begin{figure}
    \centering
    \begin{subfigure}[t]{.5\textwidth}
   \centering
   \includegraphics[width=\textwidth]{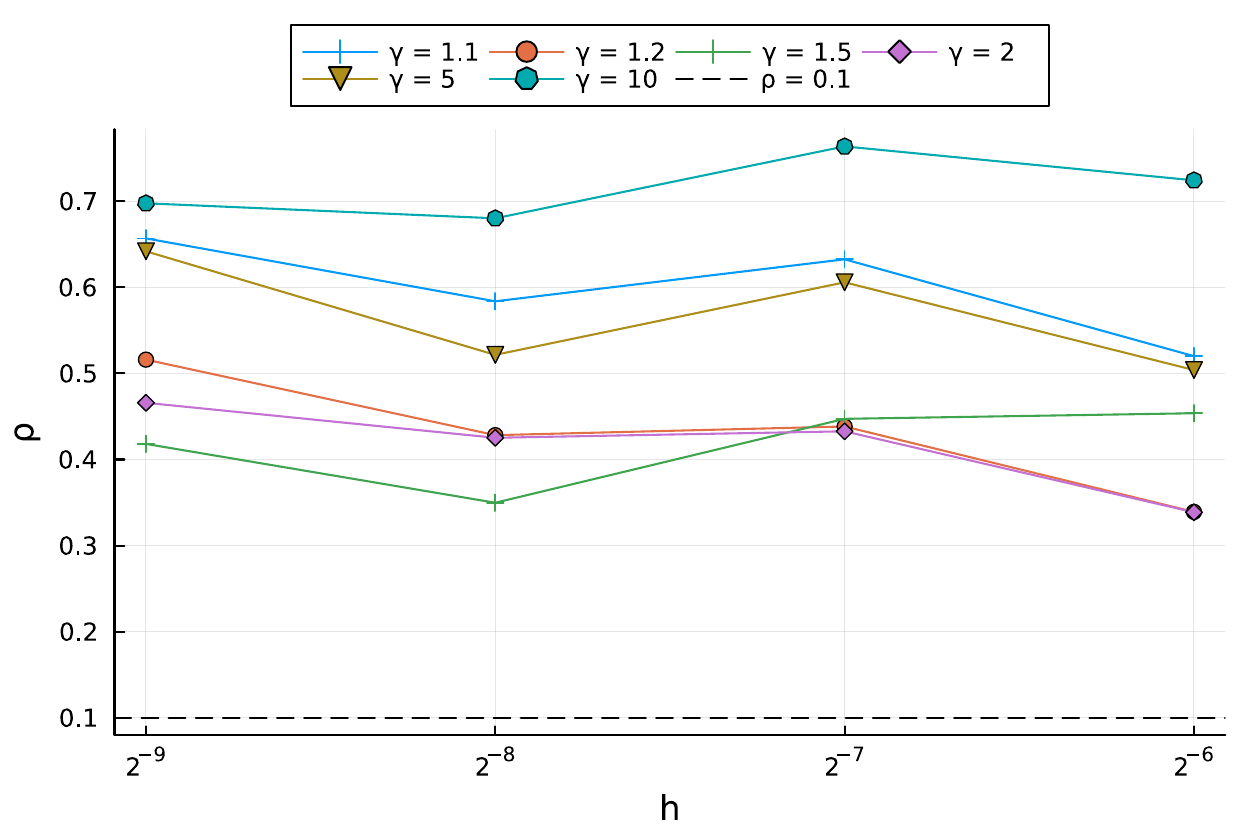}
   \caption{using a local stabilization parameter}
   \label{fig:2D_two_grid_circ_local_no_bdry_iter}
 \end{subfigure}%
 \begin{subfigure}[t]{.5\textwidth}
   \centering
   \includegraphics[width=\textwidth]{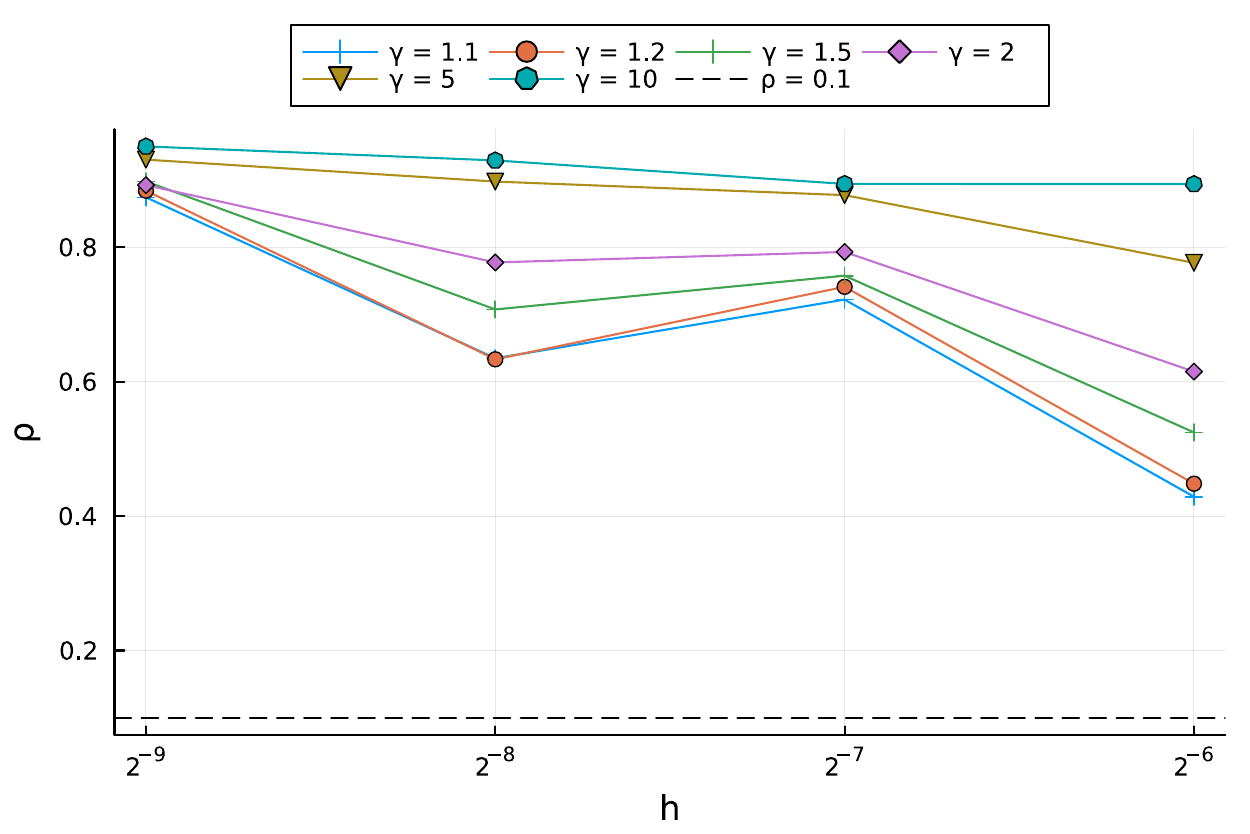}
   \caption{using a global stabilization parameter}
   \label{fig:2D_two_grid_circ_global_no_bdry_iter}
 \end{subfigure}%
    \caption{Plots of convergence factor of a \ac{tgcs} $\rho$ against mesh size $h$ on a circular domain using local and global stabilization parameters for different values of $\gamma$.}
    \label{fig:2D_circle_two_grid}
\end{figure}

 In the first test, we compute the convergence factor for a \ac{tgcs} using a local stabilization parameter. For each cut cell $K$ we use $\lambda(K) = \gamma \cdot C(K)$, where $C(K)$ is computed as described in \cref{sec:2D_optimal_lambda} and $\gamma \in  \{ 1.1,1.2,1.5,2,5,10\}$. The plots of the convergence factor against the mesh size are shown for different choices of $\gamma$ in \cref{fig:2D_two_grid_circ_local_no_bdry_iter}.

\par In the second test, we compute the convergence factor for a \ac{tgcs} using a global stabilization parameter, $\lambda = \gamma \cdot C$, where $\gamma \in  \{ 1.1,1.2,1.5,2,5,10\}$. The plots of the convergence factor against mesh size are shown in \cref{fig:2D_two_grid_circ_global_no_bdry_iter}. 

\par We observe that the \ac{mg} convergence improves when using a local stabilization parameter compared to the global one. This is because the cut cells can vary in size, and smaller cut cells require higher stabilization than larger ones. With a global parameter, we apply a larger value to all cut cells, even though they may not all require it. As a result, the convergence of the two-grid scheme is significantly affected. However, in both cases, the convergence is not optimal, that is, $\rho \gg 0.1$.
The reduced efficiency is attributed to the presence of cut cells, as the residual along the boundary is not sufficiently smoothed, causing boundary effects to propagate from the cut cells into the interior regions. To overcome the issue and improve efficiency, additional smoothing iterations are applied exclusively to the cut cells. The extra computational cost becomes negligible as $h \rightarrow 0$, since the ratio of boundary cells to internal cells converges towards zero. Moreover, this approach is computationally more efficient than applying extra smoothing steps across the entire domain, as performed, for instance, in~\cite{Kothari2021}.

\begin{figure}
    \centering
    \begin{subfigure}[t]{.5\textwidth}
   \centering
   \includegraphics[width=\textwidth]{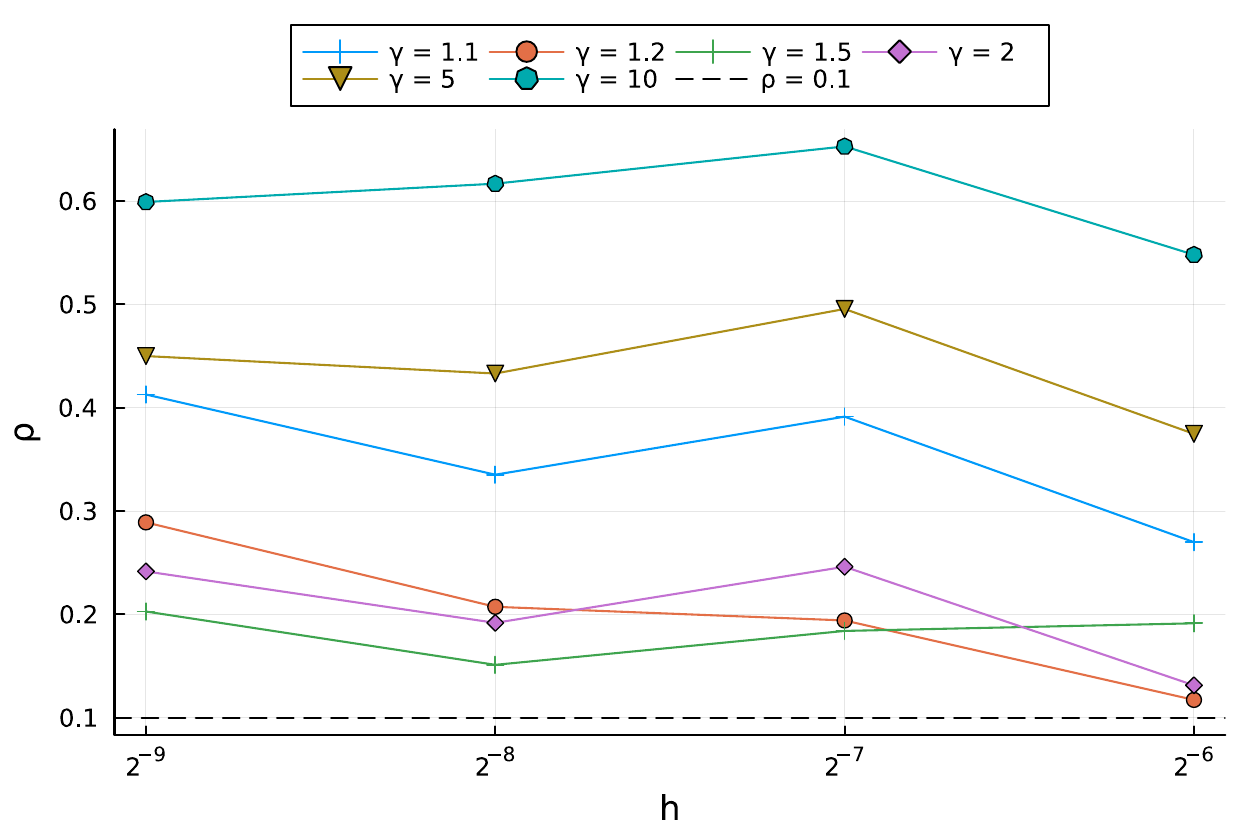}
   \caption{$\eta = 1$}
   \label{fig:2D_two_grid_circ_local_iter_1}
 \end{subfigure}%
 \begin{subfigure}[t]{.5\textwidth}
   \centering
   \includegraphics[width=\textwidth]{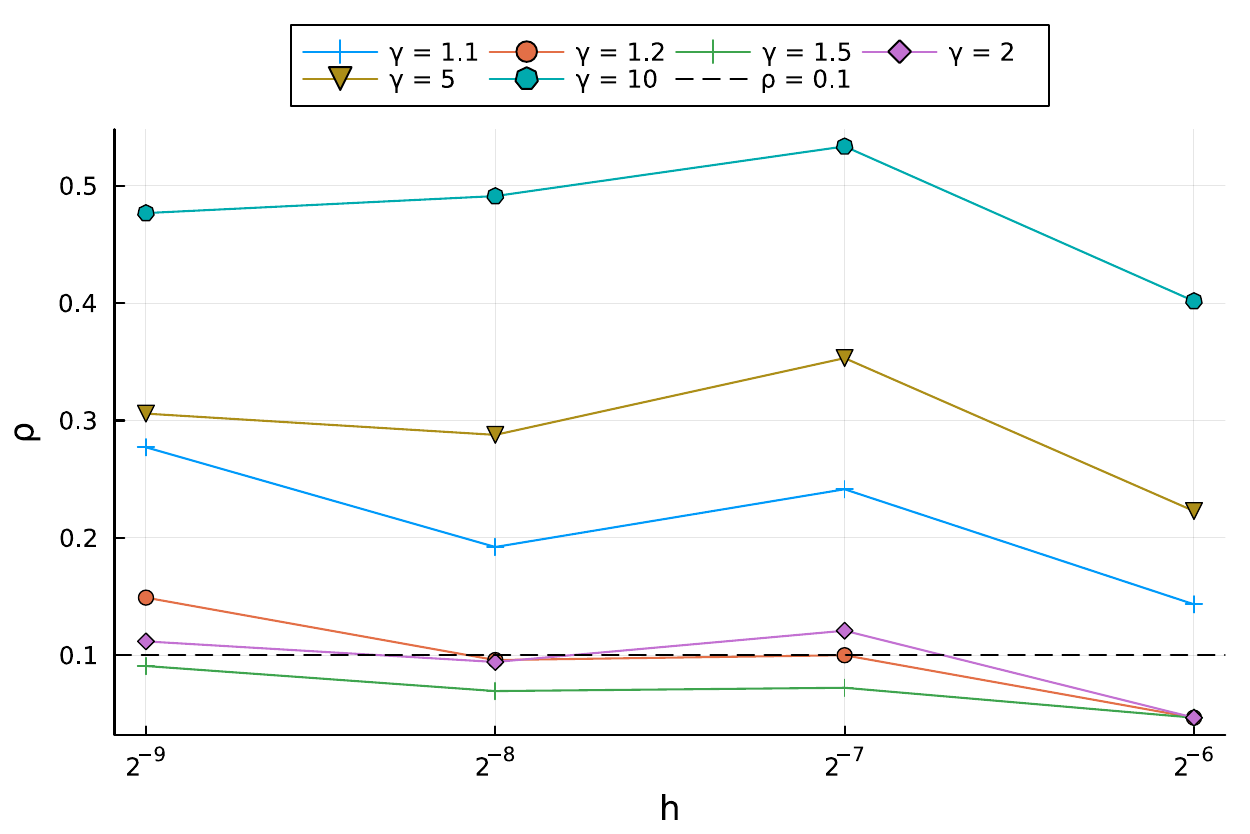}
   \caption{$\eta = 2$}
   \label{fig:2D_two_grid_circ_local_iter_2}
 \end{subfigure}%
  \vskip\baselineskip
    \begin{subfigure}[t]{.5\textwidth}
   \centering
   \includegraphics[width=\textwidth]{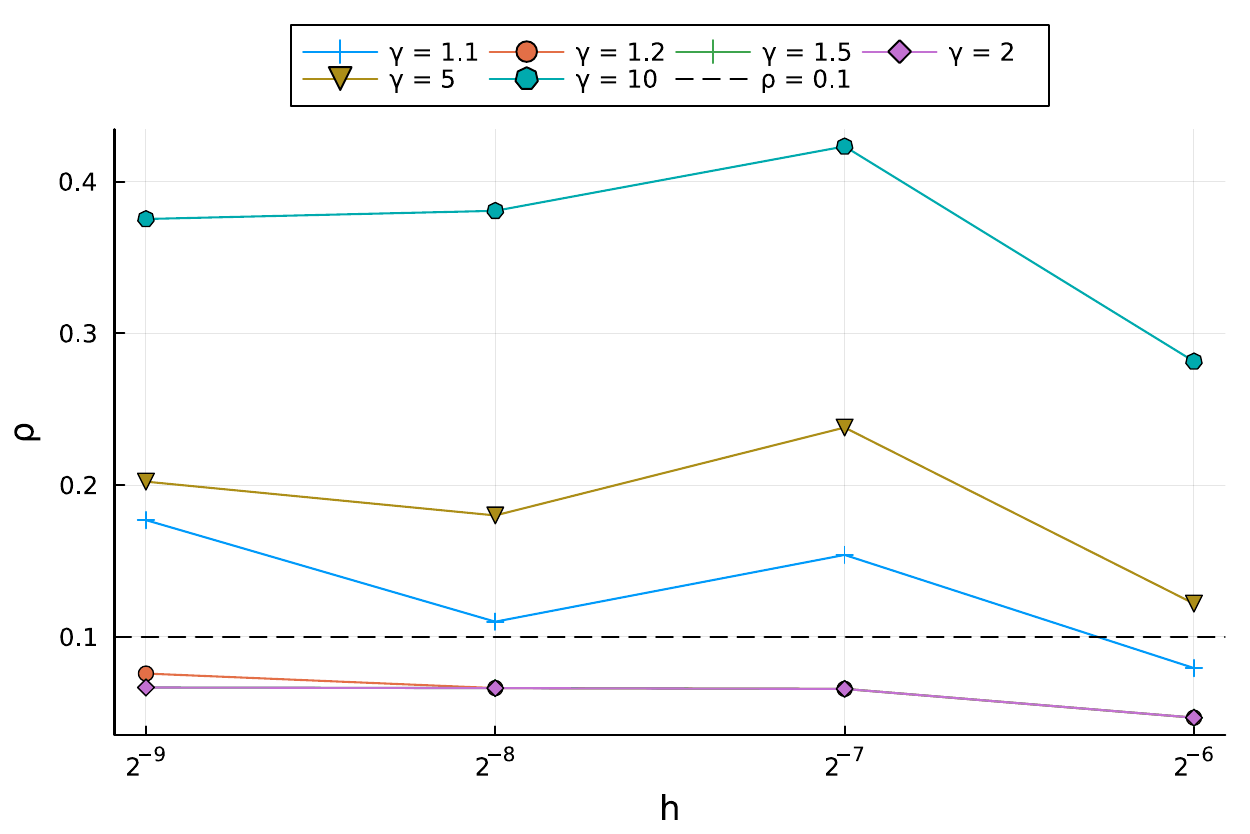}
   \caption{$\eta = 3$}
   \label{fig:2D_two_grid_circ_local_iter_3}
 \end{subfigure}%
  \begin{subfigure}[t]{.5\textwidth}
   \centering
   \includegraphics[width=\textwidth]{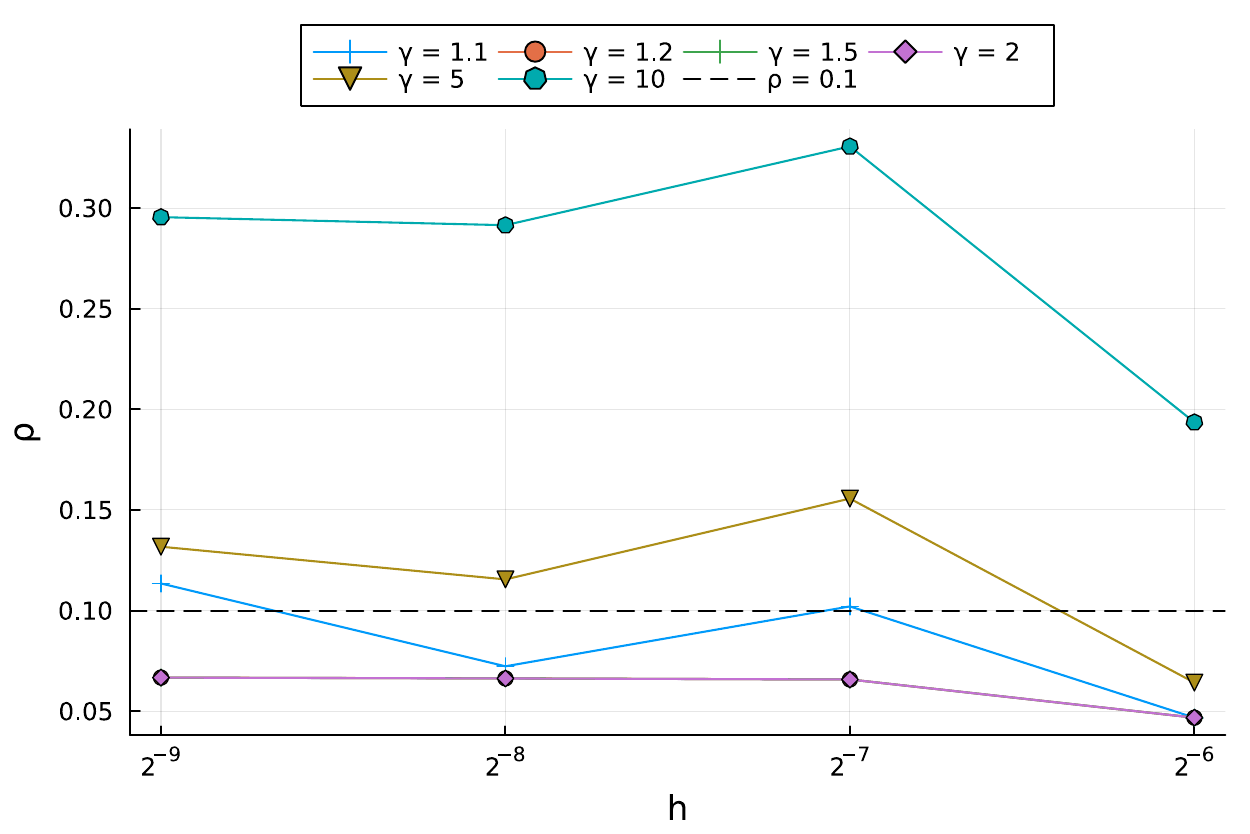}
   \caption{$\eta = 4$}
   \label{fig:2D_two_grid_circ_local_iter_4}
 \end{subfigure}%
 \caption{Plots of convergence factor of a \ac{tgcs} $\rho$ against mesh size $h$ on a circular domain using a local stabilization parameter, with $\eta$ additional iterations of smoother solely on cut cells.}
    \label{fig:2D_two_grid_circ_local}
\end{figure}

\begin{figure}
    \centering
    \begin{subfigure}[t]{.5\textwidth}
   \centering
   \includegraphics[width=\textwidth]{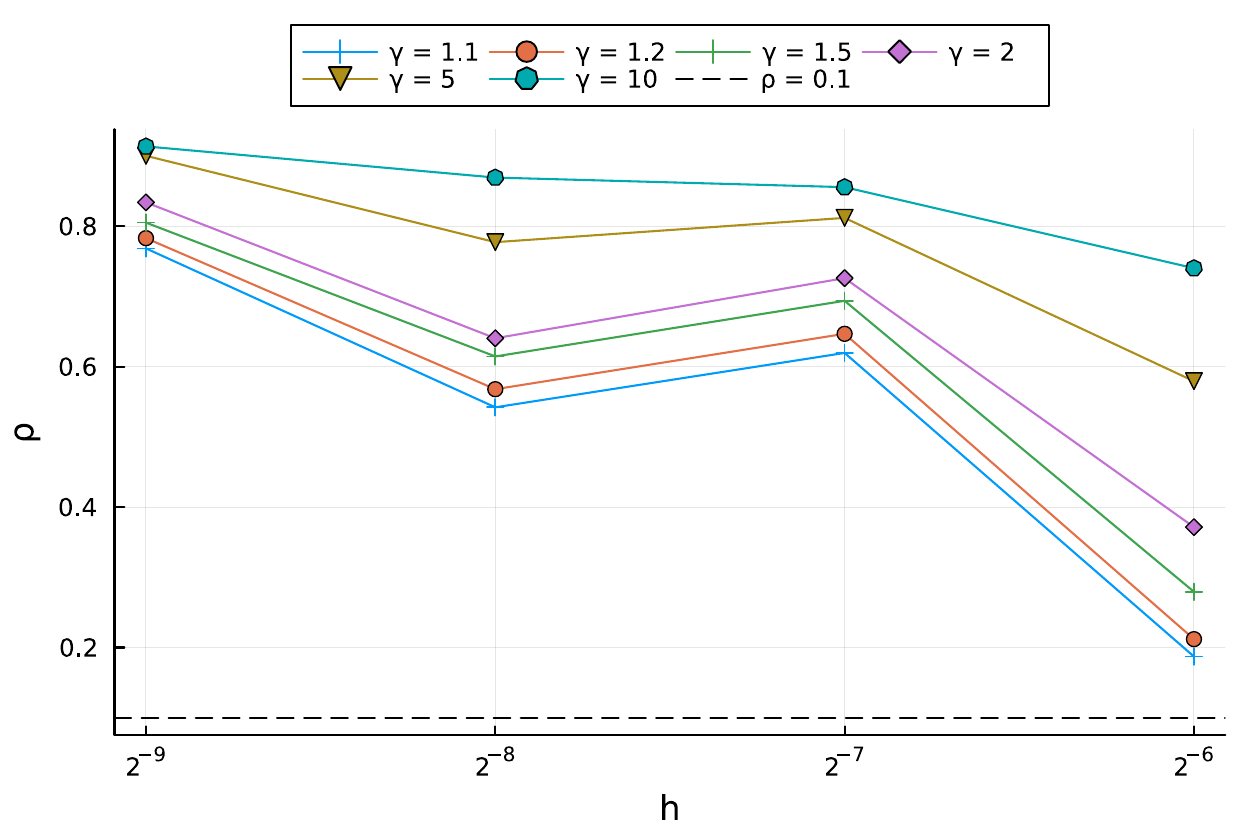}
   \caption{$\eta = 1$}
   \label{fig:2D_two_grid_circ_global_iter_1}
 \end{subfigure}%
 \begin{subfigure}[t]{.5\textwidth}
   \centering
   \includegraphics[width=\textwidth]{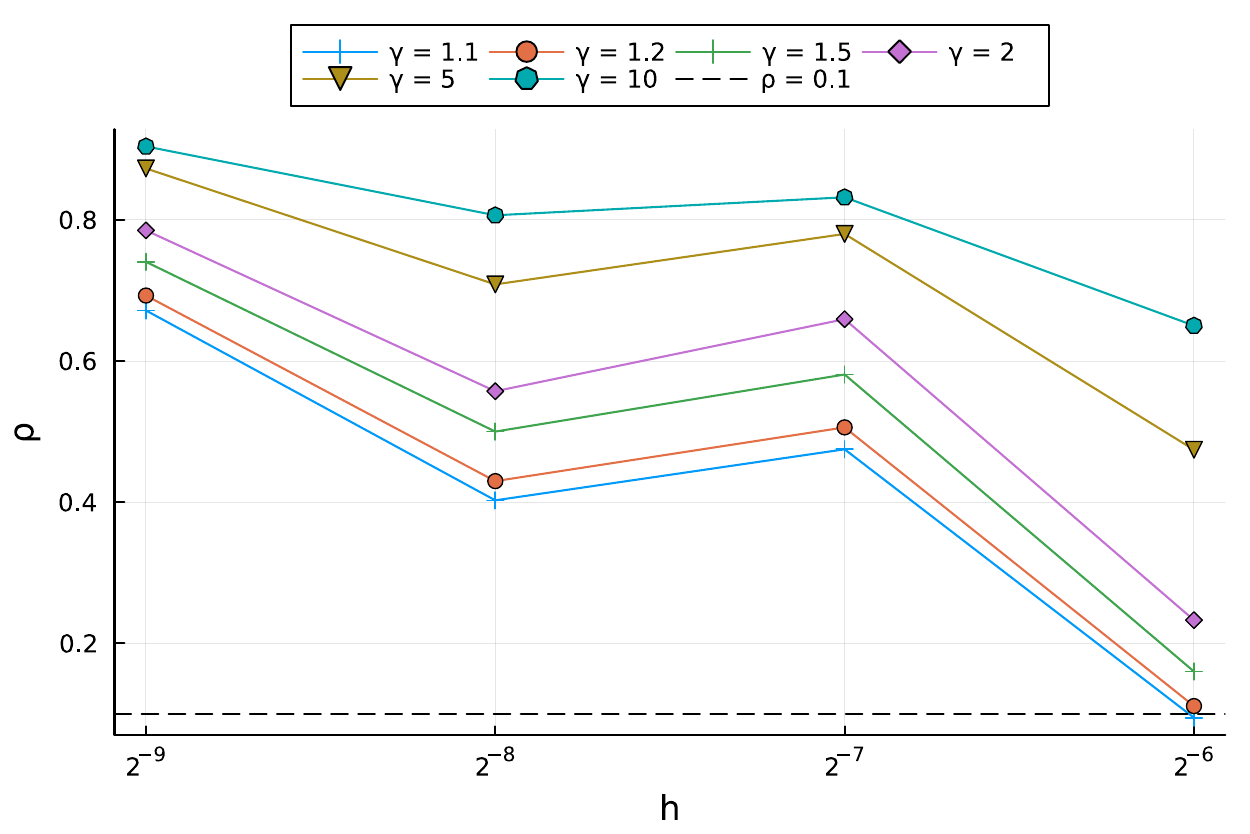}
   \caption{$\eta = 2$}
   \label{fig:2D_two_grid_circ_global_iter_2}
 \end{subfigure}%
   \vskip\baselineskip
\begin{subfigure}[t]{.5\textwidth}
   \centering
   \includegraphics[width=\textwidth]{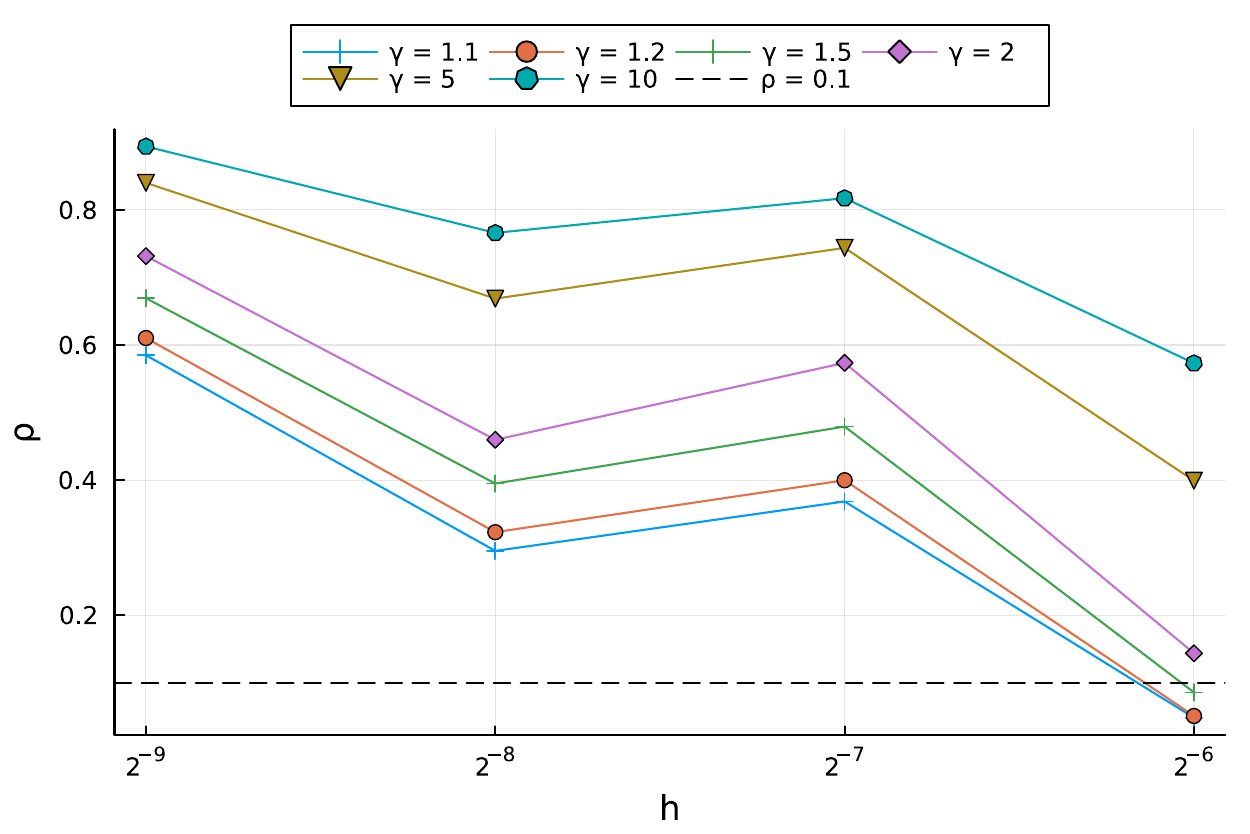}
   \caption{$\eta = 3$}
   \label{fig:2D_two_grid_circ_global_iter_3}
 \end{subfigure}%
 \begin{subfigure}[t]{.5\textwidth}
   \centering
   \includegraphics[width=\textwidth]{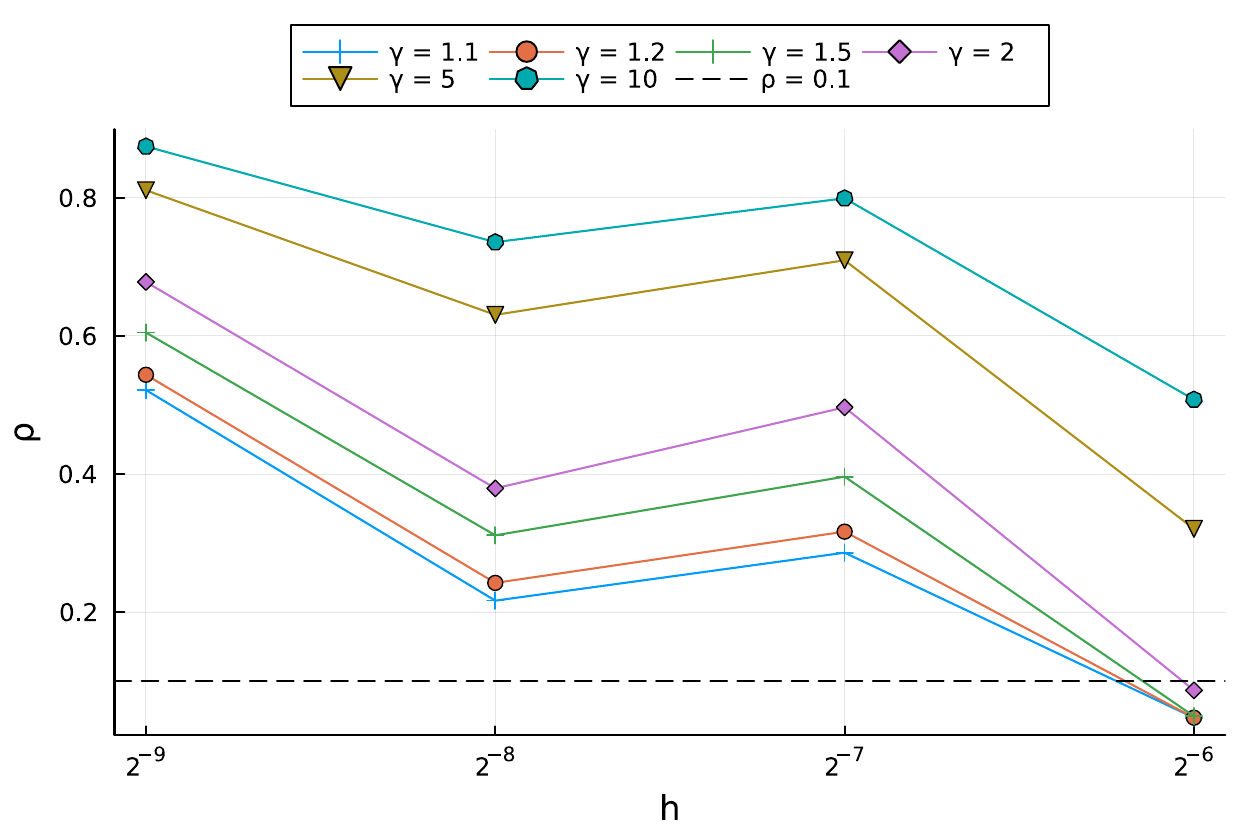}
   \caption{$\eta = 4$}
   \label{fig:2D_two_grid_circ_global_iter_4}
 \end{subfigure}%
 \caption{Plots of convergence factor of a \ac{tgcs} $\rho$ against mesh size $h$ on a circular domain using a global stabilization parameter, with $\eta$ additional iterations of smoother solely on cut cells.}
    \label{fig:2D_two_grid_circ_global}
\end{figure}

Therefore, in the next test we perform for each smoothing step $\eta$ additional iterations solely on the cut cells. The plots of the convergence factor against the mesh size is shown in \cref{fig:2D_two_grid_circ_local} using a local stabilization parameter and in \cref{fig:2D_two_grid_circ_global} using a global stabilization parameter for $\eta \in \{1,2,3,4\}$.

\par We observe that increasing the number of smoother iterations to the cut cells improves the convergence of the \ac{mg} solver. Specifically, with a choice of $\gamma = 2 $ and $\eta = 4 $ using a local stabilization parameter, optimal convergence is achieved, with $\rho \approx 0.1$. In contrast, when using a global stabilization parameter, the additional iterations enhance the convergence, but not as effectively as with the local parameter.
As a result, for the remaining tests, we will solely use local stabilization.
We point out that, with enough smoothing iterations, the \ac{mg} method can still reach optimal convergence even when using the global stabilization parameter. 

\par
\begin{figure}
    \centering
    \begin{subfigure}[t]{.5\textwidth}
   \centering
   \includegraphics[width=\textwidth]{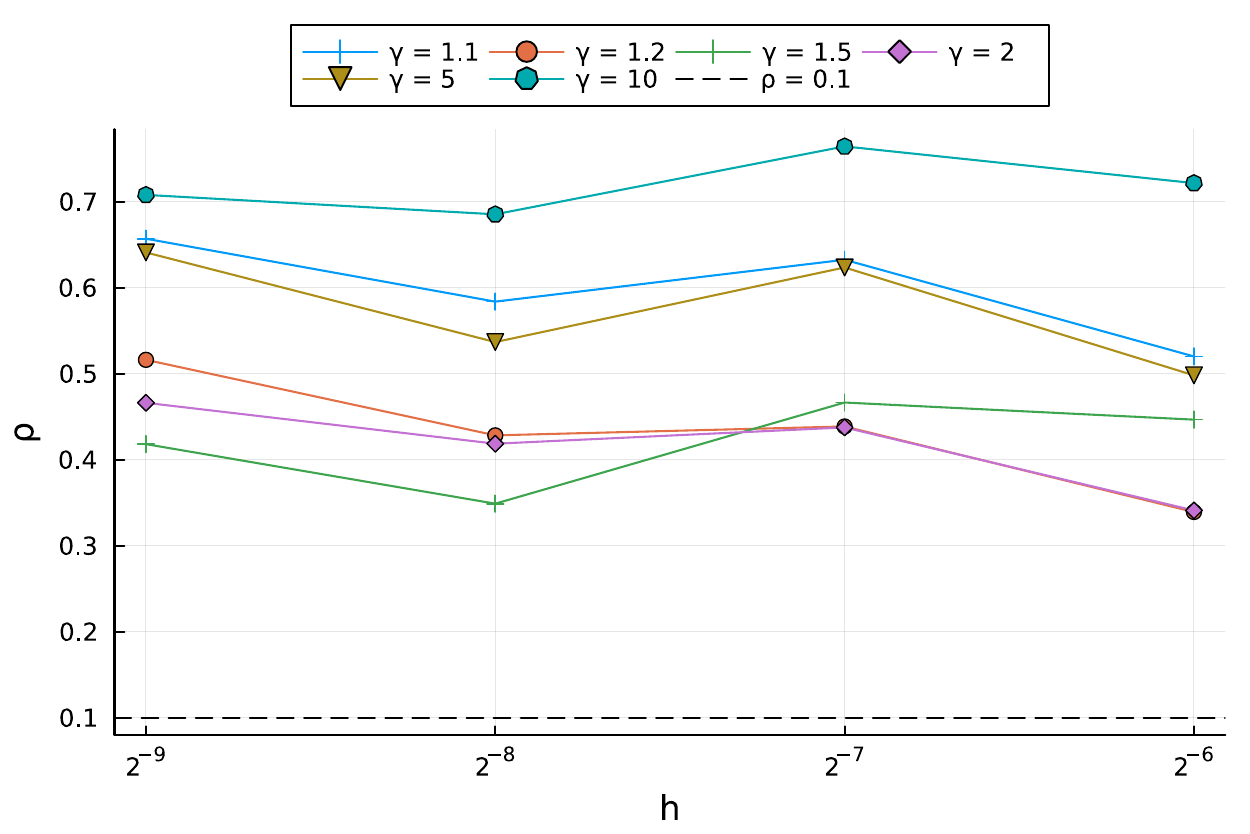}
   \caption{$\eta = 0$}
   \label{fig:2D_wcycle_iter_0}
 \end{subfigure}%
 \begin{subfigure}[t]{.5\textwidth}
   \centering
   \includegraphics[width=\textwidth]{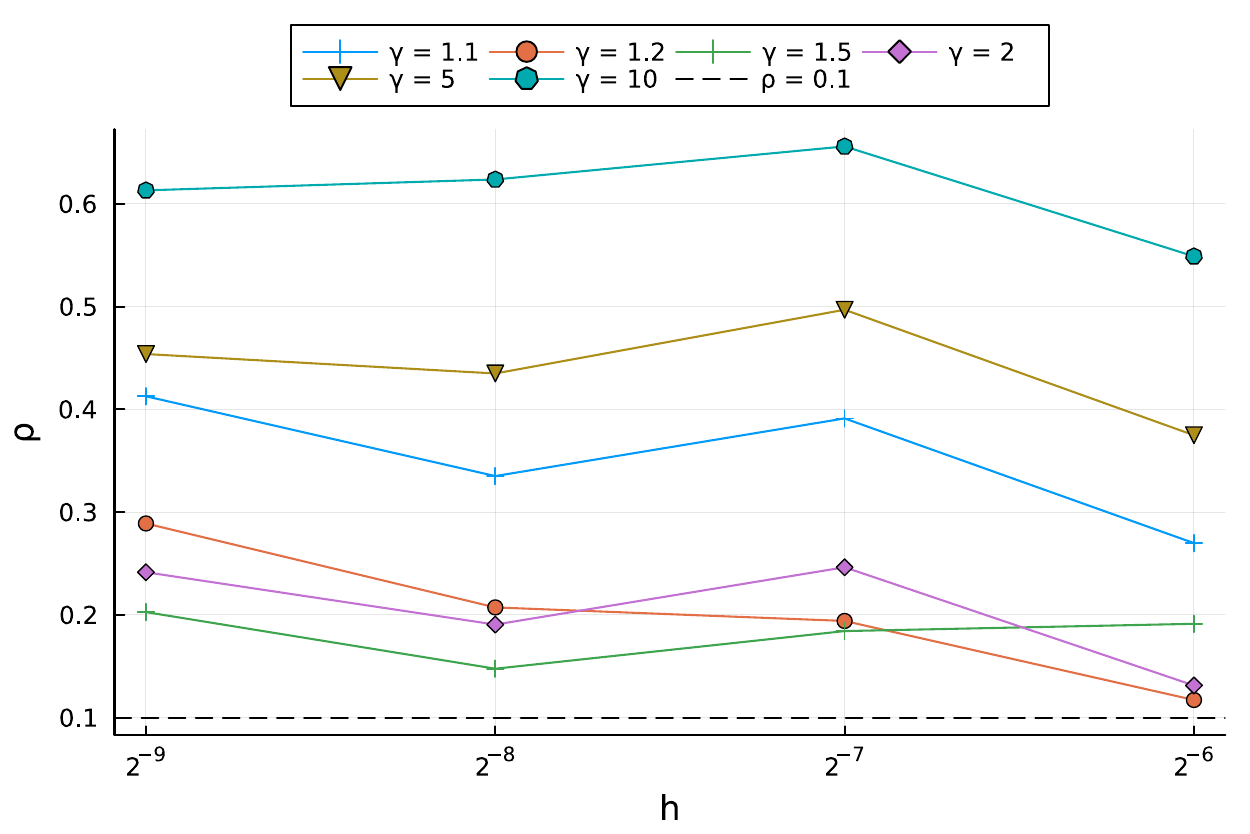}
   \caption{$\eta = 1$}
   \label{fig:2D_wcycle_iter_1}
 \end{subfigure}%
    \vskip\baselineskip
\begin{subfigure}[t]{.5\textwidth}
   \centering
   \includegraphics[width=\textwidth]{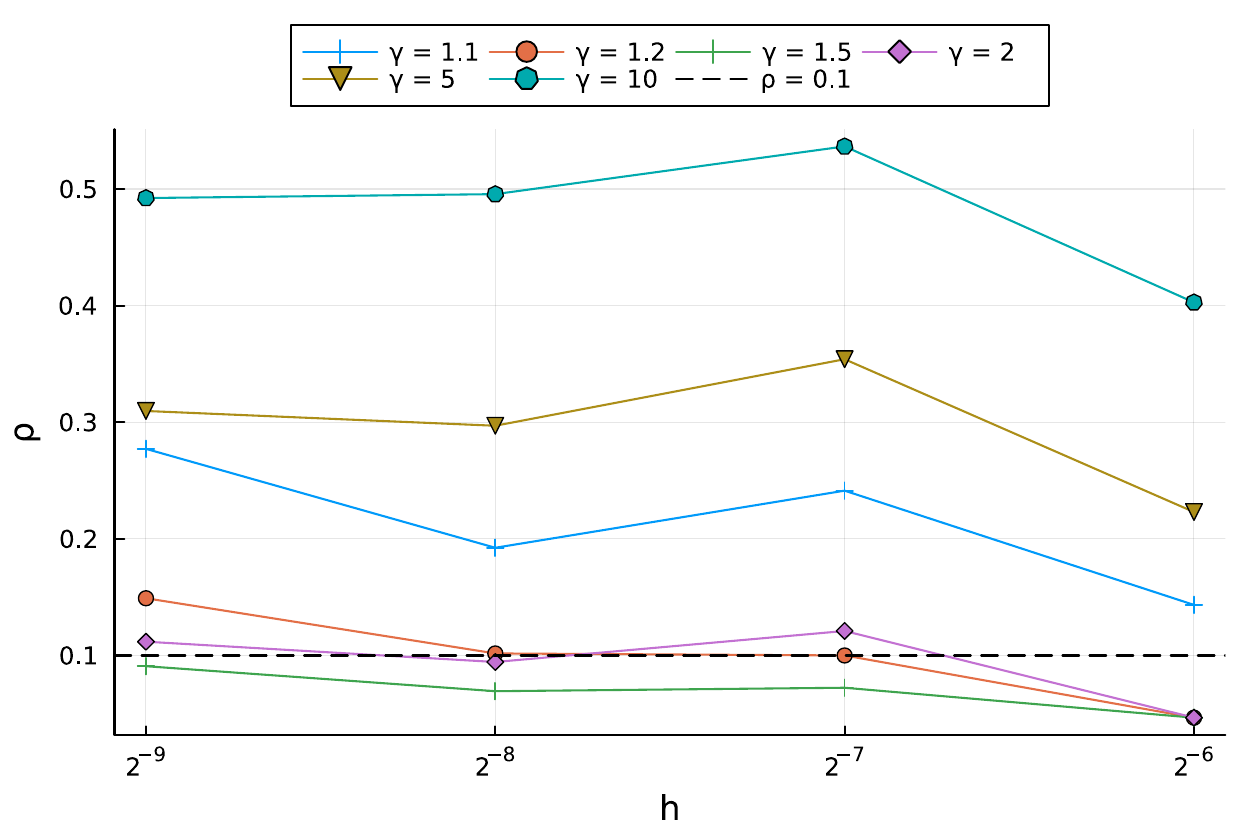}
   \caption{$\eta = 2$}
   \label{fig:2D_wcycle_iter_2}
 \end{subfigure}%
 \begin{subfigure}[t]{.5\textwidth}
   \centering
   \includegraphics[width=\textwidth]{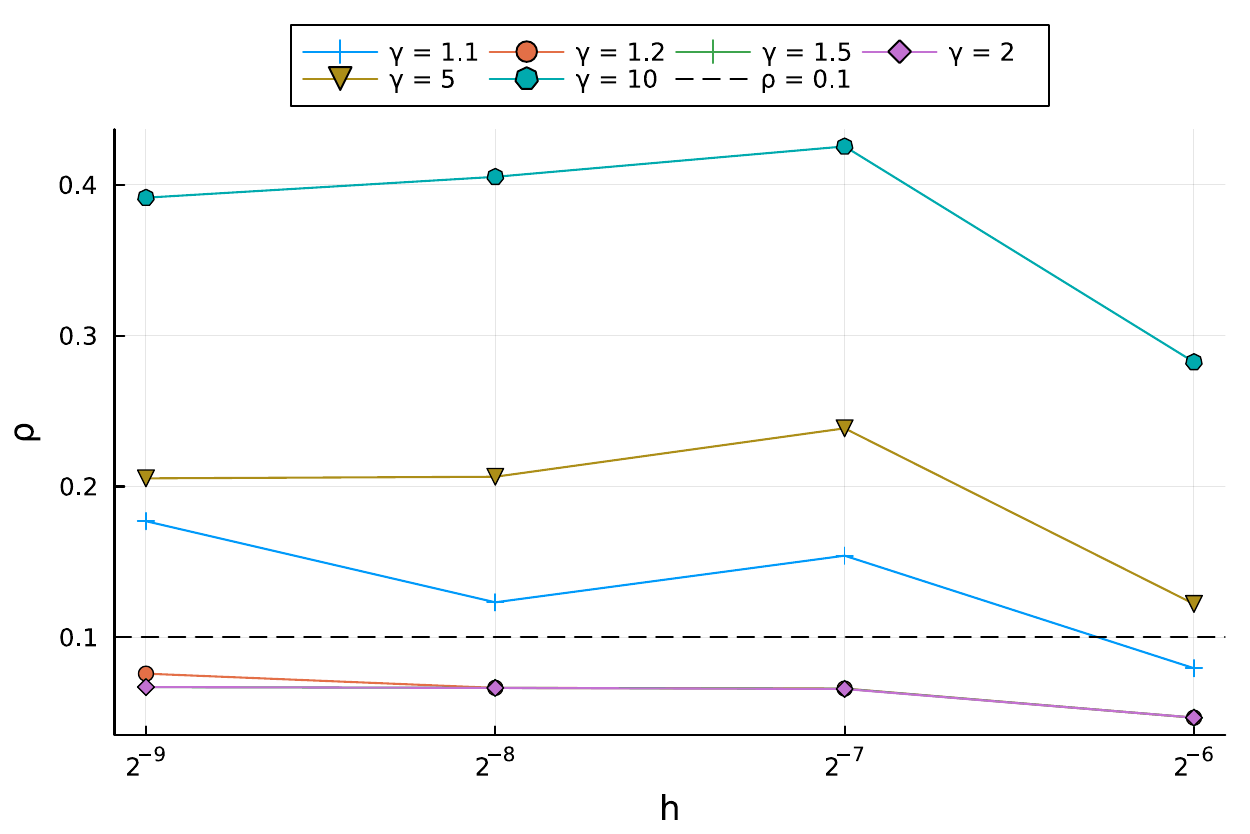}
   \caption{$\eta = 3$}
   \label{fig:2D_wcycle_iter_3}
 \end{subfigure}%
    \caption{Plot of convergence factor of a W-cycle against the mesh size on a circular domain for various values of $\gamma$ and $\eta$.}
    \label{fig:2D_wcycle}
\end{figure}

 In the final test, we assess the performance of a W-cycle. The local stabilization parameter is chosen as $\lambda(K) = \gamma \cdot C(K)$, where $C(K)$ is computed as outlined in \cref{sec:2D_optimal_lambda}, and $\gamma = \{1.1,1.2,1.5,2,5,10\}$. Additionally, we consider $\eta = \{0,1,2,3\}$.  The plots in \cref{fig:2D_wcycle} illustrate the convergence factor as a function of mesh size for different values of $\gamma$ and $\eta$. Notably, the convergence factor of the W-cycle matches that of the \ac{tgcs}, as expected~\cite{trottenberg2000multigrid}. The optimal convergence $\rho \approx 0.1$ is attained for $\gamma = 2$ and $\eta = 4$ for the \ac{tgcs}. Therefore, we adopt these values for all subsequent simulations.

\subsubsection*{Annulus-shaped domain} 
\par
\begin{figure}
    \centering
    \begin{subfigure}[t]{.5\textwidth}
   \centering
   \includegraphics[width=\textwidth]{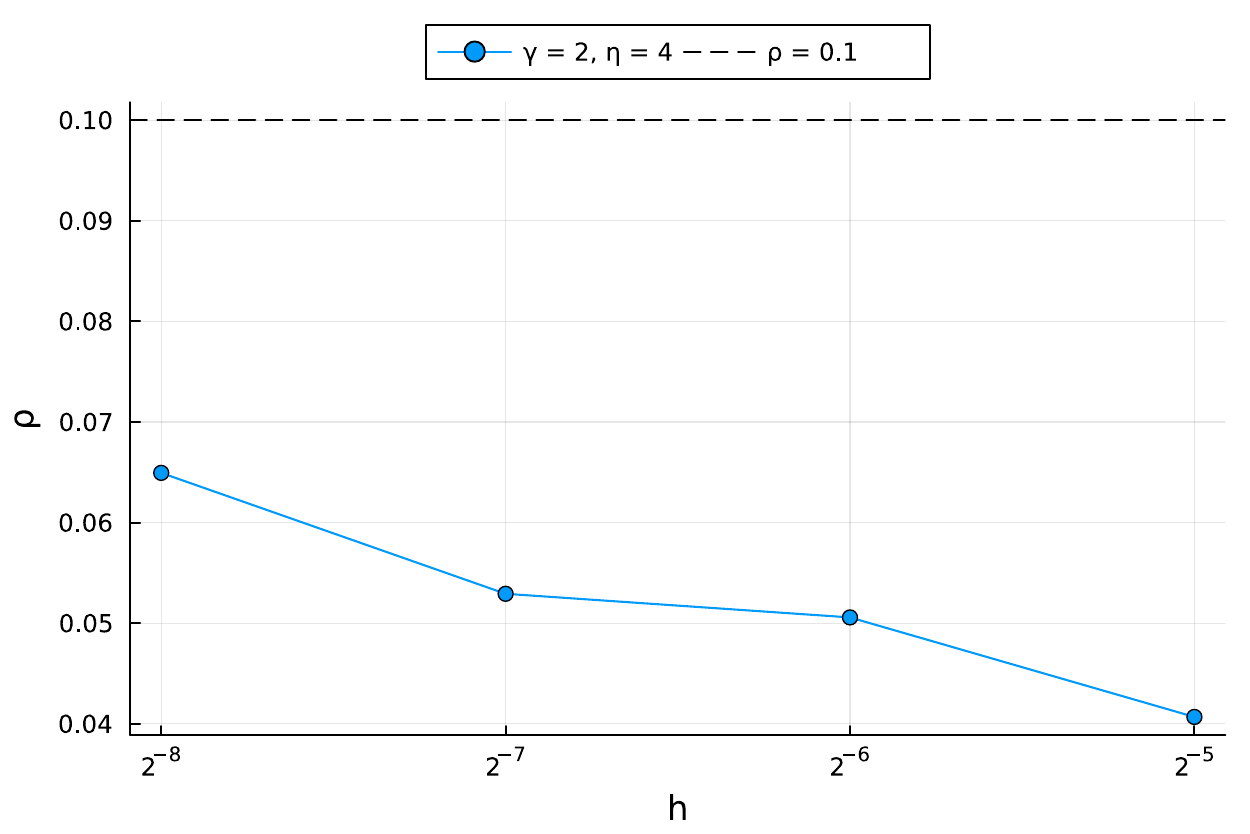}
   \caption{Annulus-shaped domain}
   \label{fig:2D_two_grid_annulus}
 \end{subfigure}%
 \begin{subfigure}[t]{.5\textwidth}
   \centering
   \includegraphics[width=\textwidth]{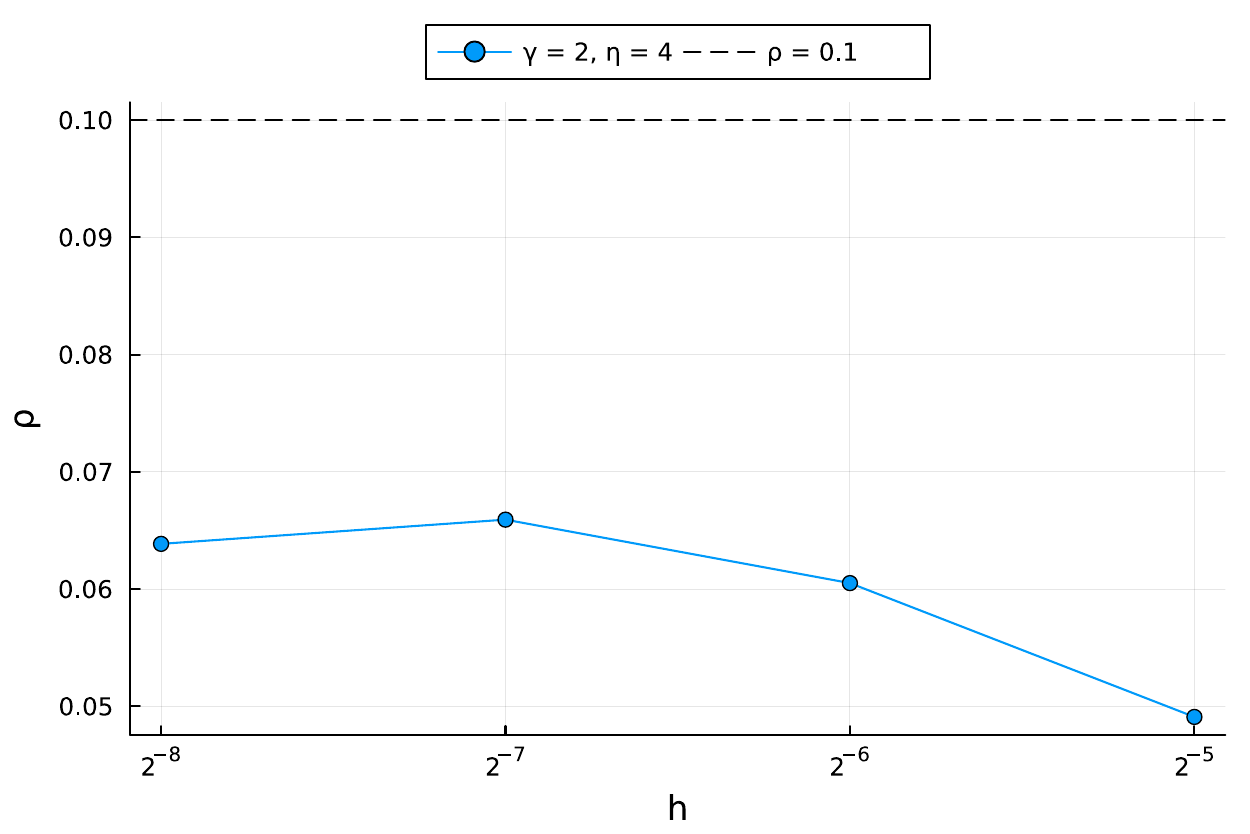}
   \caption{Flower-shaped domain}
   \label{fig:2D_two_grid_flower}
 \end{subfigure}%
\vskip\baselineskip
\begin{subfigure}{.5\textwidth}
   \centering
   \includegraphics[width=\linewidth]{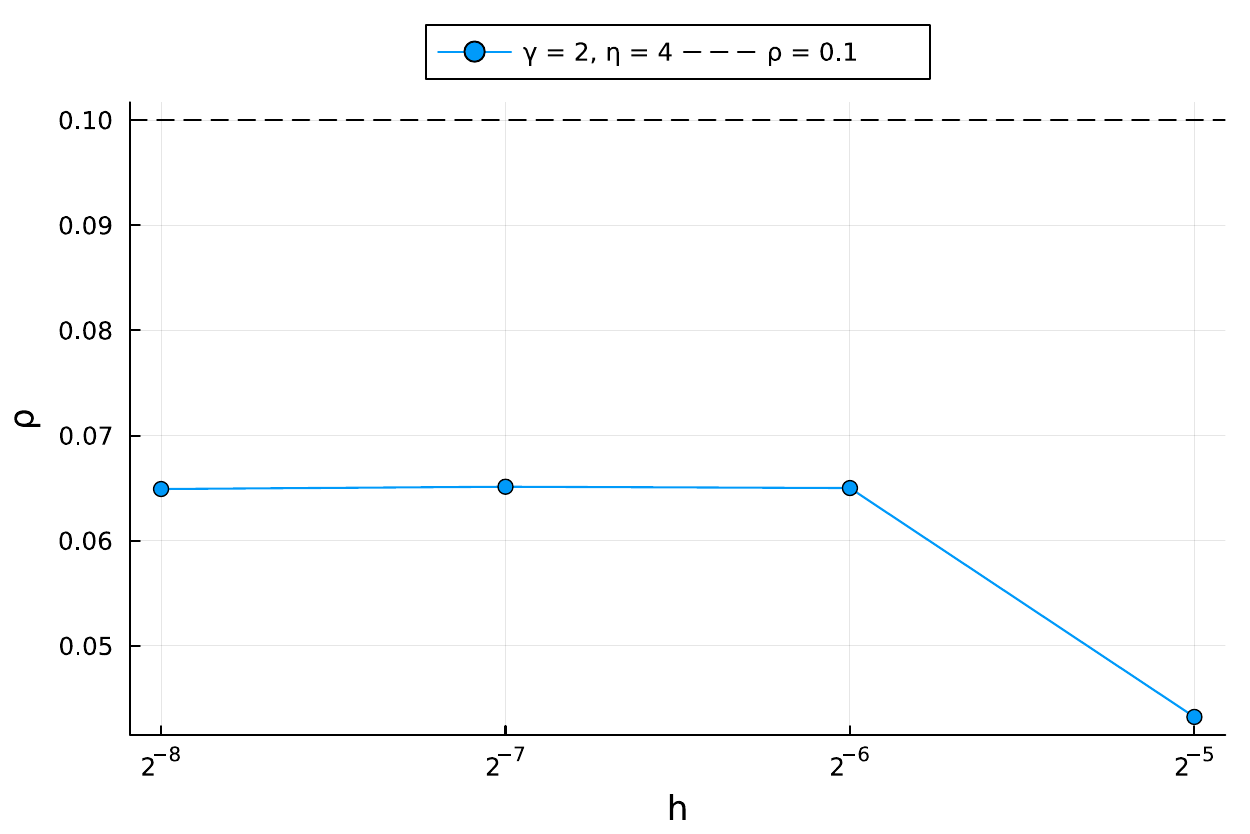}
   \caption{Leaf-shaped domain}
   \label{fig:2D_two_grid_leaf}
 \end{subfigure}%
 \begin{subfigure}{.5\textwidth}
   \centering
   \includegraphics[width=\linewidth]{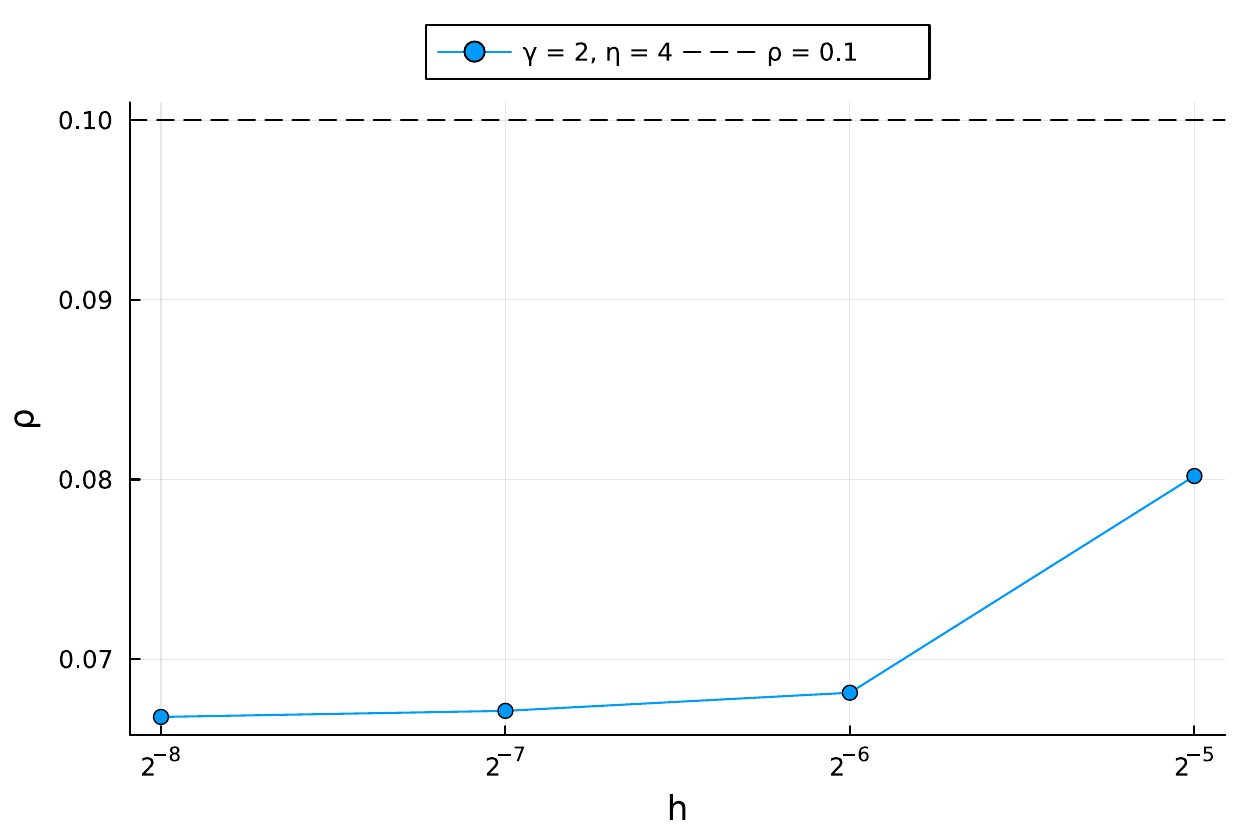}
   \caption{Hourglass-shaped domain}
   \label{fig:2D_two_grid_hourglass}
 \end{subfigure}%
    \caption{Plots of convergence factor of a \ac{tgcs} $\rho$ against mesh size $h$ for different geometries using local stabilization parameter, with $\gamma=2$ and $\eta=4$.}
    \label{fig:2D_two_grid_geos}
\end{figure}

 The annulus-shaped domain is defined as the region enclosed between two level sets (two circles): $\psi_1(x,y) = r_1^2 -(x-x_c)^2 - (y-y_c)^2 $  and $\psi_2(x,y) = (x-x_c)^2 + (y-y_c)^2 - r_2^2$, centered at $(x_c,y_c) = (0,0)$ with radii $r_1 = 0.5$ and $r_2 = 0.8$. The artificial domain is $\Omega_{art} = [-1,1]^2$ and the mesh size is $h =  2^{-m}, ~ m \in \{5,6,7,8\}$. Homogeneous Dirichlet and Neumann boundary conditions are imposed on $\Gamma_D = \{(x,y):\psi_1(x,y) = 0\}$ and $\Gamma_N = \{(x,y):\psi_2(x,y) = 0\} $, respectively. The cell-wise stabilization parameter is $\lambda(K) = 2 \cdot C(K)$. For each smoothing step, we perform $\eta = 4$ additional iterations exclusively for the cut cells. The plot of the convergence factor against different mesh size is shown in \cref{fig:2D_two_grid_annulus}.

\subsubsection*{Flower-shaped domain} 
\par The flower-shaped domain is defined by the level-set \cite{Astuto2024Comparison} 
\begin{align*}
    X &= x-0.03\sqrt{3}, \qquad Y = y- 0.04\sqrt{2}, \qquad R = \sqrt{X^2+Y^2} \\
    \psi &= R-0.52-\frac{(Y^5+5X^4Y-10X^2Y^3)}{5R^5}.
\end{align*}
The artificial domain is chosen as $\Omega_{art} = [-1,1]^2$, and homogeneous Dirichlet boundary conditions are imposed on $\Gamma_D = \{(x,y):\psi(x,y) = 0\}$. The local stabilization constant is set as $\lambda(K) = 2 \cdot C(K)$. Additionally, $\eta = 4$ extra smoothing iterations are applied exclusively to the boundary cells. The mesh size is $h =  2^{-m}, ~ m \in \{5,6,7,8\}$. The plot of the convergence factors against mesh size is illustrated in \cref{fig:2D_two_grid_flower}.

\subsubsection*{Leaf-shaped domain} 
\par The leaf-shaped domain is defined by the level set \cite{Astuto2024Comparison}
\begin{align*}
    &x_1 = -0.25 \cos(\pi/4), \quad R_1 = \sqrt{(x-x_1)^2+y^2},\\
    &x_2 = 0.25 \sin(\pi/4), \quad R_1 = \sqrt{(x-x_2)^2+y^2},\\ 
    &\psi_1 = R_1 - 0.7,  \quad \psi_2 = R_2-0.7, \quad \psi = \max\{\psi_1, \psi_2\}.
\end{align*}
The artificial domain is set as $\Omega_{art} = [-1,1]^2$ and the mesh size is $h =  2^{-m}, ~ m \in \{5,6,7,8\}$. Homogeneous Dirichlet and Neumann boundary conditions are imposed on $\Gamma_D = \partial \Omega \cap \{(x,y): x\ge 0\}$ and  $\Gamma_N = \partial \Omega \cap \{(x,y): x < 0\}$. The local stabilization constant is set as $\lambda(K) = 2 \cdot C(K)$. As before, $\eta = 4$ extra smoothing iterations are applied exclusively to the boundary cells. The plot of convergence factor against mesh size is illustrated in \cref{fig:2D_two_grid_leaf}.

\subsubsection*{Hourglass-shaped domain} 
\par The hourglass-shaped domain is defined using the level set \cite{Astuto2024Comparison}
\begin{align*}
    X = x-0.03\sqrt{3}, \quad Y = y- 0.04\sqrt{2}, \quad \psi = 256 Y^4 -16X^4 - 128Y^2 +36 X^2.
\end{align*}
The artificial domain is set as $\Omega_{art} = [-1,1]^2$, with a mesh size of $h =  2^{-m}, ~ m \in \{5,6,7,8\}$. Homogeneous Dirichlet boundary conditions are imposed on the whole boundary. The local stabilization constant is chosen as $\lambda(K) = 2 \cdot C(K)$. Furthermore, $\eta = 4$ extra smoothing iterations are applied exclusively to the boundary cells. The plot of convergence factor against mesh size is presented in \cref{fig:2D_two_grid_hourglass}. 

\par We observe that for all the different shaped domains, with an appropriate choice of the stabilization parameter, specifically, $\lambda = 2 \cdot C(K)$ and by performing four additional iterations of the smoother on the boundary, the \ac{tgcs} achieves optimal convergence, that is, $\rho \approx 0.1$. 

\section{Conclusions}\label{sec:conclusions}

\par In this paper, we have presented \ac{mg} methods for the ghost \ac{fe} approximation of elliptic \acp{pde} on 1D and 2D geometries. The ghost-\ac{fem} is second-order accurate and replacing a direct solver with a \ac{mg} method significantly improves computational efficiency.

\par We highlighted key limitations of ghost-point \ac{fdm} in the multigrid context, particularly the failure of classical smoothers such as Gauss-Seidel in certain parameter regimes. In contrast, ghost-\ac{fem} ensures convergence of standard relaxation schemes provided that the stabilization parameter $\lambda$ is chosen appropriately.
 
\par A further advantage of ghost-\ac{fem} lies in the treatment of residual transfers between grids. While the \ac{fdm}-based multigrid approach requires a splitting strategy to separate the contributions from internal equations and boundary conditions, we proved that such a strategy is unnecessary in the ghost-\ac{fem} framework. Specifically, the standard full-weighting restriction operator already captures the correct behavior by virtue of the underlying variational formulation, simplifying significantly the implementation and analysis of multigrid schemes in ghost-\ac{fem}.

\par Owing to the exact relationship between the shape functions on the fine and coarse grids, we showed that the coarse grid operator can be obtained either via the Galerkin condition or by assembling the operator directly on the coarse mesh using coarse-level shape functions.

\par The choice of the stabilization parameter is crucial for ensuring optimal \ac{mg} convergence. It must be sufficiently large to guarantee the well-posedness of the ghost \ac{fe} formulation, yet an excessively large value can degrade efficiency of the \ac{mg} solver, although the second-order accuracy is maintained. The stabilization parameter can be determined either globally or locally by solving a generalized eigenvalue problem. In one dimension, the local and global stabilization parameters are identical, as there is exactly one cut cell. However, in two dimensions, the local approach is both computationally efficient and optimal due to the varying configurations of cut cells.

\par We discussed the geometric characteristics of cut cells that arise from intersecting a structured Cartesian grid with an implicitly defined domain in two dimensions. The shape of these cut cells (typically triangles, quadrilaterals, or pentagons) is determined by evaluating the level set function at cell vertices. For each shape, we computed a lower bound for the local stabilization parameter. Additionally, performing extra iterations of the smoother exclusively on the cut cells further improves \ac{mg} convergence, without compromising the computational cost. Numerical tests confirm that the proposed \ac{mg} approach achieves optimal convergence across various geometries. 

\section*{Acknowledgements}
This work has been supported by Italian Ministerial grant PRIN 2022 “Efficient numerical schemes and optimal control methods for time-dependent partial differential equations”, No. 2022N9BM3N - Finanziato dall’Unione europea - Next Generation EU – CUP: E53D23005830006.

The work of A.C.~has been supported by the Spoke 1 Future HPC \& Big Data of the Italian Research Center on High-Performance Computing, Big Data and Quantum Computing (ICSC) funded by MUR Missione 4 Componente 2 Investimento 1.4: Potenziamento strutture di ricerca e creazione di “campioni nazionali di R \&S (M4C2-19)” - Next Generation EU (NGEU).

The work of A.C.~has been supported by Italian Ministerial grant PRIN 2022 PNRR “FIN4GEO: Forward and Inverse Numerical Modeling of hydrothermal systems in volcanic regions with application to geothermal energy exploitation”, No. P2022BNB97 - Finanziato dall’Unione europea - Next Generation EU – CUP: E53D23017960001.

H.D.~and A.C.~are members of the Gruppo Nazionale Calcolo Scientifico-Istituto Nazionale di Alta Matematica (GNCS-INdAM).

We acknowledge the CINECA award under the ISCRA initiative (ISCRA C project ``Efficient Higher Order Finite Element Schemes for Fluid Dynamics using Multigrid Methods"), for the availability of high-performance computing resources and support.

\printbibliography
\end{document}